\documentclass[12pt]{amsart}
\usepackage{amssymb}
\usepackage{epsf}
\usepackage{multirow}

\newtheorem{thm}{Theorem}[section]
\newtheorem{prop}[thm]{Proposition}

\newtheorem{cor}[thm]{Corollary}
\newtheorem{conjecture}[thm]{Conjecture}
\newenvironment{conj}{\begin{conjecture}\rm}{\end{conjecture}}
\newtheorem{question}[thm]{Question}
\newenvironment{ques}{\begin{question}\rm}{\end{question}}
\newtheorem{remark}[thm]{Remark}
\newenvironment{rem}{\begin{remark}\rm}{\end{remark}}

\numberwithin{figure}{section}
\numberwithin{equation}{section}
\numberwithin{table}{section}

\newcounter{FNC}[page]
\def\frankfootnote#1{{\addtocounter{FNC}{2}$^\fnsymbol{FNC}$%
     \let\thefootnote\relax\footnotetext{$^\fnsymbol{FNC}$#1}}}

\font\CYR=wncyr10 at 12pt   
\newcommand{\Tc}{\mbox{\CYR Ch}}
\newcommand{\DOT}{{\setlength{\unitlength}{1pt}\begin{picture}(3.5,2)(1,1)
\put(2.3,2){\circle*{2}}\end{picture}}}
\newcommand{\Fdot}{{F\!_{\DOT}}}
\newcommand{\hFdot}{{\hat{F}\!_{\DOT}}}
\newcommand{\RP}{{\mathbb R}{\mathbb P}}
\newcommand{\QED}{\quad%
\raisebox{-6pt}{\epsfysize=15pt \epsfbox{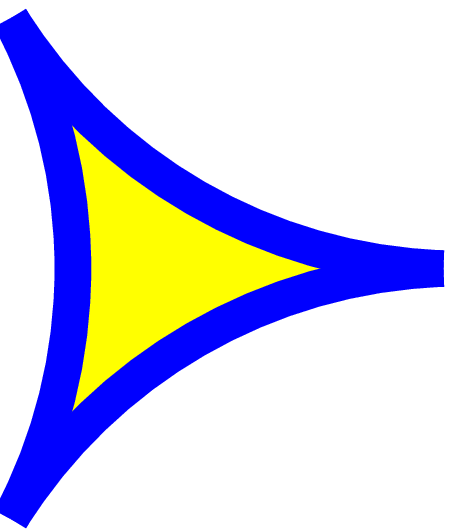}}\vspace{10pt}}
\newcommand{\Qed}{\quad%
\raisebox{-6pt}{\epsfysize=15pt \epsfbox{figures/qed.eps}}}
\newcommand{\Span}[1]{{\langle #1 \rangle}}
\newcommand{\sgn}{\mbox{\rm sgn}}
\newcommand{\calN}{{\mathcal N}}

\newcommand{\RE}{{\mathrm{Re}}}

{\catcode`\@=11
\gdef\n@te#1#2{\leavevmode\vadjust{%
 {\setbox\z@\hbox to\z@{\strut#1}%
  \setbox\z@\hbox{\raise\dp\strutbox\box\z@}\ht\z@=\z@\dp\z@=\z@%
  #2\box\z@}}}
\gdef\leftnote#1{\n@te{\hss#1\quad}{}}
\gdef\rightnote#1{\n@te{\quad\kern-\leftskip#1\hss}{\moveright\hsize}}
\gdef\?{\FN@\qumark}
\gdef\qumark{\ifx\next"\DN@"##1"{\leftnote{\rm##1}}\else
 \DN@{\leftnote{\rm??}}\fi{\rm??}\next@}}

\begin{document}

\title[Maximally inflected real rational curves]{Maximally inflected real
rational curves}

\author{Viatcheslav Kharlamov}
\address{Institut de Recherche Math\'ematique Avanc\'ee\\
        UMR 7501 de l'ULP et du CNRS\\
        7 rue Ren\'e-Descartes\\
        67084 Strasbourg Cedex\\
        France}
\email{kharlam@math.u-strasbg.fr}

\author{Frank Sottile}
\address{Department of Mathematics\\
        University of Massachusetts\\
        Lederle Graduate Research Tower\\
        Amherst, MA 01003\\
        USA}
\email{sottile@math.umass.edu}
\urladdr{http://www.math.umass.edu/\~{}sottile}

\date{1 July 2003}
\thanks{Research of second author supported in part by IRMA
Strasbourg (June 1999), IRMAR Rennes, and NSF grant DMS-0070494}

\begin{abstract}
We introduce and begin the topological study of real rational plane curves,
all of whose inflection points are real.
The existence of such curves is a corollary of results in the real Schubert
calculus, and their study has
consequences for the important Shapiro and Shapiro
conjecture in the real Schubert calculus.
We establish restrictions on the number of real nodes of such curves
and construct curves realizing the extreme numbers of real nodes.
These constructions imply the existence of real solutions to
some problems in the Schubert calculus.
We conclude with a discussion of maximally inflected curves of low degree.
\end{abstract}

\maketitle

\section*{Introduction}
In 1876, Harnack~\cite{Ha1876} established the bound of $g+1$ for the number
of ovals of a smooth curve of genus $g$ in $\RP^2$ and showed
this bound is sharp by constructing curves with
this number  of ovals.
Since then such curves with maximally
many ovals, or M-curves, have been primary objects of interest in
the topological study of real algebraic plane curves (part of
the 16th problem of Hilbert)~\cite{Wilson,Viro}.
We introduce and begin
the study of maximally inflected curves, which may be considered to be analogs
of $M$-curves among rational curves, in the sense that, like the classical
maximal condition, maximal inflection implies non-trivial restrictions on the
topology.
The case of a maximally inflected curve that we are most interested in is a
rational plane
curve of degree $d$, all of whose $3(d{-}2)$ inflection points are real.
More generally, consider a parameterization of a rational curve of degree $d$
in ${\mathbb P}^r$ with $d>r$;  a map from ${\mathbb P}^1$ to
${\mathbb P}^r$ of degree $d$.
Such a map is said to be {\it ramified}
at a point $s\in{\mathbb P}^1$ if its first $r$
derivatives do not span ${\mathbb P}^r$.
Over $\mathbb C$, there always will be
$(r{+}1)(d{-}r)$ such points of ramification, counted with
multiplicity.
A {\it maximally inflected curve}
is, by definition, a parameterized real rational
curve, all of whose ramification points are real.

We first address the question of existence
of maximally inflected curves.
A conjecture of B.~Shapiro and M.~Shapiro in the real Schubert
calculus~\cite{So00b} would imply the existence of maximally
inflected curves, for every possible placement of the ramification points.
A.~Eremenko and A.~Gabrielov~\cite{EG02} proved the conjecture
of Shapiro and Shapiro in the special case when $r=1$.
In that case, a maximally inflected curve is a real rational function, all of
whose critical points are real.

When $r>1$ and for special types of ramification (including flexes and
cusps of plane curves) and when the ramification points are
clustered very near to one another, there
do exist maximally inflected curves, by a result in the real Schubert
calculus~\cite{So99a}.
When the degree
$d$ is even and the curve has only simple flexes, but arbitrarily placed,
the existence of such curves follows from the results of Eremenko and
Gabrielov on the degree of the Wronski map~\cite{EG01b}.

Suppose $r=2$, the case of plane curves.
Consider a maximally inflected plane curve of degree $d$ with the maximal
number $3(d{-}2)$ of real flexes.
By the genus formula, the  curve has at most $\binom{d-1}{2}$
ordinary double points.
We deduce from the Klein~\cite{Klein} and Pl\"ucker~\cite{Plucker} formulas that
the number of real nodes of such a curve is however at most
$\binom{d{-}2}{2}$.
Then we show that this bound is sharp.
We use E.~Shustin's theorem on
combinatorial patchworking for singular
curves~\cite{Sh85} to construct maximally inflected curves
that realize the lower bound of 0 real nodes, and
another theorem of Shustin concerning
deformations of singular curves~\cite{Sh99} to construct maximally
inflected curves with the upper bound of $\binom{d{-}2}{2}$ real nodes.
We also classify such curves in degree 4, 
and discuss some aspects of the classification of quintics.

The connection between the Schubert calculus and rational curves in
projective space (linear series on ${\mathbb P}^1$) originated in
work of G.~Castelnuovo~\cite{Ca1889} on $g$-nodal rational curves.
This led to the use of Schubert calculus in Brill-Noether
theory (see Chapter~5 of~\cite{HaMo} for an elaboration).
In a closely related vein, the Schubert calculus features in the local
study of flattenings of curves in singularity theory~\cite{Shch82,Ka91,AVGL89}.
In turn, the theory of limit linear series of D.~Eisenbud and
J.~Harris~\cite{EH83,EH87}
provides essential tools to show reality of the
special Schubert calculus~\cite{So99a}, which gives the existence of
many types of maximally inflected curves.
The constructions we give using patchworking and gluing show the existence of
real solutions to some problems in the Schubert calculus.
In particular, they show the existence of maximally
inflected plane curves with given type of ramification when the degree 
$d$ is odd, extending the results of Eremenko-Gabrielov on the degree of the
Wronski map~\cite{EG01b}.

This paper is organized as follows.
In Section 1 we describe the connection between the Schubert calculus and
maps from ${\mathbb P}^1$ with specified points of ramification.
In Section 2, we show the existence of some
special types of maximally inflected curves and discuss the
conjecture of Shapiro and Shapiro in this context.  We restrict our
attention to maximally inflected plane curves in Section 3, establishing
bounds on the number of solitary points and real nodes.
In Section 4, we construct maximally inflected curves with extreme numbers of
real nodes and in Section 5, we discuss the classification of
quartics.
In Section 6, we consider some aspects of the classification of
quintics.
\medskip

We thank D. Pecker who suggested using duals of curves in the proof of
Corollary~\ref{cor:bounds} and E. Shustin who explained to us his work on
deformations of singular
curves and suggested the construction of Section~\ref{sec:patchworking} at the
Oberwolfach workshop ``New perspectives in the topology of real algebraic
varieties'' in September 2000.
We also thank the referee for his useful comments and for pointing out
references concerning the local study of flattenings of curves.

\section{Singularities of parameterized curves}

\subsection{Notations and conventions}\label{S:notation}
Given non-zero vectors $v_1,v_2,\dotsc,v_n$ in $\mathbb{A}^{r+1}$, we denote
their linear
span in the projective space $\mathbb{P}^r$ by
$\Span{v_1,v_2,\dotsc,v_n}$.
A subvariety $X$ of $\mathbb{P}^r$ is \emph{nondegenerate} if its linear
span is $\mathbb{P}^r$, equivalently, if it does not lie in a
hyperplane.
A rational map between varieties is denoted by a broken arrow
$X \relbar\to Y$.

Let $d>r$.
The center of a (surjective) linear projection
${\mathbb P}^d \relbar\twoheadrightarrow{\mathbb P}^r$
is a $(d{-}r{-}1)$-dimensional linear subspace
of ${\mathbb P}^d$.
This center determines the projection up to projective
transformations of $\mathbb{P}^r$.
Consequently, we identify
the Grassmannian $\textit{Grass}_{d-r-1}{\mathbb P}^d$ of
$(d{-}r{-}1)$-dimensional linear subspaces of ${\mathbb P}^d$
with the space of linear projections
${\mathbb P}^d \relbar\twoheadrightarrow{\mathbb P}^r$
(modulo  projective transformations of $\mathbb{P}^r$).

A flag $\Fdot$ in ${\mathbb P}^d$ is a sequence
 \[
   \Fdot\ \colon\ F_0\ \subset\ F_1\ \subset\ \dotsb\
           \subset\ F_d\ =\ {\mathbb P}^d
 \]
of linear subspaces with $\dim F_i=i$.
Given a flag $\Fdot$ in ${\mathbb P}^d$ and a sequence
$\alpha\colon 0\leq\alpha_0<\alpha_1<\dotsb<\alpha_r\leq d$
of integers, the \emph{Schubert cell}
$ X^\circ_\alpha\Fdot\subset\textit{Grass}_{d-r-1}{\mathbb P}^d$
is the set of all linear projections
$\pi\colon{\mathbb P}^d \relbar\twoheadrightarrow{\mathbb P}^r$
such that
 \begin{equation}\label{E:Schubert}
   \dim \pi(F_{\alpha_i})\ =\ i,\quad\text{ and, \,if $i\neq 0$, }
    \quad \dim \pi(F_{\alpha_i-1})\ =\ i-1\,.
 \end{equation}
The Schubert variety $ X_\alpha\Fdot$ is the closure of the Schubert
cell; it is obtained by replacing the equalities in~\eqref{E:Schubert}
by inequalities ($\leq$).

Replacing $\mathbb{P}^d$ by the dual projective space $\hat{\mathbb{P}}^d$
gives an isomorphism
$\textit{Grass}_{d-r-1}{\mathbb P}^d\simeq\textit{Grass}_r\hat{\mathbb{P}}^d$;
a $(d{-}r{-}1)$-plane $\Lambda$ in $\mathbb{P}^d$ corresponds to the
$r$-plane in $\hat{\mathbb{P}}^d$ consisting of hyperplanes of $\mathbb{P}^d$
that contain $\Lambda$.
Under this isomorphism, the Schubert variety $ X_\alpha\Fdot$ is mapped
isomorphically to the Schubert variety $ X_{\hat{\alpha}}\hFdot$,
where $\hFdot$ is the flag dual to $\Fdot$
($\hat{F}_i$ consists of the hyperplanes containing $F_{d-i}$) and
$\hat{\alpha}$ is the sequence
$0\leq\beta_0<\beta_1<\dotsb<\beta_{d-r-1}\leq d$ such that
 \begin{equation}\label{E:hat-def}
  \{0,1,\dotsc,d\}\ =\
   \{\alpha_0,\alpha_1,\dotsc,\alpha_r\}\ \cup
   \{d-\beta_0,d-\beta_1,\dotsc,d-\beta_{d-r-1}\}\,.
 \end{equation}

\subsection{Parameterized rational curves}\label{S:1.2}

Given $r+1$ coprime linearly independent binary homogeneous forms of the
same degree $d$, we have a morphism
$\varphi\colon\mathbb{P}^1\to\mathbb{P}^r$
whose image is nondegenerate and of degree $d$ in that
$\varphi_*([\mathbb{P}^1])=d[\ell]$, where $[\ell]$ is the
standard generator in $H_2(\mathbb{P}^r)$.
Reciprocally, every morphism
from $\mathbb{P}^1$ to $\mathbb{P}^r$ whose image is
nondegenerate
and of degree $d$
is given by $r+1$ coprime linearly
independent binary homogeneous forms of the same degree $d$.
We consider two such
maps to be equivalent when they differ by a projective
transformation of ${\mathbb P}^r$.
In what follows, we will call an equivalence class of such maps a
{\it rational curve of degree $d$ in} ${\mathbb P}^r$.

A rational curve $\varphi\colon\mathbb{P}^1\to\mathbb{P}^r$ of degree $d$
is said to be {\it ramified} at $s\in{\mathbb P}^1$ if the
derivatives
$\varphi(s),\varphi'(s),\ldots,\varphi^{(r)}(s)\in\mathbb{A}^{r+1}$
do not span ${\mathbb P}^r$.
This occurs at the roots of the Wronskian
 \[
   \det\left[\begin{array}{ccc}
           \varphi_1(s)&\cdots&\varphi_{r+1}(s)\\
           \varphi'_1(s)&\cdots&\varphi'_{r+1}(s)\\
           \vdots&\ddots&\vdots\\
          \varphi^{(r)}_1(s)&\cdots&\varphi^{(r)}_{r+1}(s)
         \end{array}\right]\ ,
 \]
of $\varphi$, a form of degree $(r{+}1)(d{-}r)$.
(Here, $\varphi_1,\dotsc,\varphi_{r+1}$ are the components of $\varphi$.)
The Wronskian is defined by the curve $\varphi$ only up
to multiplication by a
scalar; in particular its set of roots and their multiplicities are
well-defined in
the
equivalence class
of $\varphi$.
An equivalence class of such maps containing a real map $\varphi$ is a
{\it maximally inflected curve} when the
Wronskian has only real roots.

An equivalent formulation is provided by linear series on ${\mathbb P}^1$.
A rational curve $\varphi:{\mathbb P}^1\rightarrow{\mathbb P}^r$
of degree $d$ whose image is nondegenerate may be factored in the following way
 \begin{equation}\label{E:phi-def}
  {\mathbb P}^1\ \stackrel{\gamma}{\longrightarrow}\
  {\mathbb P}^d\ \stackrel{\pi}{\relbar\rightarrow}\
  {\mathbb P}^r
 \end{equation}
where $\gamma$ is the {\it rational normal curve},
which is given by
$$
  [s,t]\ \longmapsto\ [s^d,\, s^{d-1}t,\, s^{d-2}t^2,\, \ldots,\, st^{d-1},\,
  t^d]\,,
$$
and $\pi$ is a projection whose center $\Lambda$ is a ($d{-}r{-}1$)-plane which
does not meet the curve $\gamma$.
The hyperplanes of ${\mathbb P}^d$ which contain $\Lambda$
constitute
a base point
free linear series of dimension $r$ in
${\mathbb P}H^0({\mathbb P}^1,{\mathcal O}(d))={\mathbb P}^d$.
(Base point freeness is equivalent to our requirement
that the original forms
$\varphi_i$ do not have a common factor.)
The center of projection $\Lambda$ and hence the linear series depends only
upon the original map.
In this way, the space of rational curves
$\varphi:{\mathbb P}^1\rightarrow{\mathbb P}^r$ of degree $d$ is identified
with the (open) subset of the Grassmannian
$\textit{Grass}_{d-r-1}\mathbb{P}^d$ of $(d{-}r{-}1)$-planes in $\mathbb{P}^d$
which do not meet the image of the rational normal curve $\gamma$.

Let us define the {\it ramification sequence} of $\varphi$
at a point $s$  in ${\mathbb P}^1$
to be the vector $\alpha=\alpha(s)$ whose $i$th component is the smallest
number $\alpha_i$  such that the linear span
of $\varphi(s),\varphi'(s),\dotsc$, $\varphi^{(\alpha_i)}(s)$
in $\mathbb{P}^r$
has dimension  $i$.
Equivalently, ${\mathbb P}^r$ has coordinates such that
$\varphi=[\varphi_0,\varphi_1,\ldots,\varphi_r]$ with $\varphi_i$ vanishing to
order $\alpha_i$ at $s$.
We will also call this sequence $\alpha$
the {\it ramification} of $\varphi$ at $s$.
Note that $\alpha_0=0$ and $d\geq\alpha_i\ge i$ for $i=0,\dotsc,r$, and the
sequence $\alpha$ is increasing.
Call such an integer vector $\alpha$ a {\it ramification sequence} for degree
$d$ curves in $\mathbb{P}^r$ if $\alpha_r>r$.
If $s$ is a point of $\mathbb{P}^1$ with ramification sequence
$\alpha$, then the Wronskian vanishes to order
$|\alpha|:=(\alpha_1{-}1)+(\alpha_2{-}2)+\cdots+(\alpha_r{-}r)$ at $s$.
When this is zero, so that $\alpha_r=r$, we say that $\varphi$ is
{\it unramified} at $s$.
Since the Wronskian has degree $(r+1)(d-r)$,
the following fact holds.

\begin{prop}\label{prop:ram}
 Let $\varphi\colon\mathbb{P}^1\to\mathbb{P}^r$ be a rational curve of degree
 $d$.
 Then
$$
  \sum_{s\in \mathbb{P}^1} |\alpha(s)|\ =\ (r+1)(d-r)\,. \Qed
$$
\end{prop}

A collection $\alpha^1,\ldots,\alpha^n$  of ramification sequences
for degree $d$ curves in $\mathbb{P}^r$ with
$|\alpha^1|+\cdots+|\alpha^n|=(r{+}1)(d{-}r)$
will be called {\it ramification data}
for degree $d$ rational curves in ${\mathbb P}^r$.

A ramification sequence
is expressed in terms of Schubert cycles in
$\mbox{\it Grass}_{d-r-1}{\mathbb P}^d$ in the following manner.
For each $s\in{\mathbb P}^1$, let $\Fdot(s)$ be the flag of subspaces of
${\mathbb P}^d$ which
osculate the rational normal curve $\gamma$ at the point $\gamma(s)$
and $\pi$ the linear projection defining $\varphi$ as
in~\eqref{E:phi-def}.
Then the image $\pi(F_i(s))$ of the $i$-plane $F_i(s)$ in $\Fdot(s)$ is the
linear
span of $\varphi(s),\varphi'(s),\ldots,\varphi^{(i)}(s)$
in ${\mathbb P}^r$.

Comparing the definitions of ramification sequence and of Schubert cell,
we deduce
that the curve $\varphi$ has ramification
$\alpha$ at $s\in{\mathbb P}^1$
when the center of projection lies in the
Schubert cell $ X^\circ_\alpha\Fdot(s)$ of the Grassmannian.
The closure of this cell is the Schubert cycle $ X_\alpha\Fdot(s)$, which
has codimension $|\alpha|$ in the Grassmannian.
Thus, given ramification data
$\alpha^1,\ldots,\alpha^n$ and distinct points
$s_1,\ldots,s_n\in{\mathbb P}^1$,
a center $\Lambda$ provides
ramification $\alpha^i$ at point $s_i$ for each $i$
if and only if it belongs to the intersection
$$
  \bigcap_{i=1}^n  X^\circ_{\alpha^i}\Fdot(s_i)\,
$$
On the other hand,
$$
  \bigcap_{i=1}^n  X^\circ_{\alpha^i}\Fdot(s_i)\ =\
  \bigcap_{i=1}^n  X_{\alpha^i}\Fdot(s_i)\,
$$
To see this, suppose that a point in
$\bigcap_{i=1}^n  X_{\alpha^i}\Fdot(s_i)$
did not lie in a
Schubert cell $ X^\circ_{\alpha^i}\Fdot(s_i)$.
Then the corresponding rational curve would have ramification exceeding
$\alpha^i$ at the point $s_i$, which would contradict Proposition~\ref{prop:ram}.

The intersection of the above Schubert
cycles has dimension at least
$(r{+}1)(d{-}r)-\sum_i|\alpha^i|=0$.
If it had positive dimension, then it
would have a non-empty intersection
with any hypersurface Schubert variety $ X_\beta\Fdot(s)$,
where $|\beta|=1$
and $s$ is not among $\{s_1,\dotsc,s_n\}$.
This gives a rational curve violating Proposition~\ref{prop:ram}.
Thus the intersection is zero-dimensional.
Let $N(\alpha^1,\ldots,\alpha^n)$ be its degree, which may be
computed using the classical Schubert calculus of enumerative
geometry~\cite{GH78}.
Thus we deduce the following Corollary of Proposition~\ref{prop:ram}.

\begin{cor}\label{C:number_curves}
 Given $d>r$, ramification data $\alpha^1,\ldots,\alpha^n$ for degree $d$
 rational curves in ${\mathbb P}^r$, and distinct points
 $s_1,\ldots,s_n\in{\mathbb P}^1$,
 the number of nondegenerate rational curves in ${\mathbb P}^r$ of degree
 $d$ with ramification $\alpha^i$ at point $s_i$ for each $1\le i\le n$ is
 $N(\alpha^1,\ldots,\alpha^n)$, counted with
 multiplicity. \Qed
\end{cor}

Let us
consider some special cases.
A curve $\varphi$ has a {\it cusp} at $s$ if the center of projection $\Lambda$
meets the tangent line to the rational normal curve $\gamma$ at $\gamma(s)$.
The corresponding Schubert cycle has codimension $r$, and so a
rational curve of degree $d$ in ${\mathbb P}^r$ has at most
$(d{-}r)(r{+}1)/r$ cusps.
A curve $\varphi$ is {\it simply ramified}
(or has a {\it flex}) at $s$ if
$\varphi(s),\varphi'(s),\ldots,\varphi^{(r-1)}(s)$ span a hyperplane in
${\mathbb P}^r$ which contains $\varphi^{(r)}(s)$.
This occurs if the center of the projection $\Lambda$ meets the osculating
$r$-plane $F_r(s)$ in a point.
In this case the codimension is
$1$ and so a rational curve with only flexes has $(r{+}1)(d{-}r)$ flexes.

\section{Maximally inflected curves}

\subsection{Existence}
We ask the following question:

\begin{ques}\label{ques:motivate}
 Given $d>r$, ramification data $\alpha^1,\ldots,\alpha^n$ for degree $d$
 rational curves in ${\mathbb P}^r$, and distinct points
 $s_1,\ldots,s_n\in \RP^1$, are there  any
 {\it real} rational curves
 $\varphi:{\mathbb P}^1\rightarrow{\mathbb P}^r$ of degree $d$
 with ramification $\alpha^i$ at the point $s_i$
for each
 $i=1,\ldots,n$?
\end{ques}

Corollary~\ref{C:number_curves} guarantees the existence of
$N(\alpha^1,\ldots,\alpha^n)$ such \emph{complex} rational curves, counted with
multiplicity.
The point of this question is, if the ramification occurs at real points in
the domain ${\mathbb P}^1$, are any of the resulting curves real?
Call a real rational curve whose ramifications occurs only at real points
a {\it maximally inflected} curve.
The answer to Question~\ref{ques:motivate} is unknown in general.
There are however many cases for which the answer is yes.
Moreover, one, still open, conjecture
of Shapiro and Shapiro in the real
Schubert calculus would imply
a very strong resolution of
Question~\ref{ques:motivate}.
We formulate that conjecture in terms of maximally inflected curves.

\begin{conj}(Shapiro-Shapiro)\label{conj:shapiro}
 Let $d>r$ be integers and $\alpha^1,\ldots,\alpha^n$ be
 ramification data for degree $d$ rational curves in
 $\mathbb{P}^r$.
 For any choice of distinct real points
 $s_1,\ldots,s_n\in{\mathbb R}{\mathbb P}^1$, {\bf every}
 curve
 $\varphi:{\mathbb P}^1\rightarrow{\mathbb P}^r$ of degree $d$ with
 ramification $\alpha^i$ at $s_i$ for each
 $i=1,\ldots,n$ is real.
\end{conj}

Eremenko and Gabrielov~\cite{EG02} proved this conjecture in the cases
when $r$ is 1 or $d{-}2$.
It is also known to be true for a few sporadic cases of
ramification data and some cases of the conjecture imply others.
There is also substantial computational evidence in support of
Conjecture~\ref{conj:shapiro} and no known counterexamples.
For an account of this conjecture, see~\cite{So00b} or the web
page~\cite{So_shapiro-www}.

A ramification sequence
is {\it special}\/ if it is one of the following.
$$
  (0,1,\ldots,r{-}1,\,r{+}a)\qquad\mbox{ or }\qquad
  (0,1,\ldots,r{-}a,\,r{-}a{+}2,\ldots,r{+}1)\,,
$$
for some $a>0$.
When $a=1$ these coincide and give the sequence of a flex.
We can guarantee the existence of maximally inflected
curves with these special singularities.

\begin{thm}\label{thm:simpleRam}
 If $\alpha^1,\ldots,\alpha^n$ are ramification data for degree
 $d$ rational curves in ${\mathbb P}^r$ with at most one
 sequence $\alpha^i$ not special,
 then there {\bf exist}
 distinct points $s_1,\ldots,s_n\in \RP^1$ such that
 there are exactly $N(\alpha^1,\ldots,\alpha^n)$ real
 rational curves
 of degree $d$ which
 have ramification $\alpha^i$ at $s_i$
 for each $i=1,\ldots,n$.
\end{thm}

The point of this theorem is that there are the expected number of such
curves, and {\it all} of them are real.\smallskip

\noindent{\bf Proof. }
Let $\alpha^1,\ldots,\alpha^n$ be a ramification data
for degree $d$ rational curves in ${\mathbb P}^r$ with at most one
sequence $\alpha^i$ not special.
These are types of Schubert varieties in the Grassmannian of
($d{-}r{-}1$)-planes in ${\mathbb P}^d$ with at most one not special.
By Theorem~1 of~\cite{So99a}, there exist points
$s_1,\ldots,s_n\in \RP^1$ so that the
Schubert varieties $ X_{\alpha^i}\Fdot(s_i)$
defined by flags osculating the rational normal curve at points
$s_i$ intersect transversally with all points of intersection real.
Every ($d{-}r{-}1$)-plane $\Lambda$ in such an intersection is a center of
projection giving a maximally inflected real rational curve $\varphi$
with ramification sequence $\alpha^i$ at $s_i$
for each $i=1,\ldots,n$.
\QED

The proof in~\cite{So99a} gives points of ramification that are
clustered together.
More precisely, suppose that the
sequences $\alpha^2,\ldots,\alpha^n$ are
special and only (possibly) $\alpha^1$ is not special.
By
$$
 \forall s_2\ll s_3\ll\cdots\ll s_n
$$
we mean
 \begin{multline*}
  \forall s_2>0\
  \exists\calN_3>0,\mbox{ such that }\forall s_3>\calN_3\ \exists\calN_4>0,\
  \cdots\ \\
  \forall s_{n-1}>\calN_{n-1}\ \exists\calN_n>0,\mbox{ such that }\forall s_n>\calN_n\,.
 \end{multline*}
Then the choice of points $s_i$ giving all curves
real in Theorem~\ref{thm:simpleRam} is
$$
  s_1\ =\ \infty\quad\mbox{ and }\quad
  \forall s_2\ll s_3\ll\cdots\ll s_n\,.
$$
In short, Theorem~\ref{thm:simpleRam} only guarantees the existence of many
maximally inflected curves when almost all of the ramification is special, and
the points of ramification are clustered together in this way.
A modification of the proof in~\cite{So99a} (along
the lines of the Pieri homotopy algorithm in~\cite{HSS98}) shows
that there can be
two ramification indices that are not special, and
also two `clusters' of ramification points.

Eremenko and Gabrielov~\cite[Corollary 4]{EG01b} prove the following.

\begin{prop}[Eremenko-Gabrielov]\label{Prop:EG}
 Suppose $0<r<d$ and $d$ is even.
 Let $m:=\max\{r{+}1,d{-}r\}$ and $p:=\min\{r{+}1,d{-}r\}$.
 Set $M:=(r{+}1)(d{-}r)$.
 Then for any distinct points $s_1,\ldots,s_M\in\RP^1$, there exist
 at least
 \begin{equation}\label{E:E-G}
   \frac{1!2!\cdots(p{-}1)!(m{-}1)!(m{-}2)!\dotsb(m{-}p{+}1)!
   \left(\frac{pm}{2}\right)!}%
   {(m{-}p{+}2)!(m{-}p{+}4)!\cdots(m{+}p{-}2)!%
   \left(\frac{m{-}p{+}1}{2}\right)!\left(\frac{m{-}p{+}3}{2}\right)!%
   \cdots\left(\frac{m{+}p{-}1}{2}\right)!}
 \end{equation}
 maximally inflected curves of degree $d$ in ${\mathbb P}^r$ with
 flexes at the points $s_1,\ldots,s_M$.
\end{prop}

The constructions of Section~\ref{sec:sing-defs}
show that there exist maximally inflected
plane curves of any degree $d$ with any number up to
$d{-}2$ cusps for some choices of ramification points not clustered together.

\begin{rem}
 It is worthwhile to compare this number~\eqref{E:E-G} of Eremenko-Gabrielov
 to the number of (complex) curves with the same flexes, as computed by
 Schubert~\cite{Sch1886c}:
 \begin{equation}\label{eq:dmp}
  \frac{1!\,2!\cdots(p\! -\!2)!\,(p\! -\!1)!\cdot(mp)!}
  {(m)!\,(m\! +\!1)!\cdots (m\! +\!p\! -1)!}\ .
 \end{equation}
 The ratio of~\eqref{eq:dmp} to~\eqref{E:E-G} is
 \[
   r(m,p)\ :=\
     \frac{\left(\frac{m{-}p{+}1}{2}\right)!\left(\frac{m{-}p{+}3}{2}\right)!%
    \cdots\left(\frac{m{+}p{-}1}{2}\right)!\cdot (mp)!}%
    {(m{-}p{+}1)!(m{-}p{+}3)!\dotsb (m{+}p{-}1)!\cdot
    \left(\frac{pm}{2}\right)!}\ .
 \]
 According to Stirling's formula, $\log r(m,p)$ grows
  as
\[
  mp\log\frac{mp}{2}\;-\;(m-p+1)\log\frac{m-p+1}{2}\;-\;
  \dotsb \;-\;(m+p-1)\log\frac{m+p-1}{2}\,.
\]
 Fixing $p$ and letting $m$ grow (that is, fixing the target
 $\mathbb{P}^r$ and letting the degree $d$ grow),
 we obtain the asymptotic value of
 $\frac12 mp\log p$ for $\log r(m,p)$.
 Similar arguments show that, asymptotically,  the logarithm of
 Schubert's number grows like $mp\log p$.
 Thus we see that the number of real solutions guaranteed by
 Eremenko-Gabrielov is approximately the square root of the total number of
 solutions,
 in this asymptotic limit.

 This asymptotic result is reminiscent of (but different than) results of
 E.~Kostlan~\cite{Ko93} and M.~Shub and S.~Smale~\cite{SS92} concerning
 the expected number of real solutions to a system of polynomial equations.
 The set of real polynomial systems on $\mathbb{RP}^n$
 \begin{equation}\label{eq:system}
   f_1(x_0,x_1,\ldots,x_n)\ =\  f_2(x_0,x_1,\ldots,x_n)
   \ =\ \dotsb  \ =\  f_n(x_0,x_1,\ldots,x_n)\ =\  0\,,
 \end{equation}
 where $f_i$ is homogeneous of degree $d_i$ is parameterized by
 $\mathbb{RP}^{d_1}\times\mathbb{RP}^{d_2}\times\dotsm\times\mathbb{RP}^{d_n}$.
 Integrating the number of real roots of a system~\eqref{eq:system}
 against the Fubini-Study probability measure on this space of systems
 gives the expected number of real roots
 \[
    \left(d_1\cdot d_2\dotsm d_n\right)^{\frac{1}{2}}\,,
 \]
 the square root of the expected number $d_1\cdot d_2\dotsm d_n$
 of complex roots.
\end{rem}

\subsection{Deformations}\label{sec:MIC}
Deforming a maximally inflected curve
$\varphi\colon{\mathbb P}^1\to{\mathbb P}^r$ means deforming the positions
of its ramifications in $\mathbb{RP}^1$.
A set $S:=\{s_1(t),\ldots,s_n(t)\}$ of continuous functions
$s_i\colon[0,1]\to\mathbb{RP}^1$ where, for each $t$,
the points $s_1(t),\ldots,s_n(t)$ are distinct is an
{\it isotopy} between
$\{s_1(0),\ldots,s_n(0)\}$ and $\{s_1(1),\ldots,s_n(1)\}$.
Given such an isotopy $S$, suppose $\varphi$ has ramification
$\alpha^i$ at $s_i(0)$, for $i=1,\ldots,n$.
A continuous family $\varphi_t$ for $t\in[0,1]$ of maximally inflected curves
with $\varphi_0=\varphi$ and where $\varphi_t$ has ramification
$\alpha^i$ at $s_i(t)$, for $i=1,\ldots,n$, is a {\it deformation of
$\varphi$ along $S$} and $\varphi_1$ is a deformation of $\varphi$.
A maximally inflected curve $\varphi$
is said to admit {\it arbitrary deformations} if
$\varphi$ has a deformation along any isotopy
$S$ deforming the ramification points
of $\varphi$.
Since reparameterization by a projective transformation
of $\mathbb{RP}^1$ does not change the image of $\varphi$,
a basic question in the topology of maximally inflected curves is to
classify them up to deformation and reparameterization.

Experimental evidence and geometric intuition suggest that maximally inflected
curves admit arbitrary deformations, in a strong sense
that we make precise in
Theorem~\ref{thm:Implies}(2) below.
First we state a non-degeneracy conjecture.

\begin{conj}(Conjecture~5 of~\cite{So99a})\label{conj:non-degen}
 Let $d>r$ and $\alpha^1,\ldots,\alpha^n$ be ramification data
 for degree $d$ rational curves in ${\mathbb P}^r$.
 If $s_1,\ldots,s_n\in\mathbb{RP}^1$ are distinct, then there are
 exactly $N(\alpha^1,\ldots,\alpha^n)$
 {\bf complex}
 rational curves of degree $d$ in
 ${\mathbb P}^r$ with ramification
 sequence $\alpha^i$ at $s_i$ for each
 $i=1,\ldots,n$.
\end{conj}

That is, when the points $s_i$ are real,
there should be no multiplicities in Corollary~\ref{C:number_curves}.

\begin{thm}\label{thm:Implies}
\mbox{ }
\begin{enumerate}
\item Suppose Conjecture~$\ref{conj:non-degen}$ holds in all cases when
      the ramification consists only of flexes.
      Then Conjecture~$\ref{conj:shapiro}$ holds for all ramification data.
\item If Conjecture~$\ref{conj:non-degen}$ holds for all ramification data,
      then given $s_1,\ldots,s_n\in\mathbb{RP}^1$, each of the
      $N=N(\alpha^1,\ldots,\alpha^n)$ maximally inflected curves with
      ramification sequence $\alpha_i$ at $s_i$ admit arbitrary deformations,
      and the $N$ deformations along a given isotopy are distinct at each point
      $t\in[0,1]$.
\end{enumerate}
\end{thm}

\noindent{\bf Proof. }
Statement 1 is just Theorem~6 of~\cite{So99a}, adapted to maximally inflected
curves.

For the second statement, let $\{s_1(t),\ldots,s_n(t)\}$ be continuous functions
$s_i\colon[0,1]\to\mathbb{RP}^1$ with the points $s_1(t),\ldots,s_n(t)$ distinct
for each $t$.
For each $t\in[0,1]$, consider the intersection of
Schubert varieties
$$
  \bigcap_{i=1}^n  X_{\alpha^i}\Fdot(s_i(t))\quad\subset\quad
  \mbox{\it Grass}_{d-r-1}{\mathbb P}^d \times[0,1]\,.
$$
By Conjecture~\ref{conj:non-degen}, this consists of exactly
$N:=N(\alpha^1,\ldots,\alpha^n)$ points for each $t$, and by the first
statement, they are all real.
Since $N$ is the degree of such an intersection, it must be transverse,
and so the totality of these intersections define $N$ continuous and
non-intersecting arcs in the real points of
$\mbox{\it Grass}_{d-r-1}{\mathbb P}^d \times[0,1]$.
Each arc is a deformation of the corresponding curve at $t=0$, which
proves the second statement.
\QED

We remark that Conjecture~\ref{conj:non-degen} is not true if we remove
the restriction that the points $s_i$ are real.
For example, there is a unique map $\mathbb{P}^1\to\mathbb{P}^1$ of degree
$3$ with simple critical points (simple ramification)
$0,1,\omega$, and $\omega^2$, where $\omega$ is a primitive third root of unity.
Given 4 critical points in general position, there will be two such maps,
which happen to coincide for this particular choice.
If however, all critical points are real, then there will always be 2 such
maps.
(Details are found in ~\cite[Section~9]{EH83}.)

Conjecture~\ref{conj:non-degen} is known to hold whenever it has been
tested.
This includes when $r=1$ or $d{-}2$ and the ramification
consists only of
flexes~\cite{EG02} and for some other ramification data.
The case $r=2$ and $d=4$ of plane quartics with arbitrary ramification
was shown earlier, by direct computation~\cite[Theorem 2.3]{So00b}.

\subsection{Constructions using duality}
We give two elementary constructions of new maximally inflected curves
from old ones, each invoking a different notion of duality for these curves.
The first relies on Grassmann duality---the Grassmannian of
($d{-}r{-}1$)-planes in ${\mathbb P}^d$ is isomorphic to the Grassmannian of
$r$-planes in ${\mathbb P}^d$.
The second construction relies on projective duality and has new
implications for the
Shapiro conjecture.

Let $\Fdot(s)$ be the flag of subspaces in ${\mathbb P}^d$ osculating the
rational normal curve $\gamma(s)=F_0(s)$.
Then its dual flag  $\hFdot(s)$ is the flag of
subspaces osculating the rational normal curve $F_{d-1}(s)$ in the dual
projective space.
In particular, $\hat{F}_i(s)$ is dual to $F_{d-i}(s)$.
Consider a possible ramification sequence $\alpha=(0,\alpha_1,\ldots,\alpha_r)$
for degree $d$ curves in $\mathbb{P}^r$,
with the {\bf additional restriction} that $\alpha_r<d$.
Recall from Section~\ref{S:notation} that a
($d{-}r{-}1$)-plane $\Lambda$ lies in the Schubert variety
$ X_\alpha\Fdot(s)$ if and only if
the dual $r$-plane lies
in the Schubert variety $ X_{\hat{\alpha}}\hFdot(s)$.
Here, the sequence $\hat{\alpha}$ is defined from $\alpha$ by~\eqref{E:hat-def}.
Identifying ${\mathbb P}^d$ with its dual space, this gives a
bijection between the two algebraic sets
 \begin{equation}\label{eq:grass-duality}
  \left\{\begin{minipage}[c]{1.7in}
      Curves $\varphi$ of degree $d$ in ${\mathbb P}^r$ with
      ramification
     $\alpha^i$ at $s_i$ for $i=1,\ldots,n$.
   \end{minipage}\right\}
   \qquad \Longleftrightarrow \qquad
  \left\{\begin{minipage}[c]{1.85in}
      Curves $\varphi$ of degree $d$ in ${\mathbb P}^{d-r-1}$ with
      ramification  $\widehat{\alpha^i}$ at $s_i$ for $i=1,\ldots,n$.
   \end{minipage}\right\}
 \end{equation}
For a maximally inflected curve in ${\mathbb P}^r$, this
{\it Grassmann duality} gives a
maximally inflected curve of
the same degree in a possibly different projective space
ramified at the same points,
but with different ramification sequences.
\smallskip

Now let us give another construction involving 
the usual projective duality.
Given a curve $C$ in ${\mathbb P}^r$, its dual curve $\check{C}$ is the curve
in  the dual projective space $\check{\mathbb P}^r$ which is the closure
of the set of hyperplanes
osculating $C$ at general
points.
If $C$ is the image of a maximally inflected curve $\varphi$, then
$\check{C}$ is also rational with a parameterization
$\check{\varphi}$ induced from $\varphi$.
The ramification points of $\check{\varphi}$ coincide with
those of $\varphi$, but the ramification indices and degree of $\check{\varphi}$
will be different.
We compute this degree and the transformation of ramification indices
and thus show that $\check{C}$
is a maximally inflected curve.

The osculating hyperplane
at a general point $\varphi(s)$ of $C$ is determined
by the 1-form $\psi(s)$ whose $i$th coordinate is the determinant
 \begin{equation}\label{eq:matrix}
    \left(\varphi^{(a)}_b(s)\right)
      ^{a=0,1,\ldots,r-1}_{b=0,1,\ldots,\widehat{r{-}i},\ldots,r}\ ,
 \end{equation}
where $\varphi=(\varphi_1,\ldots,\varphi_r)$ and the derivatives are taken
with respect to some local coordinate of ${\mathbb P}^1$ at $s$,
and $\widehat{r{-}i}$ indicates that the index $r{-}i$ is omitted from
the list.
We may assume for simplicity
that in this local coordinate $\varphi_i$ has degree $d{-}i$.
Thus the 1-form $\psi$ has degree $r(d{-}r{+}1)$ and
its coordinates define
a linear series of dimension $r$ and degree $r(d{-}r{+}1)$ on
${\mathbb P}^1$.
In general, this linear series will have base points and so the degree of the
resulting map will be less than $r(d{-}r{+}1)$.

Let us
determine the base point divisor and the
ramification of the map determined by $\psi$.
Suppose that $\alpha$ is the ramification sequence of $\varphi$ at a point
$s$  of ${\mathbb P}^1$.
We may assume that $\varphi_i$ vanishes to order $\alpha_i$ at $s$, and a
calculation shows that the determinant~(\ref{eq:matrix}) vanishes to order
$|\alpha|+r-\alpha_{r-i}$ at $s$.
Thus $s$ has multiplicity $|\alpha|+r-\alpha_r$ in the base point divisor.
Removing this base point divisor from the map $\psi$
gives a map $\check{\varphi}$ whose $i$th coordinate
vanishes to order $\alpha_r-\alpha_{r-i}$ at $s$,
and so the ramification
sequence
of the dual curve at a point $s$ where
$\alpha=\alpha(s)$ is
 \begin{equation}\label{eq:dual-ram}
   \check{\alpha}\ :=\
  (0,\;   \alpha_r{-}\alpha_{r-1},\; \ldots,\;
   \alpha_r{-}\alpha_1,\;   \alpha_r)\;.
 \end{equation}
Thus the map determined by $\psi$, and hence the dual curve, has degree
 \begin{equation}\label{eq:deg-dual}
   r(d-r+1) - \sum_{s\in{\mathbb P}^1} (|\alpha(s)|+r-\alpha(s)_r )\,,
 \end{equation}
where $\alpha(s)$ is the ramification sequence
at $s$, which is typically
$(0,1,\ldots,r)$.
The following theorem is immediate.

\begin{thm}
  Let $d>r$
  and suppose that $\alpha^1,\ldots,\alpha^n$
  is ramification data for rational curves of degree
  $d$ in   ${\mathbb P}^r$.
  Then $\check{\alpha^1},\ldots,\check{\alpha^n}$ is
  ramification data for rational curves in ${\mathbb P}^r$ of degree
$$
   m\ :=\ r(d-r+1) -
           \sum_{s\in{\mathbb P}^1} (|\alpha(s)|+r-\alpha(s)_r)\,.
$$
\begin{enumerate}
 \item
  For any choice of distinct points $s_1,\ldots,s_n\in\RP^1$, there is a
  one to one correspondence between maximally inflected curves
  of degree $d$ with ramification $\alpha^i$ at
  $s_i$ for $i=1,\ldots,n$ and
  maximally inflected curves of degree $m$ with ramification
  $\check{\alpha^i}$ at $s_i$ for $i=1,\ldots,n$.
 \item
   Conjecture~$\ref{conj:shapiro}$ holds for $r,d$, and the sequences
   $\alpha^1,\ldots,\alpha^n$ if and only if it holds for the integers
   $r,m$, and the sequences
   $\check{\alpha^1},\ldots,\check{\alpha^n}$.\ \ \QED
\end{enumerate}
\end{thm}

We compute the degree $m$ of the dual curve when $r=2$.
Suppose the original curve has ramification indices
$\alpha^1,\ldots,\alpha^n$.
The formula~(\ref{eq:deg-dual}) becomes
 \begin{equation}\label{eq:plane-dual}
  m\ =\ 2(d-1) - \sum_i( |\alpha^i|+2-\alpha^i_2)\ =\
         4-d+\sum_i( \alpha^i_2-2)\ .
 \end{equation}
%

\section{Maximally inflected plane curves}
For plane curves, we have $r=2$.
Suppose, in the beginning, for sake of discussion that
we have a rational plane curve
whose only ramifications are cusps and flexes, and
whose only other singularities are 
ordinary double points.
Such a curve of degree $d$ with $\kappa$ cusps
has $\iota=3(d-2)-2\kappa$
flexes.
This Pl\"ucker formula follows, for example, from Proposition~\ref{prop:ram}
since a cusp has ramification
sequence $(0,2,3)$ and a flex
$(0,1,3)$, and these have weights 2 and 1 respectively.
Since the curve is of
genus zero, it must have
$g_\kappa:=\binom{d-1}{2}-\kappa$ double points.

Suppose now that the curve is real.
Each visible (in the real part of $\mathbb P^2$)
node is either a
{\it real node}---a real ordinary double point with
real tangents, or a {\it solitary point}---a real ordinary double
point with complex conjugate tangents.
All other nodes are {\it complex nodes}; they occur in complex conjugate pairs.
Let $\eta$ be the number of real nodes, $\delta$ the number of solitary
points, and $c$ the number of complex nodes.

\begin{figure}[htb]
 $$
    \epsfysize=70pt\epsfbox{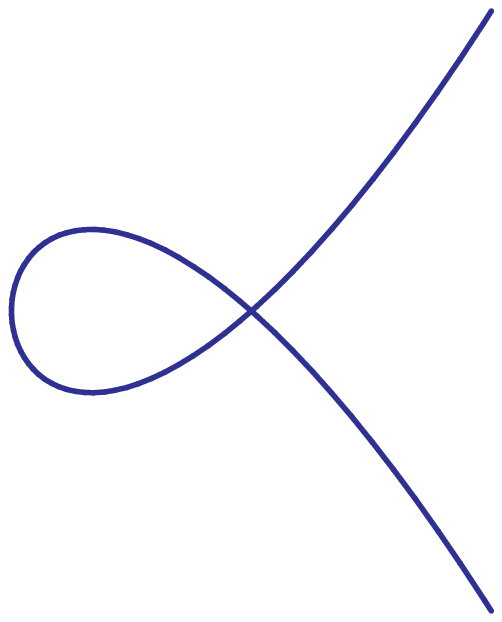}\qquad\qquad\qquad %
    \epsfysize=70pt\epsfbox{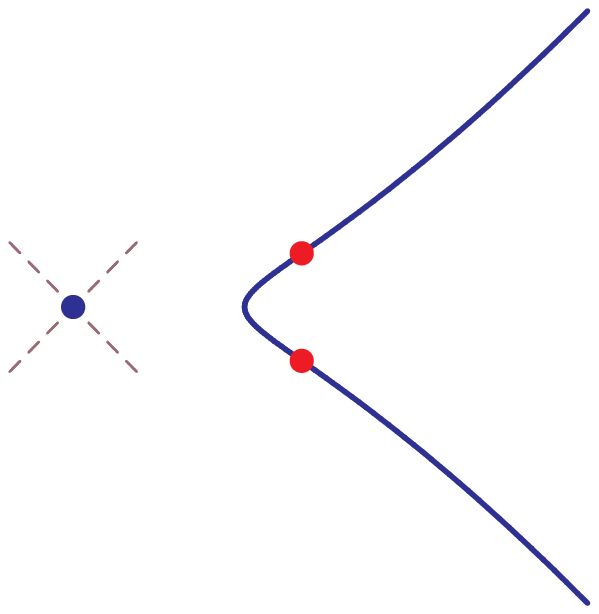}\qquad\qquad\qquad %
    \epsfysize=70pt\epsfbox{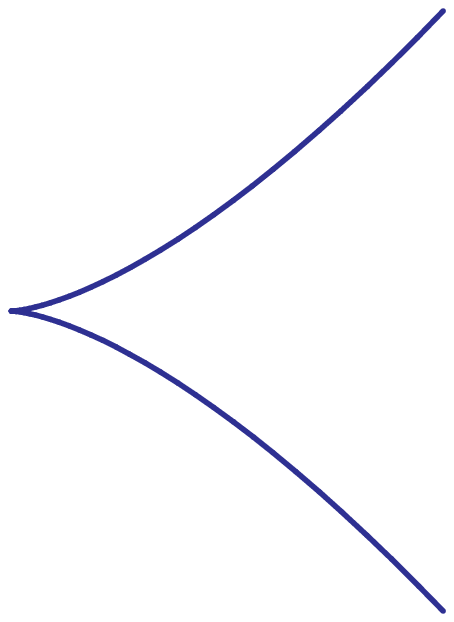}
 $$
 \caption{Real rational cubics}\label{fig:cubics}
\end{figure}

Up to projective transformation and reparametrization,
there are only three real rational plane cubic curves.
They are represented by the
equations $y^2=x^3+x^2$, \,$y^2=x^3-x^2$, and $y^2=x^3$,
and they have the shapes shown in Figure~\ref{fig:cubics}.
All three have a real flex at infinity and are singular at the origin.
The first has a real node and no other real flexes, the second has
a solitary point and two real flexes at
$(\frac{4}{3},\pm\frac{4}{3\sqrt{3}})$ (we indicate these with dots
and the complex conjugate tangents at the
solitary point with dashed lines), and the
third has a real cusp.
The last two are maximally inflected, while the
first is not.
In general, as is shown below,
the number of real nodes is restricted for maximally inflected
curves.

In what follows
we consider maximally inflected curves with arbitrary
ramifications and define a
{\it solitary point} of a real rational curve
$\varphi\colon{\mathbb P}^1\to{\mathbb P}^2$
to be a pair of distinct complex conjugate
points $s,\overline{s}\in{\mathbb P}^1$ with
$\varphi(s)=\varphi(\overline{s})$, which
necessarily represents
a point in
$\mathbb{RP}^2$.
Let $\delta$ be the number of such solitary points.
We similarly define
solitary
bitangents to be solitary points of the
dual curve, and let $\tau$ be their number.

\begin{thm}\label{thm:bounds}
Let $\varphi$ be a maximally inflected plane curve of degree $d$ with
ramification $\alpha^1,\ldots,\alpha^n$.
Then the numbers $\tau$ of
solitary bitangents and $\delta$ of  solitary points satisfy
$$
  \delta\ =\ \tau + d-2-\sum_i(\alpha^i_1-1)\,.
$$

\end{thm}

\noindent{\bf Proof. }
Let $C$ be the real points of image of the curve $\varphi$ and
$\check{C}$ be its dual curve,
the real points of the image of $\check{\varphi}$.
The generalized Klein formula due to Schuh~\cite{Sch1904} (see
also~\cite{Viro})  gives the relation
 \begin{equation}\label{eq:Schuh}
   m + \sum_{z\in C}( \mu_z-r_z)\ =\
   d + \sum_{z\in\check{C}}( \mu_z-r_z)\,,
 \end{equation}
where $\mu_z$ is the multiplicity of a point $z$ in a curve and $r_z$ is the
number of real branches of the curve at $z$.

We first evaluate the sum
 \begin{equation}\label{eq:curvesum}
   \sum_{z\in C}( \mu_z-r_z)\,.
 \end{equation}
The multiplicity $\mu_z$ at a point $z\in C$ is the local intersection multiplicity
(at $z$) of the curve $C$ with a general linear form $f$ vanishing at $z$.
This is the sum over all $s\in\varphi^{-1}(z)$ of the order of vanishing
at $s$ of the pullback $\varphi^*(f)$.
This order of vanishing is $\alpha(s)_1$, where $\alpha(s)$ is the
ramification sequence of $\varphi$ at $s$.
There are two cases to consider:
Either $s$ is real ($s\in{\mathbb R}{\mathbb P}^1$) or it is not.

In the first case, the image of a neighborhood of $s$ in
${\mathbb R}{\mathbb P}^1$ is a branch of $C$ at $z$, so the
contribution of $s$ to the
sum~(\ref{eq:curvesum}) is $\alpha(s)_1-1$.
Since $\alpha(s)_1-1$ vanishes except at some points of ramification,
the contribution of points in ${\mathbb R}{\mathbb P}^1$
to~(\ref{eq:curvesum}) is the sum
$$
   \sum_{i=1}^n(\alpha^i_1-1)\,.
$$
In the second case, the complex conjugate $\overline{s}$ of $s$ is also
in $\varphi^{-1}(z)$.
Since $\varphi$ is maximally inflected, it is unramified at these points,
so $\alpha(s)_1=\alpha(\overline{s})_1=1$.
Thus each solitary point contributes 2 to the sum~(\ref{eq:curvesum}).
Combining these observations gives
$$
  \sum_{z\in C}(\mu_z-r_z)\ = \ 2\delta + \sum_{i=1}^n(\alpha^i_1-1)\,.
$$

If $s\in{\mathbb P}^1$ and $\varphi$ has ramification $(0,a,b)$ at $s$,
then $\check{\varphi}$ has ramification $(0,b{-}a,b)$ at $s$,
by~(\ref{eq:dual-ram}).
Thus
$$
  \sum_{z\in \check{C}}(\mu_z-r_z)\ = \ 2\tau +
  \sum_{i=1}^n(\alpha^i_2-\alpha^i_1-1)\,.
$$
Substituting these expressions and
the formula~(\ref{eq:plane-dual}) for $m$ into Schuh's formula~(\ref{eq:Schuh}),
we obtain
$$
  4-d+\sum(\alpha^i_2-2) +2\delta + \sum(\alpha^i_1-1)\ =\
  d+2\tau + \sum(\alpha^i_2-\alpha^i_1-1)
$$
or $2\delta=2d-4+2\tau-2\sum(\alpha^i_1-1)$, which completes the proof.
\QED

\begin{cor}\label{cor:generalbounds}
Let $\varphi$ be a maximally inflected plane curve of degree $d$.
Let $\delta, \eta$, and $c$ be the
numbers of solitary points, real nodes, and
complex nodes (respectively).
Then,
$$
  \begin{array}{rcccl}
   d{-}2{-}\sum(\alpha^i_1-1)&\leq& \delta&\leq&
    \binom{d-1}{2}-\frac12\sum(\alpha^i_1-1)(\alpha^i_2-1)\\
   0&\leq&\eta{+}2c&\leq&\binom{d-2}{2}
   -\frac12\sum(\alpha^i_1-1)(\alpha^i_2-1)+\sum(\alpha^i_1-1).
\rule{0pt}{25pt}
  \end{array}
$$
\end{cor}

\noindent{\bf Proof. }
Since $\tau\ge 0$, the lower bound for $\delta$
follows from the formula of Theorem~\ref{thm:bounds}.
The upper bound is, in fact, an upper bound for
the number of {\it virtual double points}
which can appear, outside a given ramification,
on a curve of degree $d$.
Namely, the total number of virtual double points (including the virtual
double points accounted for by the fixed ramification points)
is equal to the genus of a
nonsingular curve of degree $d$ and the number of
virtual double points contained in a ramification
point equals $\frac12 (\mu+r-1)$,
where $\mu$ is the Milnor number and $r$ is the number of local
branches
(this formula can be found in~\cite{Ba1893};
for a more modern treatment and generalizations see~\cite{Ko76}).
Now, it remains to notice that $\mu\ge(\alpha_1{-}1)(\alpha_2{-}1)$
at a ramification point of type $(0,\alpha_1,\alpha_2)$.

Since a maximally inflected curve has genus 0, the genus formula gives
 \begin{equation}\label{eq:genus}
  \delta + \eta + 2c\ \le \  g_\alpha\ ,
 \end{equation}
where $g_\alpha$ is the genus of a generic curve
with singularities prescribed by the ramification data.
The upper bound for $\eta+2c$ now follows from the
lower bound
on $\delta$ and the upper bound on the genus
$g_\alpha\le\binom{d-1}{2}- \frac12\sum(\alpha^i_1-1)(\alpha^i_2-1)$
is given by
the above bound on the number of virtual double points.
\QED

\begin{cor}\label{cor:bounds}
 Let $\varphi$ be a maximally inflected plane curve of degree $d$ with $\kappa$
 (real) cusps and $\iota=3(d-2)-2\kappa$ (real) flexes whose remaining
 singularities are nodes.
 Let $\delta, \eta$, and $c$ be the numbers of solitary points, real nodes, and
 complex nodes (respectively) of $\varphi$, which satisfy
 $\delta+\eta+2c=\binom{d-1}{2}-\kappa =:g_\kappa$.
 If $\kappa\leq d{-}2$, then these additionally satisfy
 \begin{equation}\label{E:bound1}
  \begin{array}{rcccl}
   d{-}2{-}\kappa&\leq& \delta&\leq&g _\kappa\\
   0&\leq&\eta{+}2c&\leq&{\displaystyle \binom{d-2}{2}}\rule{0pt}{25pt}
  \end{array}
 \end{equation}
 If $d{-}2\leq\kappa$ (which is at most $3(d{-}2)/2$), then $\delta, \eta$, and
 $c$ satisfy
 \begin{equation}\label{E:bound2}
  \begin{array}{rcccl}
   0&\leq& \delta&\leq&
     {\displaystyle \binom{2d-4-\kappa}{2}}\\
   {\displaystyle g_\kappa-\binom{2d-4-\kappa}{2}}
    &\leq&\eta{+}2c&\leq&g_\kappa\rule{0pt}{25pt}\medskip
  \end{array}
 \end{equation}
\end{cor}

This may be deduced from the Klein~\cite{Klein}
and Pl\"ucker~\cite{Plucker} formulas alone.\medskip

\noindent{\bf Proof. }
Since the ramification
sequence of a flex is $(0,1,3)$ and that of a cusp is
$(0,2,3)$, we have
 \[
   \sum (\alpha^i_1-1)\ =\
    \frac{1}{2}\sum(\alpha^i_1-1)(\alpha^i_2-1)\ =\ \kappa\,,
 \]
and so~\eqref{E:bound1} is just a special case of
Corollary~\ref{cor:generalbounds}.

Now suppose that $d-2\leq \kappa$ and consider the dual curve to $\varphi$,
which has degree $m:=2(d{-}1)-\kappa$ by~(\ref{eq:plane-dual}), $\kappa$ flexes
and $\iota=3(d{-}2)-2\kappa$ cusps.
The total number of double points of this curve (bitangents of $\varphi$) is
given by the genus formula
$$
  \binom{m-1}{2}-\iota\ =\ \binom{2d-3-\kappa}{2}\,-\,
                            \bigl(3(d-2)-\kappa\bigr)\,.
$$
This expression is an upper bound for the number $\tau$ of
solitary  bitangents of $\varphi$.
This gives the upper bound for $\delta$
 \begin{eqnarray*}
  \delta &\leq& \binom{2d-3-\kappa}{2}\,-\,
                     \bigl(3(d-2)-2\kappa\bigr) +d-2-\kappa\\
   &=& \binom{2d-3-\kappa}{2} -(2d-4-\kappa)\quad=\quad \binom{2d-4-\kappa}{2}\ .
 \end{eqnarray*}
Combining the bounds for $\delta$ with the genus formula
$\delta+\eta+2c=g_\kappa$ gives the bounds for
$\eta+2c$.
\QED

\begin{rem}
A fundamental question about the statement of Corollary~\ref{cor:bounds} is
whether its hypotheses are satisfied by any
maximally inflected curve with $\kappa$ cusps and $3(d{-}2)-2\kappa$ flexes.
That is, is there a curve with $\kappa$ cusps and whose other singularities
are only ordinary double points not occurring
at points of ramification?
(If there is one such curve, then the general curve
has this property.)
As the construction in Section~\ref{sec:sing-defs} shows, this is true when
$\kappa\leq d-2$.

While we believe that a generic
maximally inflected curve has only ordinary
double points not occurring at ramification points, we do not have a proof.
A difficulty is that there are few constructions of maximally inflected
curves.
\end{rem}

Corollary~\ref{cor:bounds} raises
our first question concerning the
classification of maximally inflected curves by
their topological invariants.

\begin{ques}\label{ques:one}
 Let $d,\iota,\kappa$ be positive integers with $\iota+2\kappa=3(d{-}2)$.
\begin{enumerate}
\item Which numbers $\delta$ in the range allowed by Corollary~\ref{cor:bounds}
      occur as the number of solitary points in a maximally inflected plane
      curve of degree $d$ with $\iota$ flexes and $\kappa$ cusps, and whose
      other singularities
      are all ordinary double points?
\item Given a number $\delta$ of solitary points which occurs,
      $\eta+2c=g_\kappa-\delta$.
      Which numbers $\eta$ in the range $[0,g_\kappa{-}\delta]$ with
      the same parity as $g_\kappa{-}\delta$ occur as the
      number of real nodes of such a maximally inflected curve?
\end{enumerate}
\end{ques}

\begin{rem}\label{k4.1}
For example, when $d=5$, $\kappa=4$, and $\iota=1$, we are in the case of
$\kappa\geq d-2$
 in Corollary~\ref{cor:bounds}.
The bounds give $\delta=0$ or 1.
Question~\ref{ques:one}(1)  asks: do both values
of $\delta$  occur?
If $\delta=1$, then $\eta+2c=1$, so a curve with $\delta=1$ has a single real node.
If $\delta=0$, then $\eta+2c=2$, so there are 2 possibilities for $\eta$ of 0 or 2.
Question~\ref{ques:one}(2) asks: do both values
of $\eta$  occur?
In fact, both values for $\delta$ occur, but when
$\delta=0$, only the value of 2 for $\eta$ occurs.

The first part of the latter statement follows from the results of Section 5,
as explained in Remark~\ref{k4.3}.
Two such curves are displayed in Figure~\ref{F:4cusps}.
The impossibility of the values $\delta=\eta=0$ for a rational quintic plane
curve with four cusps and one flex is
easily explained using V.~Rokhlin's theory of
complex orientations~\cite{Ro74} as extended by N.~Mishachev~\cite{Mi75} and
V.~Zvonilov~\cite{Zv83}.
Let
$f:\mathbb{P}^1\to\mathbb{P}^2$ be a real rational plane quintic
with four cusps and no solitary points or nodes.
We first perturb
it to a new real rational quintic with
$4$ real nodes instead of the cusps, and then
smooth each real node in such a way
to obtain an oval, and get in this way a
dividing real quintic $Y$ of genus $4$.
(This smoothing of a cusp is illustrated in Figure~\ref{Fig:smooth}.)
\begin{figure}[htb]
\[
   \epsfysize=70pt\epsffile{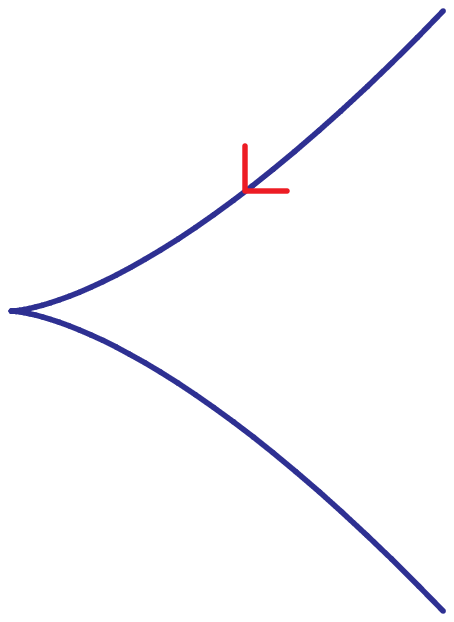}
      \quad \raisebox{32pt}{$\Longrightarrow$}\quad
   \epsfysize=70pt\epsffile{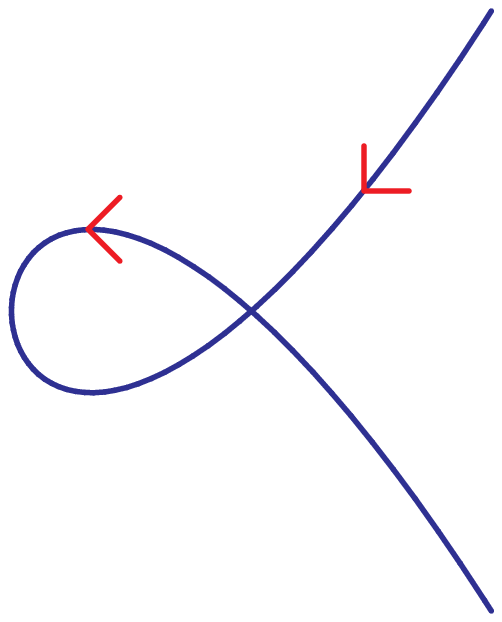}
      \quad \raisebox{32pt}{$\Longrightarrow$}\quad
   \epsfysize=70pt\epsffile{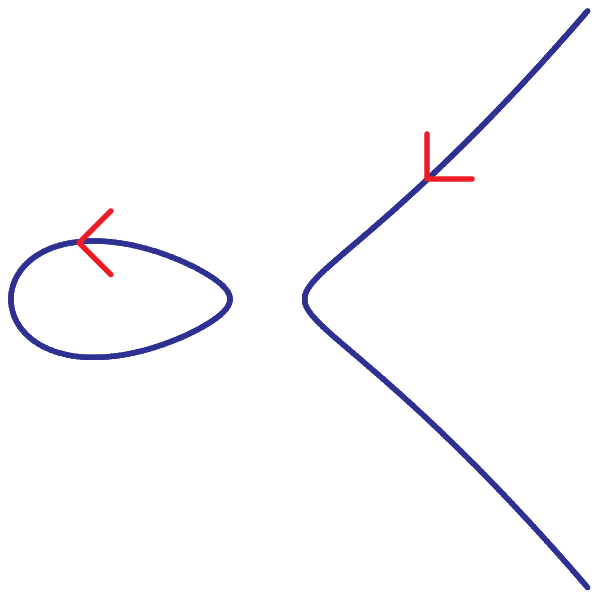}
\]
\caption{Smoothing a cusp and complex orientation}\label{Fig:smooth}
\end{figure}
Since the Rokhlin complex orientation formula extends to ``flexible'' curves,
which are only almost holomorphic near the real plane,
we do not need the
existence of such an algebraic deformation and may instead just glue
a proper local model in place of the cusps.

Choosing one component of  $Y\setminus\RE Y$
gives a complex orientation to the
ovals and the  odd branch of the curve $Y$.
This complex orientation must satisfy
 \[
    25=1+4\mu+ 4 + 4(R_+ - R_-),
 \]
where $\mu=0$ or 2 counts the intersection of the two components of
$Y\setminus\RE Y$ at the complex nodes
and $R_\pm$ counts the relative orientation of the ovals to the odd branch.
As we may see from Figure~\ref{Fig:smooth}, the ovals are all coherently
oriented with respect to the odd branch, so we have
$R_+=4$ and $ R_-=0$, which gives the contradiction.

More generally, we may ask the analog of Question~\ref{ques:one} for maximally
inflected curves with arbitrary ramification.
We leave the formulation of
this to the reader.
\end{rem}

\begin{rem}\label{k4.2}
We indicate another approach to show the impossibility of
the values $\delta=\eta=0$ for a rational quintic plane curve
with four cusps and one flex.
Consider the dual curve.
It is a rational quartic with 4 flexes and 1 cusp. By formula (12),
$\tau\ge 1$ and, thus, the dual curve has at least one
solitary point.
Trace a straight line through such a solitary point and the cusp.
By  B\'ezout's theorem, this line contains no other points
of the curve.
Hence, choosing a nearby line as the line at infinity, we
may assume the real part of the quartic
lies in an affine part of $\RP^2$.

According to the Fabricius-Bjerre formula~\cite{FB77},
$\iota+2\eta+2\kappa=2(t_+-t_-)$, where $\iota$ is the number of flexes,
$\eta$ the number of nodes, $\kappa$ is the number of cusps,
and $t_+$ (respectively, $t_-$) is the number of one-sided
(respectively, two-sided) double tangents.
Here, a line going through a cusp and tangent at another point
is counted as a double tangent as well.
Because
this
quartic curve
is connected, any double tangent has
both germs of the curve at the tangencies on the same side, that is, it is a
one-sided double tangent.
The projection from a cusp is 2-sheeted and, thus the number of
such tangents through the cusp is at most $2$.
All this together implies that
there is at least one ordinary
(i.e., not going through a cusp) double tangent. Hence, the quintic
has a real node.
\end{rem}

\begin{rem}\label{k4.3}
 Yet another proof of this impossibility of $\eta=0$ for a maximally inflected
 quintic with 4 cusps and one flex uses results of
 Section~\ref{sec:lowdegree}.
 As we have seen, the dual of such a curve is a maximally inflected rational
 quartic in $\RP^2$  with 4 flexes and a single cusp.
 By Theorem~\ref{T:quart-isotopy}, the possible topological types of such
 curves and 
 of arrangements of 
 their bitangents
 with respect to the curve
 are exhibited by the curves in the second column of
 Table~\ref{table:two} (the bitangents are not drawn there, but may be
 inferred).
 One curve has one bitangent, while the other (nodal) curve has two
 bitangents.
 Thus our original quintic must have one or two real nodes, so that $\eta=0$
 is impossible.
\end{rem}

\section{Two Constructions of Maximally Inflected Curves}
In the study of real plane curves, examples are typically constructed either
by deforming a reducible curve
(for example,
in the constructions of Harnack~\cite{Ha1876} and Hilbert~\cite{Hi1891}) or by
Viro's patchworking method deforming a curve in a reducible toric
variety~\cite{Vi83,Vi84} (see also~\cite{Ri}),
sometimes combined with Cremona transformations.
These methods were used initially to build smooth curves.
Gudkov and, then, Shustin extended these methods
to use and to obtain singular curves~\cite{Sh85,Sh99}.
We use such extensions
to construct maximally inflected plane curves of
degree $d$ with no complex nodes and the extreme numbers of 0 and
$\binom{d-2}{2}$ real nodes.

Gudkov and Shustin have established a number of
theorems extending the classical result of Brusotti~\cite{Br1921} which state
that given a singular curve satisfying a numerical condition and local models
for certain deformations of the singular points of the curve,
there exists a deformation of the singular curve where each
singularity is deformed according to the corresponding local model.
It is this approach that we use in the proof of the following theorem.
(We give the details in our proof, as the results in the literature are vastly
more general than we need and use subtle hypotheses.)

\begin{thm}\label{thm:classical}
 \mbox{ }
 \begin{enumerate}
  \item For any $d$ and $\kappa$ with $0\leq \kappa\leq d{-}2$, there exists a
        maximally inflected plane curve with $\kappa$ cusps,
        $3(d{-}2){-}2\kappa$ flexes, and $\binom{d{-}2}{2}$ real nodes
        (hence $d{-}2{-}\kappa$ solitary points).
  \item For any $d$ there exists a maximally inflected plane curve of degree
        $d$ with $3(d{-}2)$ flexes and $\binom{d{-}1}{2}$ solitary points.
        (Hence without real or imaginary nodes.)
 \end{enumerate}
\end{thm}

We prove the first statement in Section~\ref{sec:sing-defs} and the second in
Section~\ref{sec:patchworking}.
These realize the maximum and minimum possible numbers of real nodes $\eta$
and solitary points $\delta$ allowed by Corollary~\ref{cor:bounds} for
degree $d$ curves with $3(d{-}2)$ flexes.
In Section~6 we discuss some implications of the constructions in
Section~\ref{sec:sing-defs}.

\begin{rem}
  Theorem~\ref{thm:classical} asserts the existence of maximally inflected
  curves with the  given ramification, for some choices of ramification.
  Such an existence is something new when $d$ is odd.
  Eremenko and Gabrielov's result (Proposition~\ref{Prop:EG}) does not
  guarantee the existence of any maximally inflected curves
  of odd degree
  with given ramification.
  From the constructions below, we can see that the ramification points are
  not `clustered together', as in the discussion following
  Theorem~\ref{thm:simpleRam},
  so Theorem~\ref{thm:simpleRam} also does not apply.

  It follows from the relation between maximally inflected curves and the real
  Schubert calculus that these classical constructions of curves from
  Theorem~\ref{sec:sing-defs} imply the existence of real solutions to some
  problems in the Schubert calculus.
  When $d$ is odd and for the ramifications of Theorem~\ref{sec:sing-defs},
  this result is new and gives further evidence in support of the conjecture
  of Shapiro and Shapiro.

 Let us notice also that whatever is the value of $d$, even or odd, the proof
 of the Theorem gives not
just
existence but some explicit maximally inflected
curves with well controlled topology.
\end{rem}

\subsection{Deformations of Singular Curves}\label{sec:sing-defs}
For sufficiently small $\epsilon>0$ the equation
$$
   12 xy((x+y-1)^2-xy)\;-\;\epsilon(x+y)^2(x+y-1)^2\ =\ 0
$$
defines a rational quartic curve $C(\epsilon)$ with 6 flexes, one real node,
and 2 solitary points.
The curve $C(0)$ is the union of the coordinate axes and the
ellipse shown on the left below.
The curve $C(\frac{1}{20})$ is displayed on the right below.
$$
   \epsfysize=90pt\epsffile{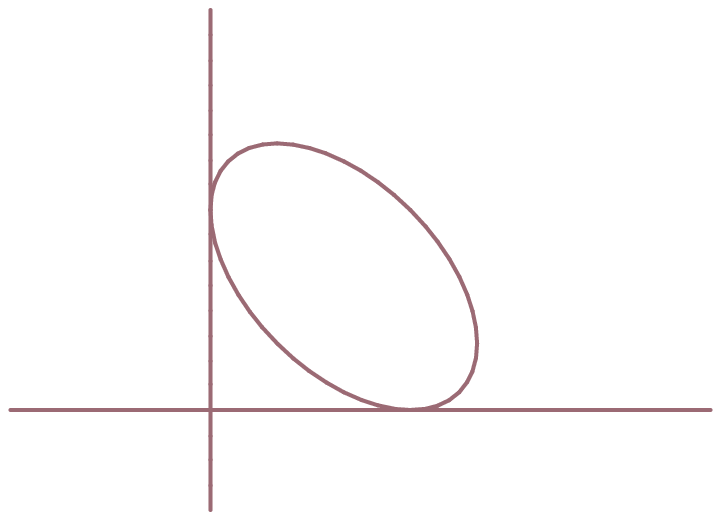}
   \qquad\qquad
   \epsfysize=90pt\epsffile{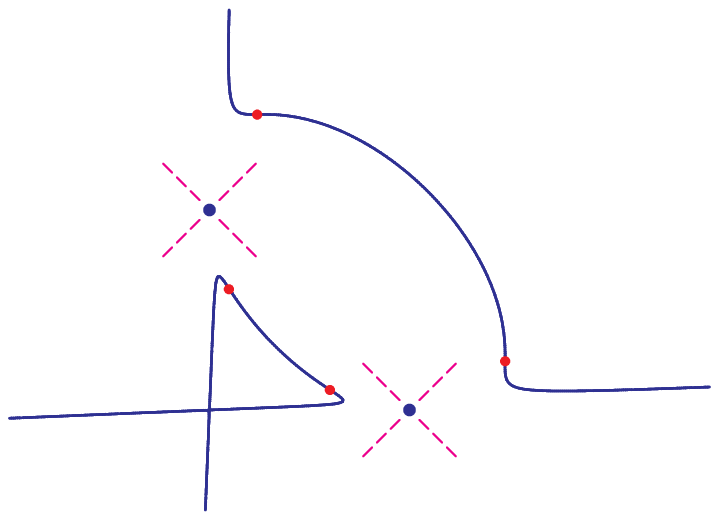}
$$

The quartic curve on the right has 2 solitary points at $(1,0)$ and $(0,1)$
and a node at the origin.
By the genus formula the curve is rational.
It also has 6 flexes.
Four are indicated by circles, and there is one more along
each of the branches close to the coordinate axes, as
the oriented geodesic curvature (which is given by the Wronskian
$\det (\phi,\phi',\phi'')$ with respect to the local orientation
defined by $\phi$, where $\phi$ is the
parametrization of the curve)
of the segment of the curve along such a branch
takes values of opposite sign at the extremal points
(close to the initial
tangency points) of these branches.
(The local orientation
changes with respect to an affine one when
the branch crosses
infinity.)
Thus any parameterization of this curve is a maximally inflected quartic with
6 flexes, 1 node, and 2 solitary points.\medskip

\noindent{\bf Proof of Theorem~\ref{thm:classical}(1).}
Fix an integer $d>2$ throughout and let $P_0$ be the union of a
nonsingular conic and any
$d{-}2$ distinct lines tangent to the conic.
Then $P_0$ is a reducible curve of degree $d$.
Each pair of tangents meet and no three meet in a point as the dual curve to
the conic is another nonsingular conic.
Thus $P_0$ has $\binom{d-2}{2}$ real nodes and $d{-}2$ other singularities
at the points of tangency.
We deform those tangency singularities while preserving the nodes.

In a neighborhood $V$ of each point of tangency, $P_0$
is isomorphic to the reducible curve $I_0$ given by the equation
$$
  y\,(y-x^2)\ =\ 0\,,
$$
in some neighborhood $U$ of the origin.
For each $t\in{\mathbb R}$, let $I_t$ be the deformation of $I_0$ defined in
$U$ by
 \begin{equation}\label{eq:Kt}
  y\,(y-x^2)\; +\; t\,x^2\ =\ 0\,.
 \end{equation}
For $t$ positive but sufficiently small, $I_t$ has a solitary point at the
origin and 2 flexes near the parabola $y=x^2$, one flex along each
branch and within $U$.
Moreover, $I_t$ lies above the $x$-axis.
Counting the Whitney index by means of the Gauss map shows the existence
of two such flexes.
\bigskip

To construct curves with cusps we replace
the local deformation model
$I_t$~(\ref{eq:Kt}) of the tangent points by the local model $K_t$ defined by
 \begin{equation}\label{eq:Lt}
  y\,(y-x^2) \;+\; t\,x^3\ =\ 0\,.
 \end{equation}
For $t>0$,  $K_t$ has a cusp at the origin and one flex near the origin.

Suppose that we can deform each tangency singularity (a tacnode) according the
local model $I_t$ while preserving the nodes so that we get
a deformation $P_t$ of the curve $P_0$ for $t\in (0,\epsilon)$ such that
(1)
$P_t$ has degree $d$, (2) $P_t$ has a node in
a neighborhood of each node of $P_0$, and (3) in a neighborhood $V$ of
each point of tangency of $P_0$, $P_t$ is isomorphic to $I_t$, in the
neighborhood $U$ of the origin.
For $0<t<\epsilon$, the curve $P_t$ has $\binom{d-2}{2}$ nodes and $d{-}2$
solitary points, by the construction.
Since it has degree $d$, it is rational.
Each solitary point contributes 2
flexes, accounting for $2(d{-}2)$ flexes.
Furthermore, there is an additional flex along each asymptote (the original
tangent lines) as the concavity of $P_t$  changes while passing through
infinity.
Thus any parameterization of the curve $P_t$ gives a maximally inflected
curve of degree $d$ with $3(d{-}2)$ flexes, $\binom{d-2}{2}$ nodes, and
$d{-}2$ solitary points.

If
we replace some local deformation models $I_t$ by $K_t$,
every tangency point where we use the local perturbation $K_t$
gives us a cusp, no solitary point, and only one flex
(the flex along the tangent
line does not appear, as the concavity of $K_t$ does not change along
that line).
Doing this for $\kappa$ of the $d{-}2$ tangency points of $P_0$ gives
a maximally inflected curve with  $\kappa$ cusps,
$3(d{-}2)-2\kappa$ flexes,
$d{-}2{-}\kappa$ solitary points, and $\binom{d{-}2}{2}$ real nodes,
proving statement (1) (the curve is rational because
$g=0$).

This simultaneous deformation of the tacnodes according to
to arbitrary independent local models
while preserving the nodes follows from the transversality
of the equisingularity strata of the singularities.
The equisingularity stratum of a node is smooth and its tangent space is given
by polynomials
vanishing
at the node.
The equisingularity stratum of a tacnode
is also smooth and its tangent space is given by polynomials
vanishing at the tacnode
that satisfy
two additional conditions:
their first derivative at a tangent direction
of the branches is zero, and the second derivatives along both branches
(taken in the same direction) are equal.

So, to check the transversality it is sufficient merely
to count the dimension of the intersection of the tangent spaces.
In our case, the intersection is
contained in the linear span
of polynomials $L_0'Q_0L_1\dotsb L_{d-3}+
L_1'Q_0L_0\dotsb L_{d-3}+\dotsb +L_0L_1\dotsb L_{d-3}Q_0'$,
where $Q_0$ is our initial nonsingular conic, $L_0,\dotsc,L_{d-3}$
our tangent lines, and $L_0',L_1',\dotsc, Q_0'$ are
equations for arbitrary lines and a conic.
So, the
condition
on the second derivatives
takes the form of a triangular linear system and, hence,
the subspace we are looking for is of dimension at most
$2(d-2)+6-(d-2)=d+4$, which implies transversality.
\QED

\begin{figure}[htb]
$$
  \epsfysize=80pt\epsfbox{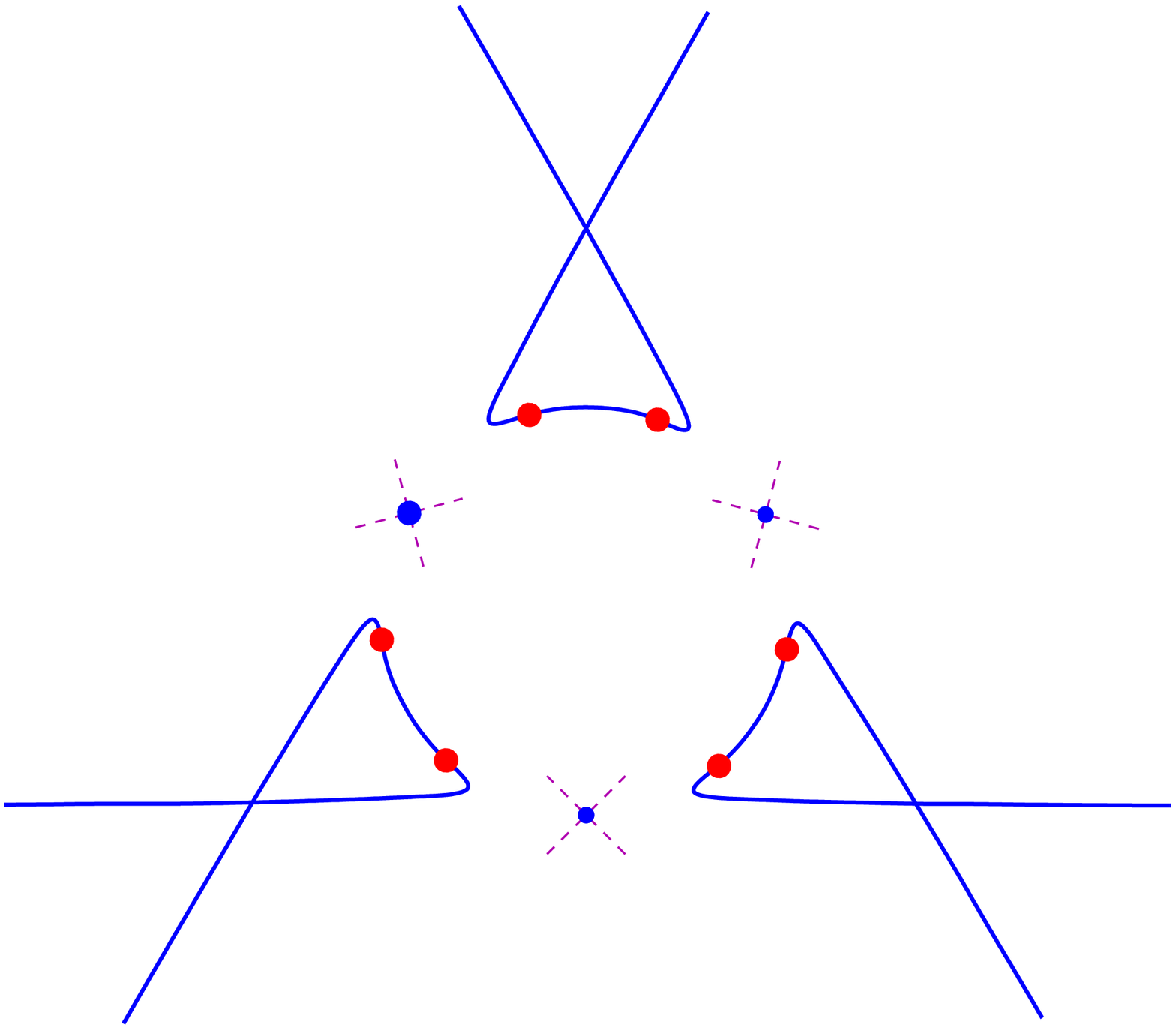}\quad
  \epsfysize=80pt\epsfbox{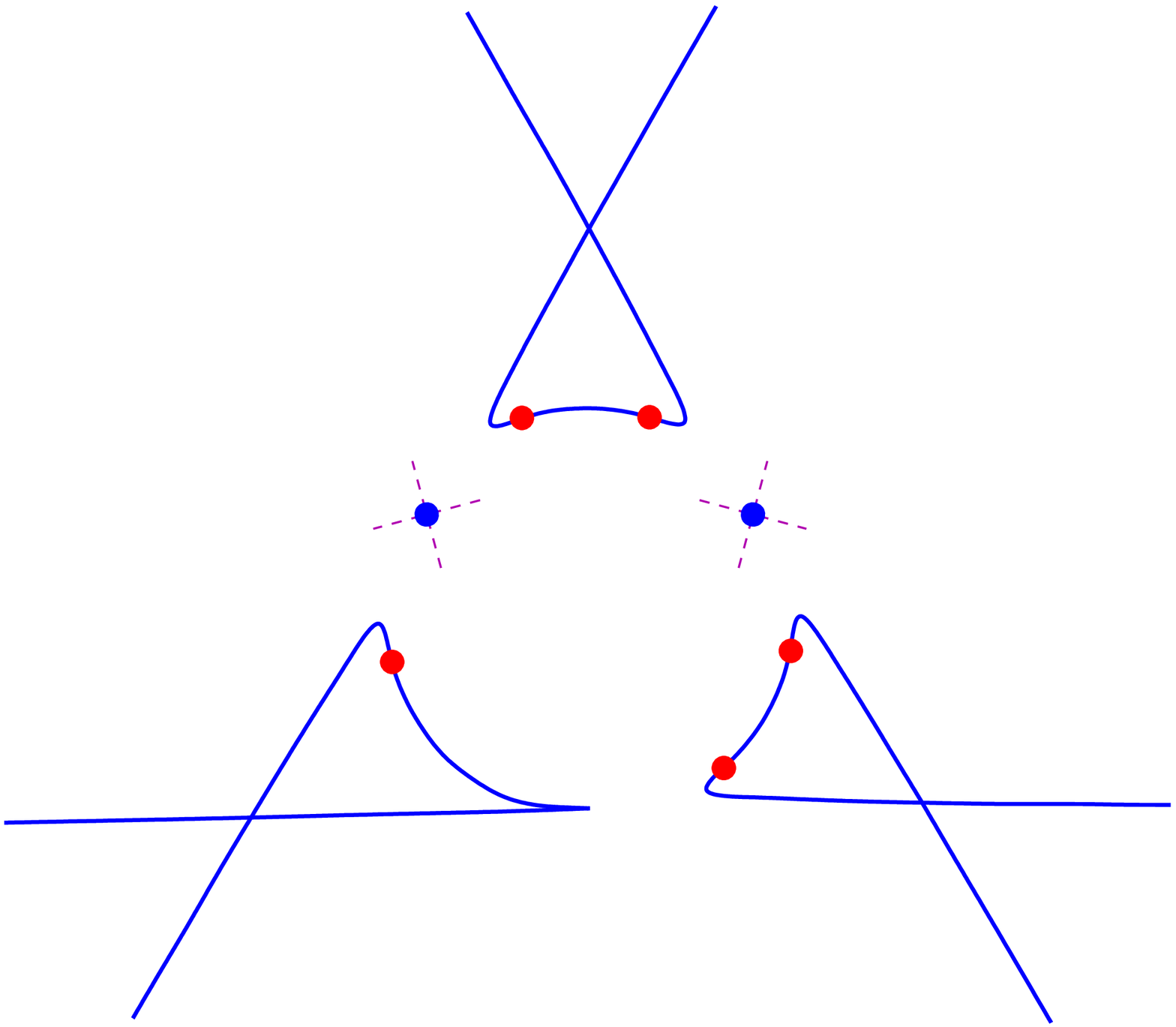}\quad
  \epsfysize=80pt\epsfbox{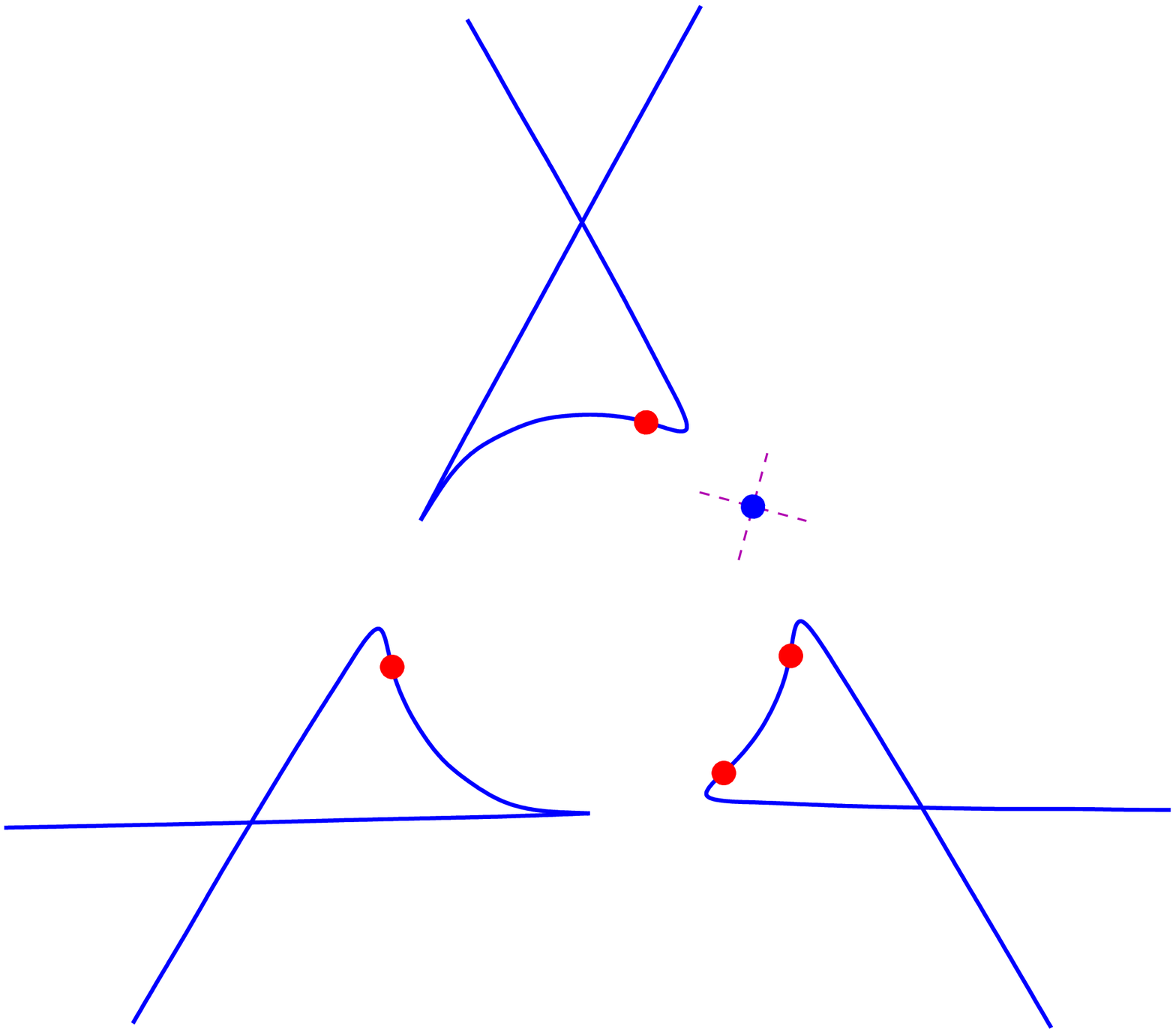}\quad
  \epsfysize=80pt\epsfbox{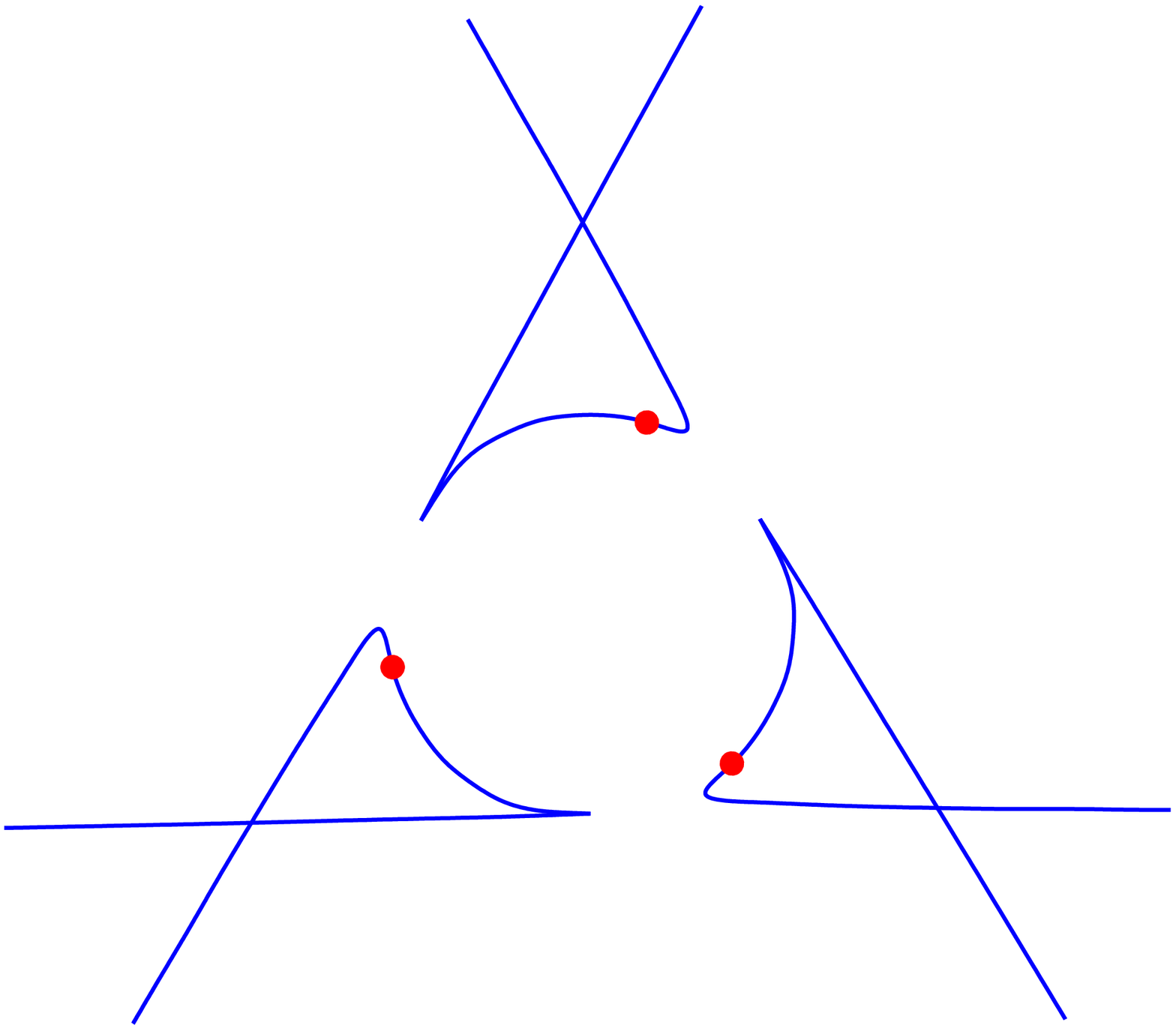}
$$
\caption{Maximally inflected degree 5 curves with 3 real nodes}
\label{fig:shustin}
\end{figure}

Figure~\ref{fig:shustin} shows these curves for $d=5$ and $\kappa=0,1,2,3$.

\begin{rem}\label{rem:alt-model}
For $t>0$, the cusp in the curve $K_t$~(\ref{eq:Lt}) is on the branch to the
left of the origin.
Had we instead used the perturbation $K'_t$ given by
$$
  y\,(y-x^2) \;-\; t\,x^3\ =\ 0\,,
$$
then the cusp is now on the right branch for $t>0$.
In Section 6.1, we use this to study Question~\ref{ques:two} concerning
possible necklaces for a given set of ramifications and numbers of solitary
points.
\end{rem}

\subsection{Patchworking of Singular Curves}\label{sec:patchworking}
Viro's method for constructing real plane curves with prescribed
topology~\cite{Vi83,Vi84} (see also~\cite{Ri}) begins with a
subdivision of the simplex
$$
  \Delta_d\ :=\ \{(i,j)\in({\mathbb Z}_+)^2 \mid i+j\leq d\}\,
$$
defined by a piecewise linear convex  {\it lifting function}.
Reflecting this subdivision in the coordinate axes and in the origin gives a
subdivision of the region
$$
  \lozenge_d\ :=\ \{(x,y)\in{\mathbb R}^2 \mid
                |x|+|y|\leq d\}\,.
$$
Gluing opposite edges of $\lozenge_d$ gives a
PL-space homeomorphic to $\mathbb{RP}^2$.

The other ingredient is, for each facet $F$ of the subdivision, a real
polynomial $f_F$ whose Newton polytope is $F$ and such that $f_F=0$ defines a
smooth curve in the torus $\left({\mathbb C}^\times\right)^2$.
These polynomials additionally satisfy a compatibility condition:
For an edge $e$ common to two facets $F$ and $G$,
we have $f_F|_e=f_G|_e$, where the latter expressions
denote
the truncations of the polynomials $f_F$ and $f_G$
to the edge $e$ (i.e.,
those
monomials whose
exponent vector
is
in $e$).
Furthermore, this common
truncation
has no multiple factors (except $x$ and $y$).

The real points of the real curve $f_F=0$ lie
naturally in the four copies of the facet
$F$ in $\lozenge_d$ (with the boundary points
representing the asymptotic behavior of solutions, or, equivalently,
the zeros of the restriction $f_F$ to the corresponding edge).
The pair consisting of these facets and curves
will be called the {\it Newton portrait}
of  $f_F=0$.
By the compatibility condition, the Newton
portraits of the facet curves glue
naturally together, giving a  topological curve $\Gamma$ in
$\mathbb{RP}^2$.
Viro's theorem asserts that there exists a curve $C$ of degree $d$
in $\mathbb{RP}^2$ such
that the pair $(\mathbb{RP}^2,C)$ is homeomorphic
to the pair $(\mathbb{RP}^2,\Gamma)$.
The homeomorphism
preserves each coordinate axis and each quadrant.
The complex points of $C$ are smooth, and $\Gamma$ and $C$ meet
each coordinate axis in the same number of points.\smallskip

Shustin~\cite{Sh85} (see also~\cite{Sh93,Sh98}) showed how to modify this
construction when the facet curves ($f_F=0$) have singularities in
$\left(\mathbb{C}^\times\right)^2$.
If a numerical criterion is satisfied
(such a criterion is given in Theorem 1.7 of~\cite{Sh93}), there exists a curve
of degree $d$ in $\mathbb{RP}^2$ whose singularities are the disjoint union of
the singularities of the facet curves, and whose topology is glued from that of
the facet curves as before.
(See Theorem~1.8 of~\cite{Sh93}, originally proven in~\cite{Sh85}.)
In particular, Shustin shows the following.

\begin{prop}[See ~\cite{Sh93} and ~\cite{Sh85}] \label{prop:Shustin}
 If the singularities of the facet curves are only nodal, then there exists a
 curve of degree $d$ in $\mathbb{RP}^2$ whose only (complex) singularities
 are the disjoint union of the singularities of the facet curves, and whose
 topology is given by gluing the facet curves as in the Viro construction.
\end{prop}

We use this to prove Statement~(2) of Theorem~\ref{thm:classical}.
(Another approach is indicated in a remark in the end
of Section.)\smallskip

\noindent{\bf Proof of Theorem~\ref{thm:classical}(2).}
Consider the subdivision of $\Delta_d$ given by the
piecewise linear convex lifting function $w$
which we define
for some the vertices of $\Delta_d$.
Set $ w(0,0) = w(d,0) = 0$, and

$$
   \left. \begin{array}{rcl}
     w(0,2i{+}2)& =& (2i+2)^2 \\
     w(d{-}1{-}2i,2i{+}1)& =& (2i+1)^2\rule{0pt}{13pt}
   \end{array} \right.
     \mbox{\quad for }i=0,\ldots,
  \left\lfloor\frac{d{-}1}{2}\right\rfloor\,.
$$
Here is the resulting subdivision of $\Delta_4$ and the
values of the lifting function.
$$
  \begin{picture}(133,76)(-40,0)
   \put(14,6){\epsfxsize=60pt \epsfbox{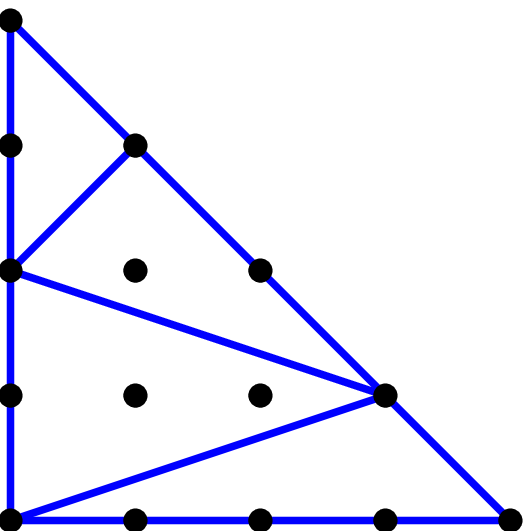}}
   \put(0,62){16}   \put(34,51){9}
   \put(6,32){4}    \put(62,21){1}
   \put(6, 0){0}    \put(76, 0){0}
   \put(100,15){\vector(-1,0){25}}\put(102,11){$H_4$}
   \put(85,36){\vector(-1,0){25}}\put(90,33){$Q_0$}
   \put(-40,19){$P_0$}\put(-25,22){\vector(1,0){25}}
   \put(-40,49){$P_1$}\put(-25,50){\vector(1,0){25}}
  \end{picture}
$$

This regular subdivision of $\Delta_d$ has three types of triangles
\begin{itemize}
 \item[(i)]
      The triangle $H_d$ with vertices $\{(0,0), (d,0), (d{-}1,1)\}$.
 \item[(ii)]
      The triangle $P_i$ with vertices $\{(0,2i),(0,2i+2),(d{-}1{-}2i,2i+1)\}$
      for each $i$ from 0 to $\lfloor\frac{d{-}2}{2}\rfloor$.
 \item[(iii)]
      The triangle $Q_i$ with vertices
      $\{(0,2i+2),(d{-}1{-}2i,2i+1),(d{-}3{-}2i,2i+3)\}$
      for each $i$ from 0 to $\lfloor\frac{d{-}3}{2}\rfloor$.
\end{itemize}

For each facet, we give polynomials $f_F$ that define curves with only
solitary points.
These will not necessarily satisfy the compatibility
condition,
but rather a weaker one:
A common edge $e$ between two facets $F$ and $G$ contains only two lattice
points, and after possibly multiplying one of the facet polynomials by $-1$,
the signs of the monomials from the two facet polynomials agree.
This weak compatibility implies
that there are monomial transformations
with positive coefficients of the
facet polynomials which do satisfy the compatibility criteria
after adjusting the sign of one of the two polynomials.
Since these monomial transformations do not change the geometry
(number of solitary points, topology of the glued curve $\Gamma$),
and the dual graph to the triangulation is a chain,
it will suffice to construct polynomials satisfying the weaker criteria
and giving the desired topology.

We describe the monomial transformations.
Consider first a common edge $e$ between adjacent facets $F$ and $G$ of the
triangulation.
Since $e$ has no interior lattice points, the restrictions of the facet
polynomials to $e$ will be binomials of the form
$$
  f_F|_e\ =\ A x^ay^b + B x^cy^{b+1}\qquad\qquad
  f_G|_e\ =\ C x^ay^b + D x^cy^{b+1}\,.
$$
(Not only do the exponents of $y$ differ by 1, but one of 
the exponents
$a$ or $c$ is zero.)
For $z\neq 0$, let $\sgn(z):=z/|z|$.
We compute\medskip

\quad $\sgn\left(\frac{A}{C}\right)
         \left|\frac{C}{A}\right|^{\frac{-c}{a-c}}
         \left|\frac{D}{B}\right|^{-b}
         f_F\left( \left|\frac{C}{A}\right|^{\frac{1}{a-c}}\,x,\
         \left|\frac{D}{B}\right|\,y\right)$\vspace{-5pt}
\begin{eqnarray*}
\qquad\qquad&=& {\textstyle \sgn\left(\frac{A}{C}\right)
         \left|\frac{C}{A}\right|^{\frac{-c}{a-c}}
         \left|\frac{D}{B}\right|^{-b}
    \left( A\left|\frac{C}{A}\right|^{\frac{a}{a-c}}
           \left|\frac{D}{B}\right|^{b}x^ay^b +
          B \left|\frac{C}{A}\right|^{\frac{c}{a-c}}
            \left|\frac{D}{B}\right|^{b+1}x^cy^{b+1}\right)}\\
&=& {\textstyle \sgn\left(\frac{A}{C}\right)
    \left( A\left|\frac{C}{A}\right|x^ay^b +
          B \left|\frac{D}{B}\right|x^cy^{b+1}\right)}\\
&=& C x^ay^b + D x^cy^{b+1}\,,
\end{eqnarray*}
as the weak compatibility criteria ensures that
$\sgn(\frac{A}{C})=\sgn(\frac{D}{B})$.

Since the dual graph of the triangulation is a chain,
we encounter no obstructions when
transforming the facet polynomials so that they
satisfy the compatibility condition.
More precisely, given the facet polynomial for $H_d$, we transform the facet
polynomial for $P_0$ as above, then the facet polynomial for $Q_0$, then for
$P_1$, and {\it et cetera}.\medskip

We now give the facet polynomials.
The reader is invited to check that the weak compatibility conditions
are satisfied.
The facet $H_d$ is the Newton polytope of the polynomial
$h_d:=x^{d-1}y - (1-x)(2-x)\cdots(d-x)$.
Here are the Newton portraits of $h_3$ and $h_4$.
 \[
  \epsfysize=.6in \epsfbox{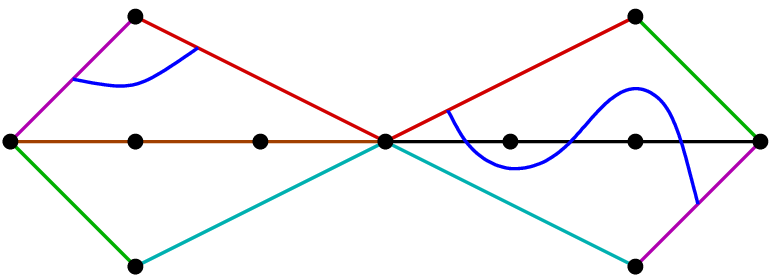}\qquad
  \epsfysize=.6in \epsfbox{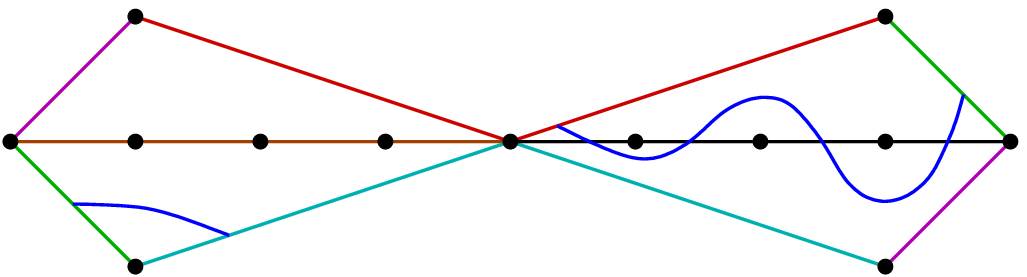}
 \]

The remaining facet polynomials are based on an idea of
Shustin~\cite[p.~849]{Sh93}.
Recall that the Chebyshev polynomials $\Tc_k(x)$ (defined recursively by
$$
  \Tc_0 := 1, \quad \Tc_1 := x, \quad\mbox{and for $k>1$,}\quad
  \Tc_k(x) := 2x\Tc_{k-1}(x) - \Tc_{k-2}(x)\,.)
$$
have the property that $\Tc_k(x)$ has exactly $k$ roots in the interval
$(-1,1)$ and $k{-}1$ local extrema in this interval with values $\pm 1$, and
for $|x|>1$, $|\Tc_k(x)|>1$.
Lastly, the leading term of $\Tc_k(x)$ is $2^{k-1}x^k$.
Then the polynomial
$f_k(x):=y^2-2y\Tc_k(x+2)+1$ has Newton polytope the triangle
 \[
  \mbox{Conv}\{(0,0), \ (1,k), \ (0, 2)\}\,,
 \]
which is a translate of the polytope $P_i$ by $(0,-2i)$ when $d-1-2i=k$.
The curve $f_k=0$ in ${\mathbb R}^2$ has 2 connected components and
$k-1$ solitary points.
We display these curves for $k=2, 3$, and $4$, scaling the $y$-axis by the
transformation $y\mapsto \mbox{sign}(y)|y|^{1/k}$.
$$
{\setlength{\unitlength}{0.8pt}
  \begin{picture}(118,95)(-8,-2)
   \put( 0, 0){\epsfxsize=90pt \epsfbox{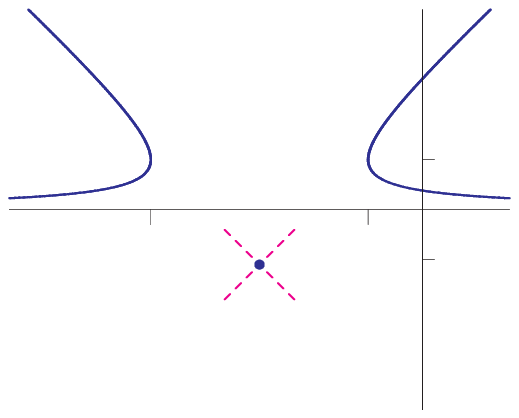}}
   \put( 33, 1){$f_2=0$}   \put(-3,32){$x$}    \put(76,82){$y$}
  \end{picture}
\qquad
  \begin{picture}(118,93)(-8,0)
   \put( 0, 0){\epsfxsize=90pt \epsfbox{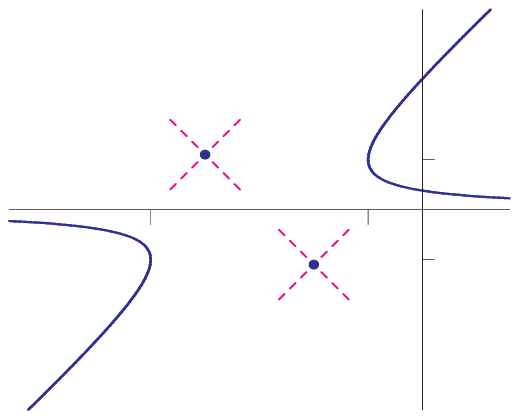}}
    \put( 33, 1){$f_3=0$}  \put(-3,32){$x$}    \put(76,82){$y$}
  \end{picture}
\qquad
  \begin{picture}(118,93)(-8,0)
   \put( 0, 0){\epsfxsize=90pt \epsfbox{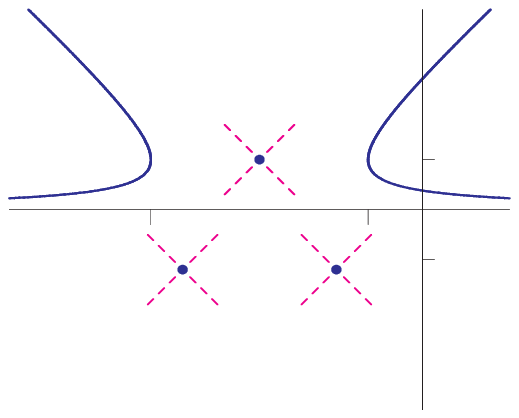}}
   \put( 33, 1){$f_4=0$}   \put(-3,32){$x$}    \put(76,82){$y$}
  \end{picture}}
$$
Here are the Newton
portraits of these curves.
We omit the interior lattice
points in the triangles.
$$
\epsfysize=63pt \epsfbox{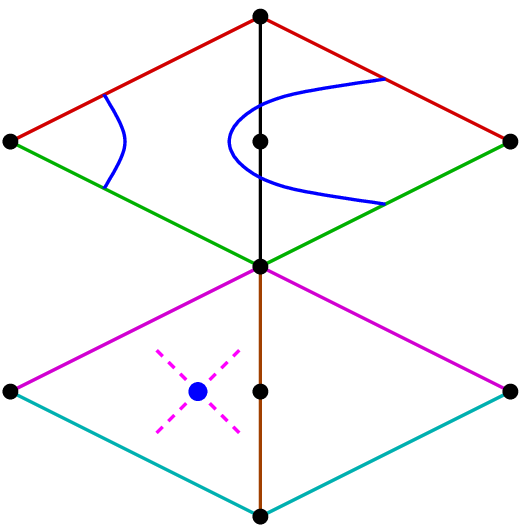}\qquad
\epsfysize=63pt \epsfbox{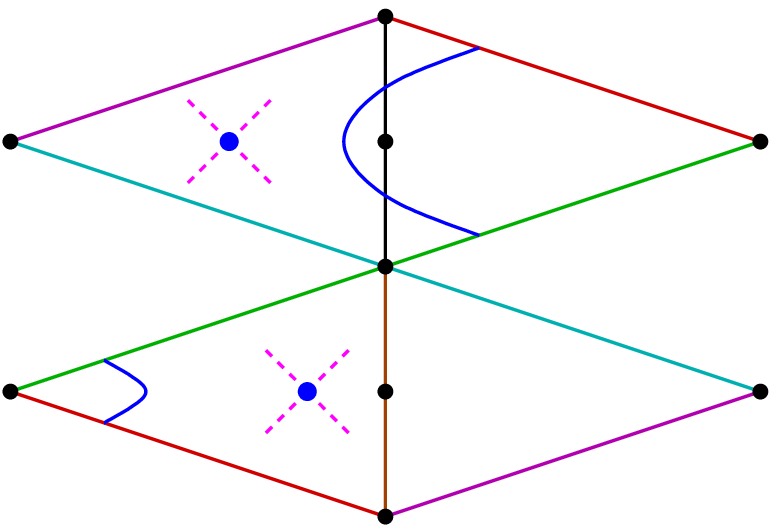}\qquad
\epsfysize=63pt \epsfbox{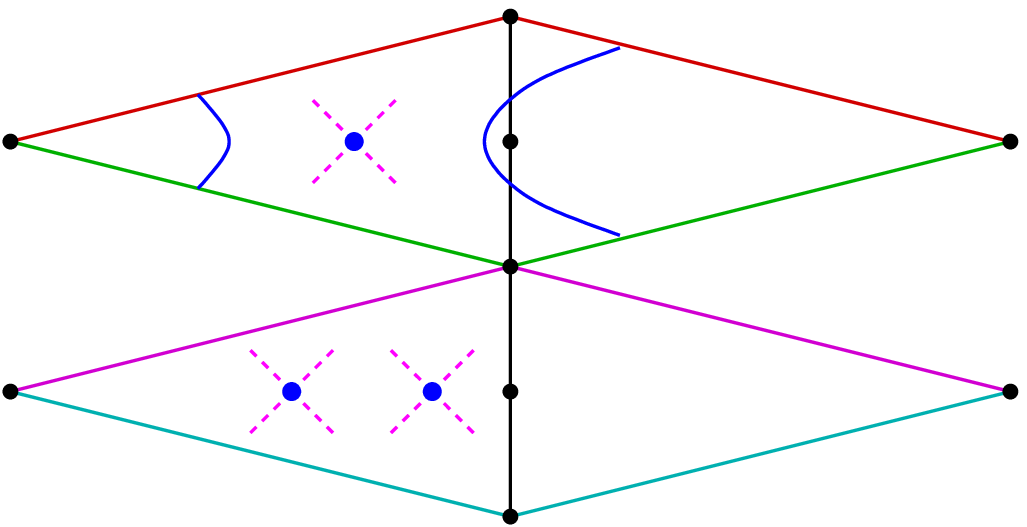}
$$
Let $y^{2i}f_{d-1-2i}$ be the facet polynomial for the facet $P_i$.

Finally, set
$$
  g_k(x,y) := f_k\left(-\frac{1}{x},\,(-1)^k\frac{y}{x}\right)\ = \
       \frac{y^2}{x^2} - (-1)^k 2\frac{y}{x}
                 \Tc_k\left(2-\frac{1}{x}\right) + 1\,.
$$
The Newton polytope of $g_k$ has vertices $\{(0,0),  (-2,2), \ (-k-1, 1)\}$.
We display the Newton portraits
of $g_1$, $g_2$, and $g_3$.
(For this, we first translate their Newton polytopes by $(k+1,0)$, placing it
into the positive quadrant.)
$$
\epsfysize=70pt \epsfbox{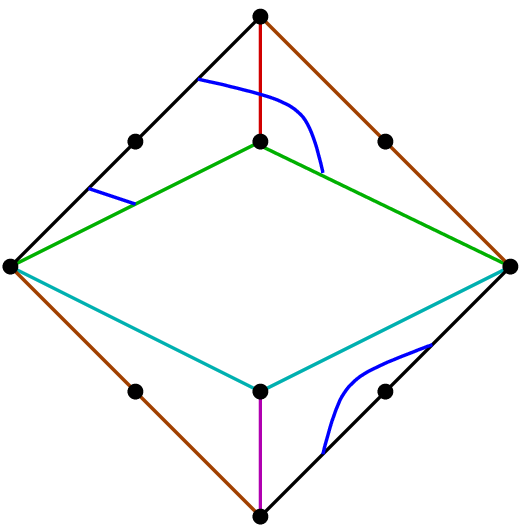}\qquad
\epsfysize=70pt \epsfbox{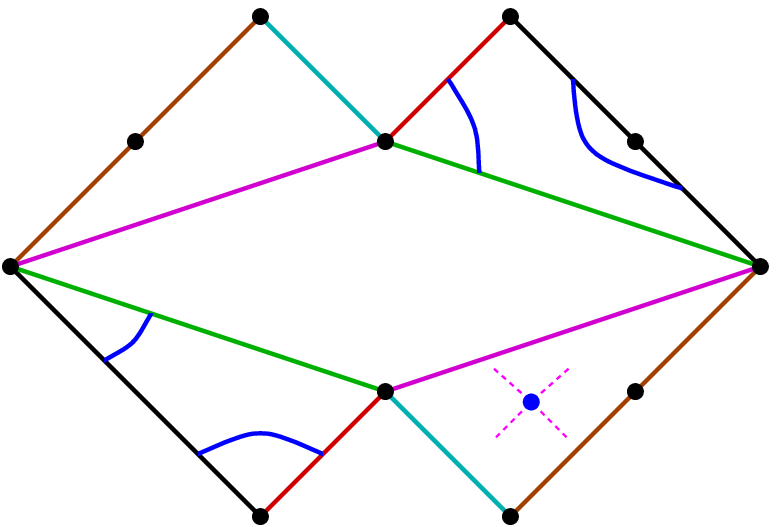}\qquad
\epsfysize=70pt \epsfbox{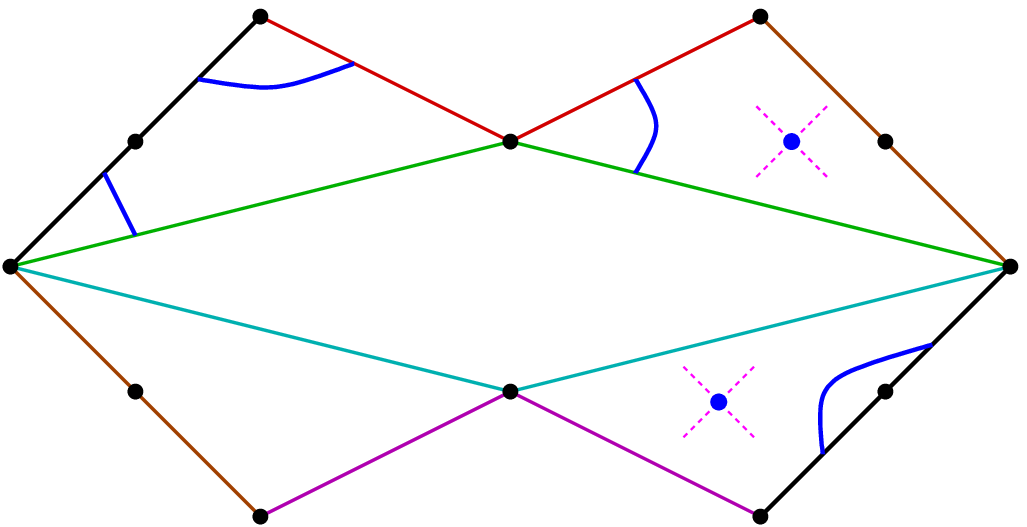}
$$
Let $x^{d-1-2i}y^{1+2i}g_{d-2-2i}$ be the facet polynomial of the facet
$Q_i$.\medskip

The curve $C_d$ constructed from these data by Proposition~\ref{prop:Shustin}
has $3(d{-}2)$ flexes and $\binom{d-1}{2}$ solitary points as claimed.
First, since the facet curves $f_k$ and $g_k$ each have $k{-}1$
solitary points,
$C_d$ has $(d{-}2)+(d{-}3)+\cdots+1+0 = \binom{d-1}{2}$ solitary points, and
so is rational (in fact, the solitary points
correspond to internal integer points of the Newton polygon).
{}From the Newton portrait
of $h_d$, we see that $C_d$ meets the $x$-axis in $d$
points.
Each facet $P_k$ contributes 2 points of intersection of $C_d$ with the $y$-axis.
When $d$ is even, this gives $d$ points of intersection, and when $d$ is odd,
$d{-}1$ points of intersection.
When $d$ is odd, the last facet $Q_1$
(corresponding to $g_1$)
contributes an additional point of
intersection with the $y$-axis.
Finally, each facet $Q_k$ contributes 2 points of intersection of $C_d$ with the
$z$-axis, giving $d{-}2$ points of intersection when $d$ is even and $d{-}1$ points
of intersection when $d$ is odd.
The facet $H_d$ contributes an additional point, and when $d$ is even, the
facet $P_{\frac{d-2}{2}}$
(corresponding to $f_2$)
contributes one point.

As a result, the curve $C_d$
has three separate segments, each intersecting
a different coordinate axis in $d$ consecutive points going in the same order on
the curve and on the axis.
Thus, by counting the Whitney indices by means of the
Gauss map we find at least $d{-}2$ flexes on each segment.
Hence, the curve constructed has $3(d{-}2)$ flexes.
\QED

\begin{figure}[htb]
$$
   \epsfysize=140pt \epsfbox{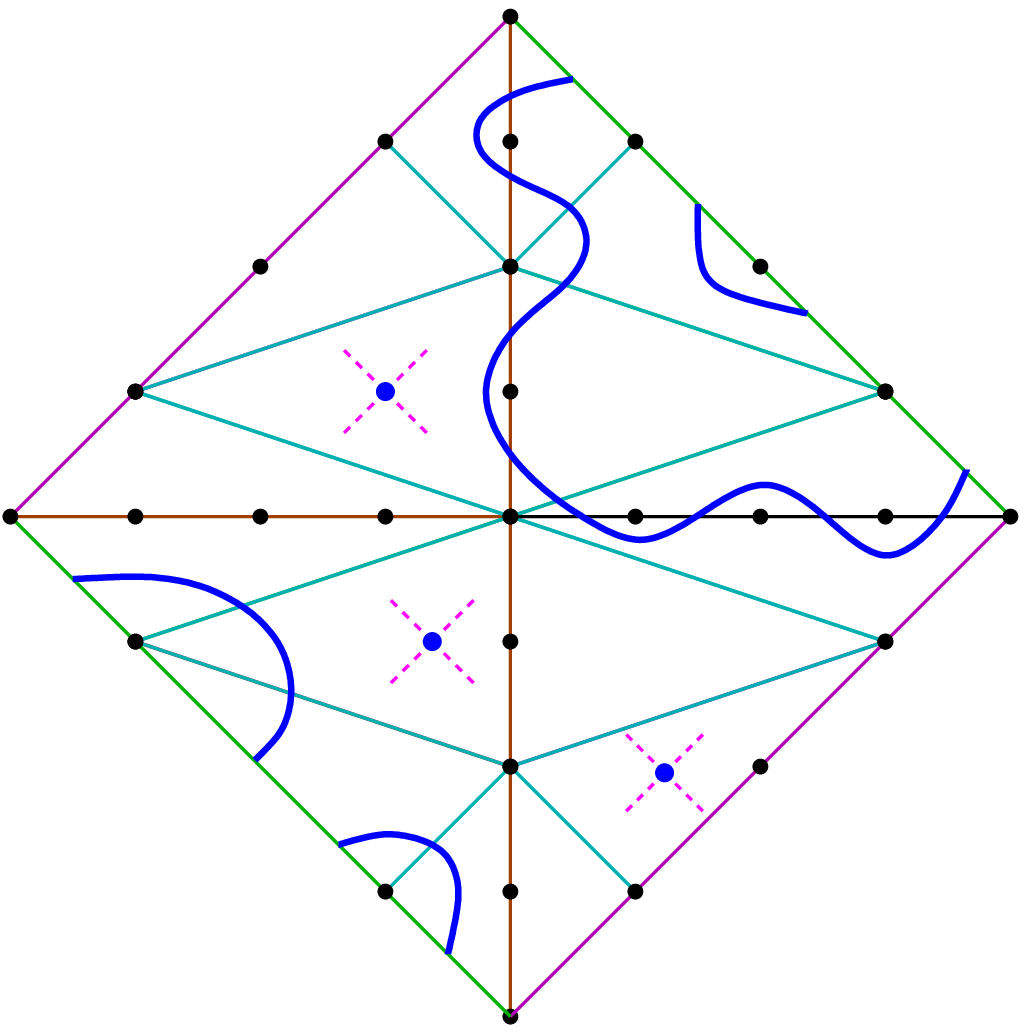} \qquad
   \epsfysize=140pt \epsfbox{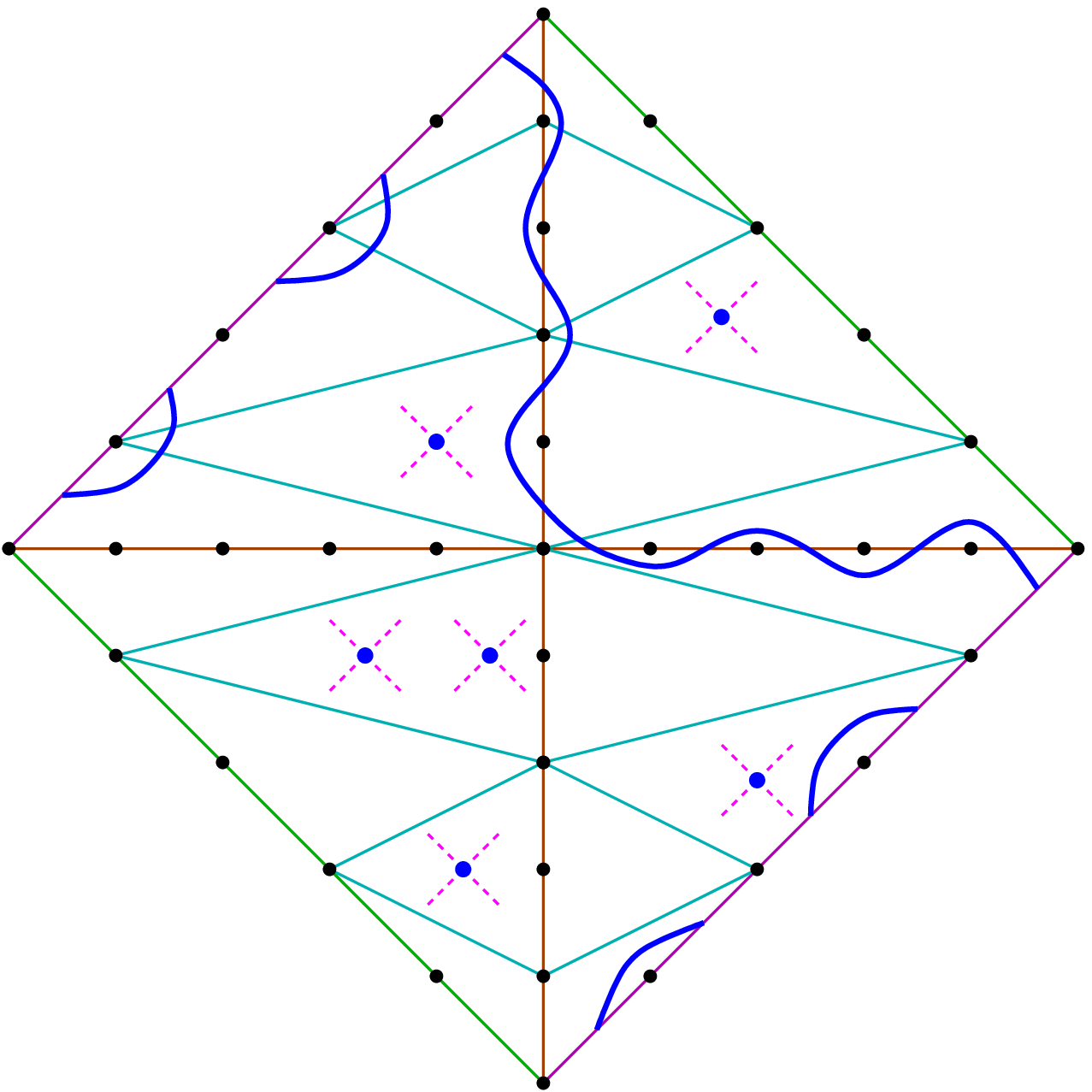}
$$
\caption{The curves $C_4$ and $C_5$.}
\label{fig:harnack}
\end{figure}

Figure~\ref{fig:harnack} shows the curves $C_4$ and $C_5$.
The curve $C_4$ has the topological type of the quartic in
Table~\ref{table:two} with six flexes and no real nodes.

\begin{rem}
Another patchworking is via gluing parameterizations of the
facet curves $f$ and $g$ (which are rational).
Since the dual graph of the triangulation is a chain,
it is sufficient to have a gluing construction for
a pair of rational plane curves intersecting transversally.
For that purpose, one can pick parameterizations $F$ and $G$ such
that $F(0)=G(0)$ and consider, for generic $\lambda, \mu$
and sufficiently small $\epsilon>0$ the rational curves $H_\epsilon$
given by $\lambda F(u)+\mu G(v), uv=\epsilon$.
The flexes and the nodes are preserved under this patchworking construction,
since they are stable under small deformations.
\end{rem}

\section{Maximally inflected plane curves of degrees three and four}
\label{sec:lowdegree}

\subsection{Cubic Plane Curves}
In Figure~\ref{fig:cubics} we saw that a plane cubic with 3 real flexes
has a solitary point and no nodes and a plane cubic with one real cusp has no
nodes.
These are the only possible maximally inflected cubics and
they
exhaust all the possibilities allowed by Corollary~\ref{cor:bounds}.

\subsection{Quartic Plane Curves}
Consider now quartics whose only ramifications are flexes and cusps.
The upper bound for $\eta+2c$ allowed by Corollary~\ref{cor:bounds} for quartics
is 1, so maximally inflected quartics have no complex nodes and either 1 or 0 real
nodes.
Table~\ref{table:one} summarizes the possible numbers $\eta$ of real nodes
and $\delta$ of solitary points.
\begin{table}[htb]
\begin{tabular}{|c||c|c|c|c|c|c|c|}
\hline
$\kappa$&\multicolumn{2}{c|}{0}&\multicolumn{2}{c|}{1}&
                   \multicolumn{2}{c|}{2}&3\\\hline
$\iota$&\multicolumn{2}{c|}{6}&\multicolumn{2}{c|}{4}&
                   \multicolumn{2}{c|}{2}&0\\\hline
$\eta$&0&1&0&1&0&1&0\\\hline
$\delta$&3&2&2&1&1&0&0\\\hline
\end{tabular}\medskip
\caption{Topological invariants of quartics allowed by Corollary~\ref{cor:bounds}.\label{table:one}}
\end{table}
Clearly, for quartics, these
numbers determine the real part of the image
up to homeomorphism. In fact, for quartics, they determine it up
to isotopy in $\mathbb{RP}^2$
(see Theorem \ref{thm:quartic-isotopy}, Remark \ref{on-isotopy},
and Theorem \ref{T:quart-isotopy}).

\begin{thm}\label{thm:4-class}
For every quadruple $(\kappa,\iota,\eta,\delta)$ given in
Table~\ref{table:one}, there
is a (real) plane quartic with $\kappa$ cusps, $\iota$ inflection points,
$\eta$ nodes, and $\delta$ solitary points.
\end{thm}

\noindent{\bf Proof. }
This may be deduced the classification of real rational quartics given by
D.A. Gudkov, et. al.~\cite{GUT66} (see also F. Klein~\cite{Klein1876b},
H.G.~Zeuthen~\cite{Ze1874}, and C.T.C.~Wall~\cite{Wa95}).
However, we will instead exhibit examples of curves with each possible
ramification and singularity.
In Table~\ref{table:two}, we display a maximally inflected plane quartic
curve for each quadruple $(\kappa,\iota,\eta,\delta)$
of Table~\ref{table:one}.
\begin{table}[htb]
\begin{tabular}{|c|c||c|c|c|c|} \cline{3-6}
\multicolumn{2}{c||}{ }&\multicolumn{4}{c|}{\large$\kappa,\iota$}\\ \cline{3-6}
\multicolumn{2}{c||}{ }& 0,6 & 1,4 & 2,2 & 3,0\\\hline\hline
\multirow{2}{15pt}{\large$\eta$}&
\raisebox{30pt}{0}
 &\epsfysize=70pt \epsfbox{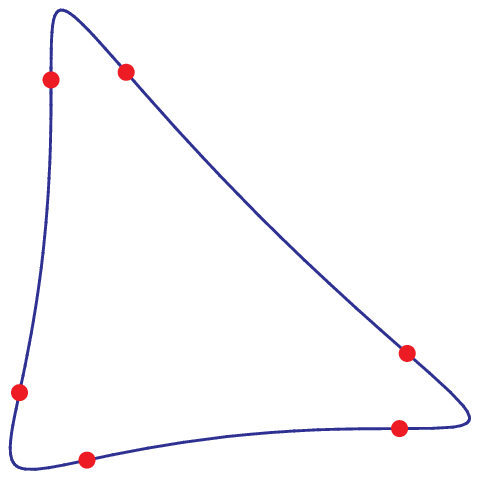}
 &\epsfysize=70pt \epsfbox{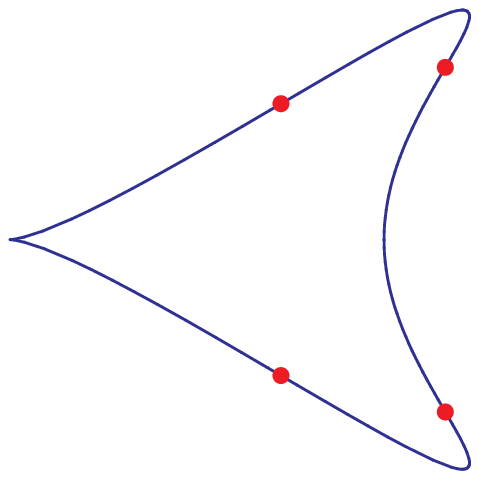}
 &\epsfysize=70pt \epsfbox{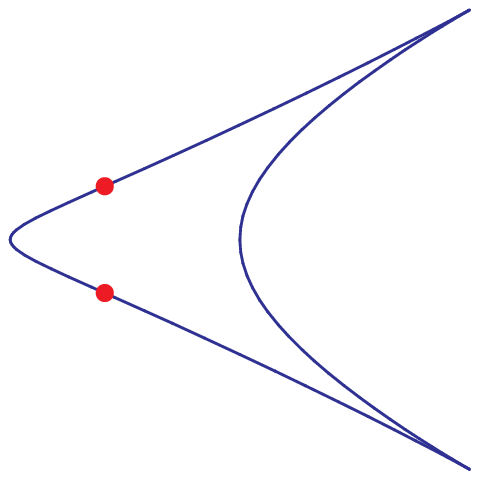}
 &\epsfysize=70pt \epsfbox{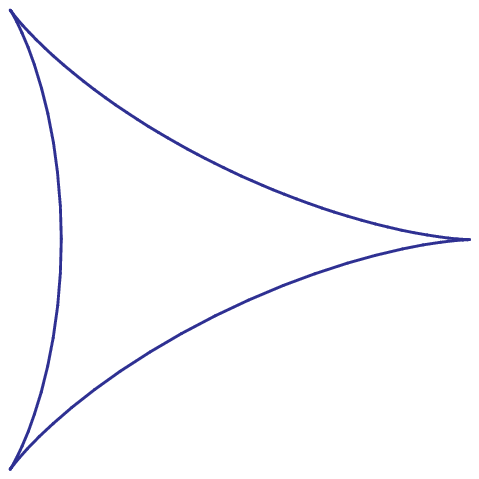}\\\cline{2-6}
&\raisebox{30pt}{1}
 &\epsfysize=70pt \epsfbox{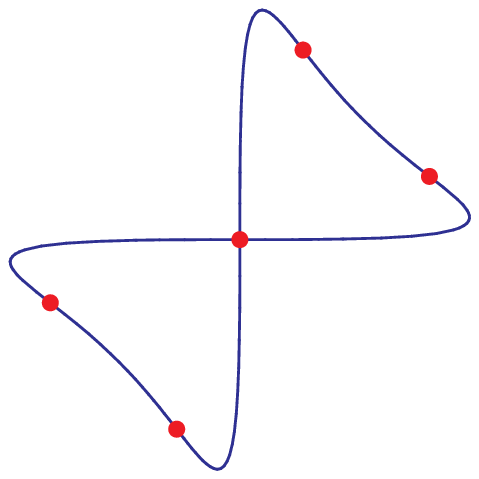}
 &\epsfysize=70pt \epsfbox{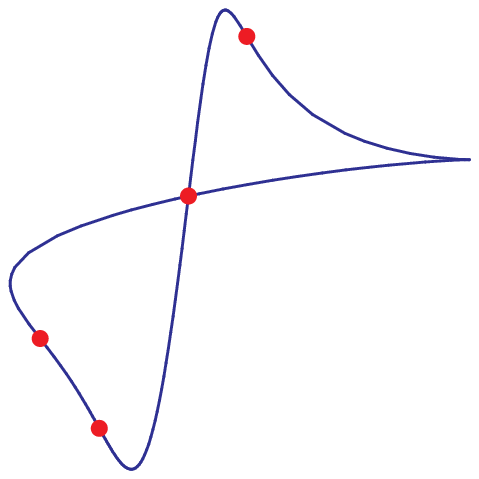}
 &\epsfysize=70pt \epsfbox{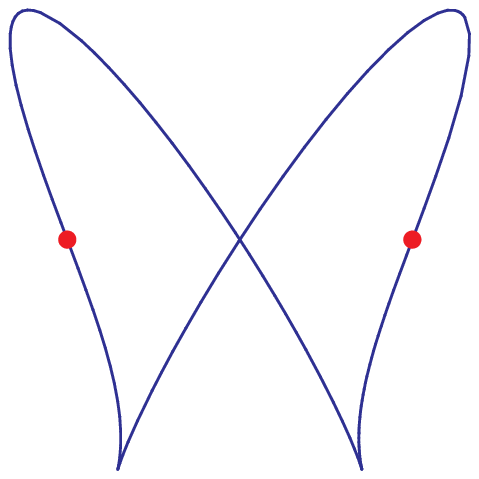}
 &\raisebox{30pt}{
   $\begin{array}{c}\textrm{not}\\ \textrm{allowed}\end{array}$}\ \  
\\\hline
\end{tabular}\medskip
\caption{Quartics realizing all topological types allowed by
Corollary~\ref{cor:generalbounds}.\label{table:two}}
\end{table}
These were generated using a computer calculation, by first solving for the
centers of projection (as described, for example in~\cite[Section 2]{So00b} or
in~\cite[Section 2]{So_shapiro-www}), and then drawing the resulting
parameterized curve using MAPLE.
(This method is used to draw most of the curves we display.)
The positions of the flexes are marked with dots and the curve with 6 flexes
and one node has 2 flexes at
its node.\QED

Recall from Section~\ref{sec:MIC} that Conjecture~\ref{conj:non-degen}
holds for plane quartics and so we have the strong information of
Theorem~\ref{thm:Implies} about deformations of plane quartics.
We explore some consequences of that fact.
There is a single isotopy class of sextuples of distinct points
on $\mathbb{RP}^1$.
Thus, given any
sextuple $S=\{s_1,\ldots,s_6\}$
of distinct (non ordered) points, there is
an isotopy from the sextuple $\{-3,-1,0,1,3,\infty\}$ to $S$.

\begin{figure}[htb]
\[
   \epsfbox{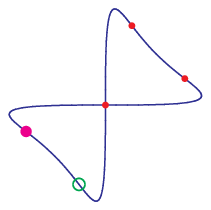}\qquad
   \epsfbox{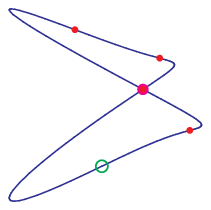}\qquad
   \epsfbox{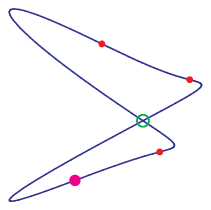}\qquad
   \epsfbox{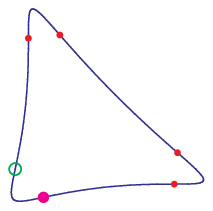}\qquad
   \epsfbox{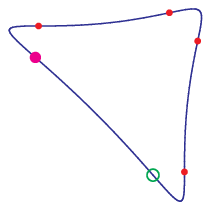}
\]
   \caption{The 5 curves with flexes at
            $\{-3,-1,0,1,3,\infty\}$\label{fig:6flex}}
 \end{figure}

The Schubert calculus gives 5 rational quartics with 6 given points
of inflection, thus each of the 5 maximally inflected quartics with flexes
at $S$ are deformations of one of the 5 maximally inflected curves with flexes at
$\{-3,-1,0,1,3,\infty\}$, which we display in Figure~\ref{fig:6flex}.
(Each nodal curve has 2 flexes at its node.)
We indicate the differences in the parameterizations of these curves,
labeling the flex at $-3$ by the larger dot and the flex at $-1$ by the
circle.

\begin{figure}[htb]
   \epsfysize=130pt  \epsfbox{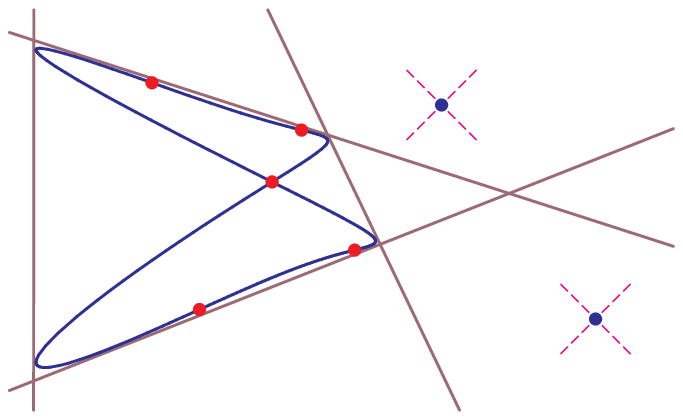}\qquad \qquad
   \epsfysize=130pt\epsfbox{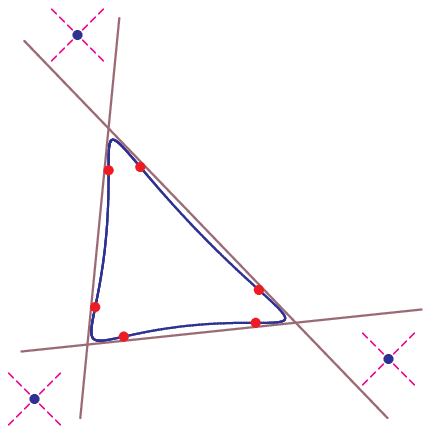}
   \caption{Bitangents and solitary points.\label{fig:bitangents}}
 \end{figure}

\begin{thm}\label{thm:quartic-isotopy}
 For each sextuple $S=\{s_1,\ldots,s_6\}$ of distinct points in $\mathbb{RP}^1$,
 there are exactly five maximally inflected quartics with flexes at $S$.
 Of these five, two have three solitary points and no real nodes, while
 three have two solitary points and one real node.
 Moreover, the possible arrangements of solitary points and bitangents are as
 indicated in Figure~\ref{fig:bitangents} for each of these types of
 curves.
\end{thm}

\noindent{\bf Proof. }
Since the five quartics in Figure~\ref{fig:6flex} show that the statement is valid
for the choice of flexes at $\{-3,-1,0,1,3,\infty\}$,
it suffices to show
that the number of solitary points does not change under a deformation of such
a quartic curve.

Any deformation of a quartic with one real node must have one real node; in
passing to a curve without a real node, the deformation would include a curve
with some other singularity, which would necessarily be a ramification point
that is not a flex.
Since the curve already contains six flexes, this would violate
Proposition~\ref{prop:ram}.
Thus every deformation of the first three curves has a single real node.

To see that every deformation of
any of the five curves has only ordinary double points,
consider the arrangements of bitangents and solitary points
on an example.
Figure~\ref{fig:bitangents} shows the second and fourth curves with their
bitangents and solitary points.
Here, each solitary point is separated
from other real, solitary or not, points of the curve by some number
of real double tangents.
The same phenomena takes place for
each of the five maximally inflected curves.

Since the curves have degree 4, the bitangents cannot meet the curves in
additional points, by B\'ezout's theorem.
Similarly, a bitangent cannot be tangent to a flex.
Thus in a deformation, the flexes are confined to the arcs of the curve
between two points of tangency to bitangents, and the number of bitangents
does not change.
Similarly, the solitary points cannot meet
 another point of the curve, to do so,
they would first have to meet a bitangent, which is not possible,  by
B\'ezout's theorem.
This completes the proof of the theorem.
\QED

\begin{rem}
 For reference, in Figure~\ref{F:Param},
 we give the solitary points and bitangents, as well as
 parameterizations $[\varphi_1(s,t),\varphi_2(s,t),\varphi_3(s,t)]$ for the
 two curves in Figure~\ref{fig:bitangents}.
 The decimals are numerical approximations.
\begin{figure}[htb]
\begin{center}
 \begin{minipage}[t]{220pt}
   \hspace{1cm}{\sc Nodal Curve}\smallskip

   \hspace{.1cm}{Parameterization}
  \begin{eqnarray*}
   \varphi_1(s,t)&=&{\textstyle s^4-6t^2s^2+9t^4}\\
   \varphi_2(s,t)&=&{\textstyle \frac{3}{4}s^4+ts^3 -\frac{3}{2}t^2s^2
                    + \frac{9}{4}t^4}\\
   \varphi_3(s,t)&=&{\textstyle s^4+ts^3-3t^2s^2-2t^3s+\frac{15}{2}t^4}
  \end{eqnarray*}
   \hspace{.1cm}{Solitary Points}
  \begin{eqnarray*}
    &(1.514769,\ 0.854076)&\qquad\\
    &(2.088892,\ 0.040735)&\qquad
  \end{eqnarray*}
   \hspace{.1cm}{Bitangents}
  \begin{eqnarray*}
    x&=&0\\
    y&=&{\textstyle 3 -\frac{25}{12}\,x}\\
    y&=&1.07221014   - 0.31889744\,x\\
    y&=&-.17859312   + 0.39336553\,x
  \end{eqnarray*}
 \end{minipage}
\qquad
 \begin{minipage}[t]{180pt}
   \hspace{1cm}{\sc Anodal Curve}\smallskip

   \hspace{.1cm}{Parameterization}
  \begin{eqnarray*}
   \varphi_1(s,t)&=&{\textstyle ts^3-\sqrt{6}t^2s^2 -\frac{3\sqrt{6}}{2}t^4}\\
   \varphi_2(s,t)&=&{\textstyle t^3s+\frac{\sqrt{6}}{2}t^4}\\
   \varphi_3(s,t)&=&{\textstyle s^4+6t^2s^2+9t^4}
  \end{eqnarray*}
   \hspace{.1cm}{Solitary Points}
  \begin{eqnarray*}
    {\textstyle \bigl(\frac{7}{16}-\frac{\sqrt{6}}{32},}&&\hspace{-12pt}
      {\textstyle -\frac{7}{48}+\frac{13\sqrt{6}}{288}\bigr)}\\
    {\textstyle \bigl(-\frac{7}{16}-\frac{\sqrt{6}}{32},}&&\hspace{-2pt}
      {\textstyle \frac{7}{48}+\frac{13\sqrt{6}}{288}\bigr)}\\
    {\textstyle \bigl(-\frac{\sqrt{6}}{4},}&&\hspace{-12pt}
      {\textstyle -\frac{\sqrt{6}}{36}\bigr)}
  \end{eqnarray*}
   \hspace{.1cm}{Bitangents}
  \begin{eqnarray*}
    y&=&{\textstyle 0.03402069-\frac{1}{3}\,x}\\
    y&=&1.4984552+3.2996598\,x\\
    y&=&-.0015448+   .0336735\,x
  \end{eqnarray*}
 \end{minipage}
\end{center}
\caption{Parameterizations, solitary points, and bitangents for the
         curves of Figure~3\label{F:Param}}
\end{figure}

 Interestingly, the two solitary points of the nodal curve lie on the line
 $y=3-\frac{17}{12}\,x$, and this line also meets the two
 points of intersection of the bitangents to the right of and above the
 quartic.
 Similarly, the node lies at the intersection of the pairs of lines through the
 other points of intersection of the bitangents.
 The other two nodal quartics with flexes at $\{-3,-1,0,1,3,\infty\}$ shown in
 Figure~\ref{fig:6flex} also have this property.
 That is clear for the second asymmetric nodal quartic whose image in $\RP^2$
 is isomorphic to this quartic.
 The symmetric nodal quartic has its solitary points at $[1,0,0]$ and
 $[0,1,0]$, and the corresponding pairs of bitangents are parallel, so all four
 points lie on the line at infinity.
 The statement about its node is also clear from symmetry.

 By the formula of Theorem~\ref{thm:bounds} for such quartics,
 $\delta=\tau+2$, where $\delta,\tau$ count the solitary points and solitary
 bitangents.
 Thus the nodal curve has no solitary bitangents, while the other curve has
 a single solitary bitangent, which is the line at infinity.
 To see this, note that $\varphi_3$ factors as
 $(s^2+3t^2)^2$, and so the line at infinity is a bitangent, with tangencies
 the points of the curve where $[s,t]=[\pm\sqrt{-3},1]$.
 Evaluating, we see these are at $[-1\pm2\sqrt{-2},1,0]$.
\end{rem}

\begin{rem}\label{on-isotopy}
 Theorem~\ref{thm:quartic-isotopy} shows that the isotopy type of the
 embedding  of a quartic with 6 flexes
 into $\mathbb{RP}^2$ does not change under an isotopy of the positions
 of the flexes in $\mathbb{RP}^1$.
 In fact something stronger is true.
 The space of all possible positions of flexes
 modulo reparameterizations preserving orientation (i.e., the quotient of
 $(\mathbb{RP}^1)^6$
 minus all the diagonals by $\mathit{SL}(2,\mathbb{R})$)
 is contractible  (to see this, consider fixing the position of one flex),
 and because of the confinement property of flexes and the
 tangent points of the bitangents (see the proof of
 Theorem~\ref{thm:quartic-isotopy}),
 there is no way to deform
 one of the five quartics in Figure~\ref{fig:6flex} into any
 other.

 However, one may allow two flexes to come together, for example, the flexes
 in Figure~\ref{fig:6flex} represented by the thickened dot and the open
 circle.
 When the positions of these two flexes collide, the second, third, and fourth
 curve will develop a cusp, while the first and fifth will develop a
 {\it planar point} with ramification sequence (0,1,4).
 Suppose that we fix
 five of the ramification points, and let the sixth
 move along $\mathbb{RP}^1\simeq S^1$.
 At every position of the sixth point, we get a maximally inflected curve
 which has 6 flexes, except when the 6th point collides with one of the
 fixed points, and then we get a curve with either a cusp or a planar
 point\frankfootnote{Motion pictures of families of curves with such moving
 ramification points may
 be found at {\tt www.math.umass.edu/\~{}sottile/pages/inflected}.}.

 We have used symbolic methods to calculate what happens when the sixth point
 moves.
 The number of nodes is always preserved, and when the sixth point returns to
 its original position, we get a curve different than the original one.
 In fact this monodromy action cyclically permutes the three nodal quartics in
 Figure~\ref{fig:6flex} and interchanges the two quartics without nodes.
 This can also be inferred by visualizing the effect of moving the the flexes
 labeled by the open circles in Figure~\ref{fig:6flex}.
\end{rem}

\begin{rem}\label{on-Zeuthen}
The proof we gave of Theorem~\ref{thm:quartic-isotopy}
is based on the following property going
back to Zeuthen~\cite{Ze1874}: each
of the components of the complement in $\RP^2$
of the double and solitary tangents contains at most $1$ connected
compact component of the real locus. It is valid for any quartic, and
a similar proof shows that every maximally inflected
rational quartic has only ordinary double
points, apart from the singularities at the ramification points, and thus
we obtain the analog of Theorem~\ref{thm:quartic-isotopy} for all maximally
inflected rational quartics, which we state below after we introduce a
further, necessary notion.
\end{rem}
\medskip

Underlying a deformation of a maximally inflected curve is
the isotopy type of its ramification.
Consider the curves with 2 flexes and 2 cusps in
Table~\ref{table:two}.
In the curve with no real nodes, $\eta=0$, the cusps ($\kappa$) and
flexes ($\iota$) occur along $\mathbb{RP}^1$ in the order
$\kappa\kappa\iota\iota$,  while for the curve with a
single real node, the order is $\kappa\iota\kappa\iota$.
There is another curve with 2 cusps, 2 flexes, and one node.
$$
   \epsfysize=90pt \epsfbox{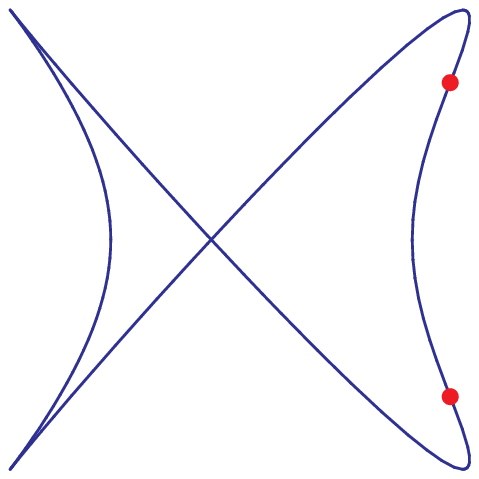}
$$
Here, the ramification occurs in the order
$\kappa\kappa\iota\iota$.
Thus quartics with the ramification $\kappa\kappa\iota\iota$ may have
both possibilities of zero or one real node.
Interestingly, all quartics with order $\kappa\iota\kappa\iota$
have one real node, which we explain below.

\begin{figure}[htb]
   \begin{picture}(72,72)
     \put(0,2){\epsfysize=70pt \epsfbox{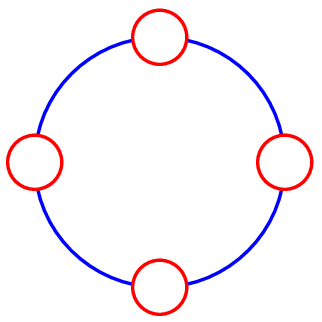}}
     \put(4,35){$\kappa$} \put(59,35){$\kappa$}
     \put(33,7){$\iota$} \put(33,62){$\iota$}
   \end{picture}
   \qquad
   \begin{picture}(72,72)
     \put(0,2){\epsfysize=70pt \epsfbox{figures/necklace4.eps}}
     \put(4,35){$\kappa$} \put(61,35){$\iota$}
     \put(31,7){$\kappa$} \put(33,62){$\iota$}
   \end{picture}
\caption{Necklaces for $\kappa=\iota=2$.\label{fig:necklaces}}
 \end{figure}

Questions of classifying maximally inflected plane rational curves can thus be
refined to take
into account the isotopy type of the ramification in $\mathbb{RP}^1$.
The different isotopy types of the placement of ramification data
$\alpha^1,\ldots,\alpha^n$ on $\mathbb{RP}^1$ are encoded by
combinatorial objects called necklaces:
$n$ `beads' with `colors' $\alpha^1,\ldots,\alpha^n$ are strung along
$S^1=\mathbb{RP}^1$ to make a necklace.
Given a maximally inflected curve
$\varphi\colon\mathbb{P}^1\to\mathbb{P}^2$, we may reverse the
parameterization of $\mathbb{RP}^1$ to obtain another maximally inflected curve
whose necklace is the mirror image of the necklace of the original curve.
Thus we identify two necklaces that differ only by the orientation of
$\mathbb{RP}^1$.
For example, Figure~\ref{fig:necklaces} displays the two necklaces with
2 beads each of color $\kappa$ and $\iota$.
To a maximally inflected curve with ramification $\alpha^1,\ldots,\alpha^n$, we
associate a necklace with beads of colors $\alpha^1,\ldots,\alpha^n$ where
the bead with color $\alpha^i$ is placed at the corresponding point of
ramification on $S^1=\mathbb{RP}^1$.

To state the promised generalization of
Theorem~\ref{thm:quartic-isotopy} let us denote
by $C(\Omega)$ the space of maximally inflected rational
quartics with a given necklace $\Omega$ and by $P(\Omega)$
the space of the placements of the necklace in $\mathbb{RP}^1$.
This 
generalization
may be proven in a manner similar to the proof of
Theorem~\ref{thm:quartic-isotopy}
(cf. Remark \ref{on-Zeuthen}), and so we omite its proof.

\begin{thm}\label{T:quart-isotopy}
 Maximally inflected rational quartic curves in $\RP^2$ admit arbitrary
 deformations, and the isotopy type of the image
 (together with the bitangents and solitary points)
 in $\RP^2$ is a deformation invariant.
 Moreover,
 for any necklace $\Omega$ the canonical projection $C(\Omega)\to P(\Omega)$
 is a trivial covering.

 In particular, the only 
isotopy types of the image are those indicated
 in Table 2 and for any given necklace
 the number of maximally inflected quartics
 having a given isotopy type
 of the image in $\mathbb{RP}^2$
 does not depend on the placement of the necklace.
\end{thm}

(Note that there are rational real quartics with one real
node and nested loops, but such a quartic cannot
have solitary points, and, as it follows from \ref{cor:bounds}, it cannot be
maximally inflected.)

We refine Question~\ref{ques:one}.

\begin{ques}\label{ques:two}
 Given a necklace $\Omega$ with beads of color $\alpha^1,\ldots,\alpha^n$,
 which are ramification data for degree $d$ plane curves, what are the
 possibilities for the numbers $\delta,\eta,c$ of solitary points, real nodes,
 and complex nodes of a maximally inflected curve of degree $d$ whose
 associated necklace is $\Omega$?
 Are any of these (or their positions) a deformation invariant?
\end{ques}

A weaker problem is to give bounds better than those of
Corollary~\ref{cor:bounds} for these possibilities.
For example, $\eta=0$ is not possible for the necklace on
the left in Figure~\ref{fig:necklaces}.
Indeed, the two curves for this necklace with
ramification points at $\{-1,0,1,\infty\}$ are obtained from each other
by reversing the parameterization, and so their images in
$\RP^2$ are equal.
In fact, this common image, which has a single node, is displayed the second
row of the third column of Table~\ref{table:two}.
Invoking Theorem~\ref{T:quart-isotopy} completes the proof.

We give
a more direct
proof that $\eta=0$ is not possible.
To see this, suppose such a curve has no real nodes.
Then the real locus is a topologically embedded
two-sided circle.
The inflection points are the points where the
concavity (which can be represented by
an oriented curvature covector taking zero on the tangent
direction) of the curves changes.
Note that the concavity does not change at a cusp,
and because of a flex between the cusps,
one cusp is pointed outward and another inward, as indicated below.
\[
  \epsfysize=70pt \epsfbox{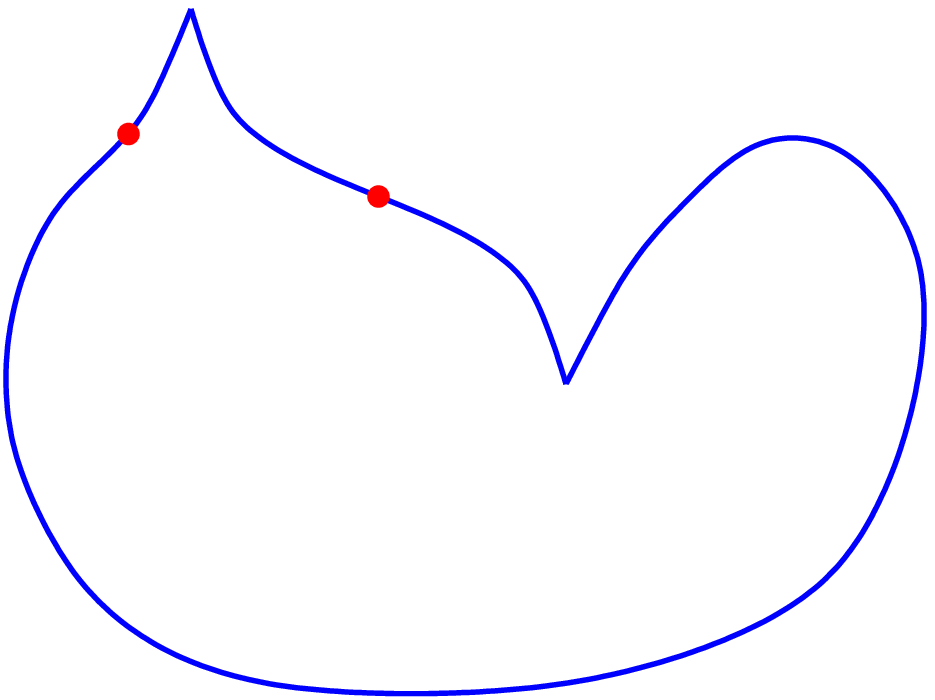}
\]
Consider now the line through the cusps.
It must meet the curve in at least 1 additional point, or have a local
intersection number of 3 with a cusp,
and thus it has intersection number
at least 5 with the quartic.
This contradicts B\'ezout's theorem, and proves the impossibility of such a
curve.\medskip

We now consider the possible positions of the nodes in a maximally
inflected curve.
There are 4 possibilities for the positions of the node with respect to the
flexes in a maximally inflected quartic with 6 flexes and one real node.
We display all 4:
$$
  \epsfysize=70pt \epsfbox{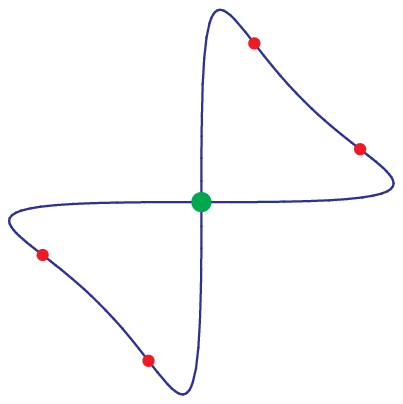}
  \qquad
  \epsfysize=70pt \epsfbox{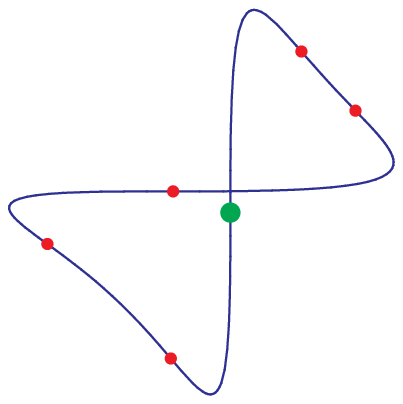}
  \qquad
  \epsfysize=70pt \epsfbox{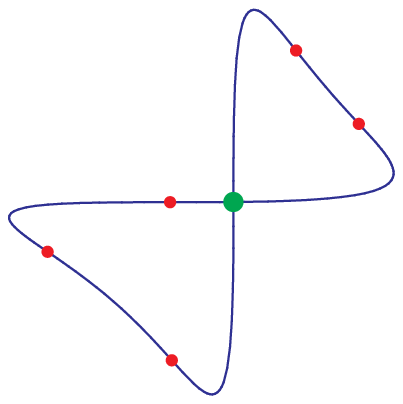}
  \qquad
  \epsfysize=70pt \epsfbox{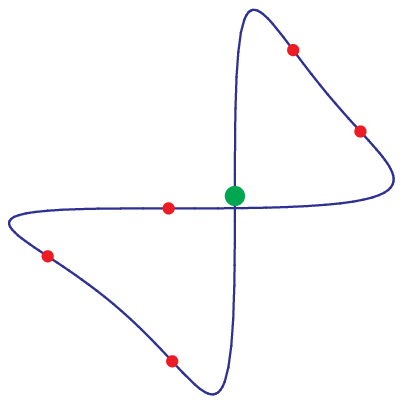}
$$
The position of the node may be represented in the associated necklace by
drawing a chord joining the two points whose images coincide:
$$
   \setlength{\unitlength}{.9pt}
   \begin{picture}(92,82)
     \put(0,0){\epsfysize=72pt \epsfbox{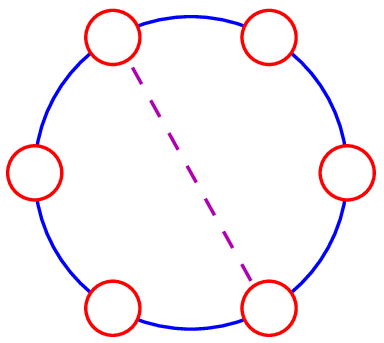}}
     \put(5,37){$\iota$} \put(80,37){$\iota$}
     \put(24,5){$\iota$} \put(24,69){$\iota$}
     \put(61,5){$\iota$} \put(61,69){$\iota$}
   \end{picture}\qquad
   \begin{picture}(92,82)
     \put(0,0){\epsfysize=72pt \epsfbox{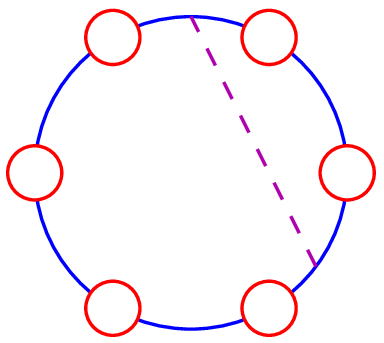}}
     \put(5,37){$\iota$} \put(80,37){$\iota$}
     \put(24,5){$\iota$} \put(24,69){$\iota$}
     \put(61,5){$\iota$} \put(61,69){$\iota$}
   \end{picture}\qquad
   \begin{picture}(92,82)
     \put(0,0){\epsfysize=72pt \epsfbox{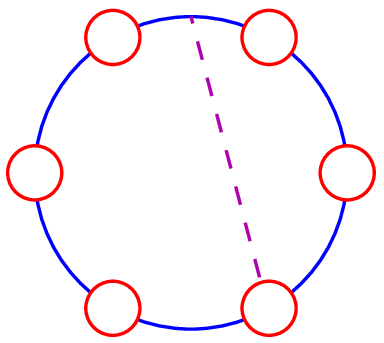}}
     \put(5,37){$\iota$} \put(80,37){$\iota$}
     \put(24,5){$\iota$} \put(24,69){$\iota$}
     \put(61,5){$\iota$} \put(61,69){$\iota$}
   \end{picture}\qquad
   \begin{picture}(92,82)
     \put(0,0){\epsfysize=72pt \epsfbox{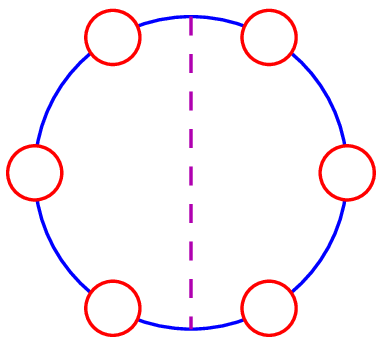}}
     \put(5,37){$\iota$} \put(80,37){$\iota$}
     \put(24,5){$\iota$} \put(24,69){$\iota$}
     \put(61,5){$\iota$} \put(61,69){$\iota$}
   \end{picture}
$$
These are the only 4 possibilities:
Since a bitangent to a quartic can neither be tangent at a flex nor meet another
point of the curve, the flexes and nodes of such a quartic are constrained to
lie on the arcs of the curve in between the  contacts of its bitangents.
Since, by Theorem~\ref{thm:quartic-isotopy}, every maximally inflected quartic with a
single node can be deformed into one of those shown in Figure~\ref{fig:6flex},
this constraint rules out the other possibilities for the chord in a necklace of
such a quartic.

We further refine Questions~\ref{ques:one} and~\ref{ques:two}.

\begin{ques}\label{ques:three}
 Which necklaces $\Omega$ with beads of color $\alpha^1,\ldots,\alpha^n$ and
 $\eta$ chords occur as maximally inflected curves?
 Given such a chord  diagram, what are the possible numbers of solitary
 points (and hence complex nodes)?

 Ignoring the beads, we obtain a circle with $\eta$ chords, and the same
 questions may be asked of these diagrams.
 Such pure chord diagrams encode, together with
 the number of solitary points, the topology of the image as
 an abstract, not embedded in $\mathbb{RP}^2$, topological space,
 and therefore the classification of chord diagrams is a
 necessary part of
 any topological classification.
\end{ques}

We make one final observation concerning the orientation of the cusp of a
maximally inflected quartic with four flexes and one cusp.
The examples in Table~\ref{table:two} both have the cusp pointing into the
unbounded region in the complement of the curve.
This is necessarily the case.
If a quartic had a cusp pointing into a bounded region, then it could
not have
real nodes or solitary points, as the line joining a cusp with such a real
double point would meet the curve in at least two additional points, and
thus contradicts B\'ezout's theorem.
But Corollary~\ref{cor:bounds} requires such a maximally inflected curve to
have at least one solitary point.
While we are presently unable to formulate a general question concerning
restrictions on the disposition of cusps (and other ramification), it is
likely there are further topological restrictions.

Further pictures of maximally inflected quartics with more general
ramification and
many quintics may be found on
the web\frankfootnote{See
{\tt www.math.umass.edu/\~{}sottile/pages/inflected}}.

\section{Maximally inflected plane quintics}

Unless explicitly stated otherwise, all quintics in this section are assumed to
have cusps and flexes as their only ramification points and to have ordinary
double points as their only other singularities.
In Section~\ref{S:kkkiii}, we additionally assume that all quintics 
have three real nodes.

\subsection{Quintics with cusps and flexes having three nodes}\label{S:kkkiii}
The construction of Theorem~\ref{thm:classical} for quintics gives maximally
inflected quintics whose number $\kappa$ of cusps is $0,1,2$, or $3$.
We obtain a quintic with the maximal number of 4 cusps by taking the dual
of a maximally inflected quartic with one cusp and 4 flexes, whose existence
was addressed in Theorem~\ref{thm:4-class}.
The existence of such curves is also a consequence of Theorem~\ref{thm:simpleRam},
as the ramification indices of cusps and flexes are special, in the sense of
that theorem.

Since Theorem~\ref{thm:simpleRam} is an existence result and gives no
information about the geometry of the resulting curve, it is very instructive to
look at specific examples coming from constructions.
As an example, consider the possible necklaces of such curves, which is
interesting only for two or three cusps.
There are three possible necklaces of maximally inflected quintics with
2 cusps and 5 flexes, 
 \[
   \begin{picture}(72,72)
     \put(0,2){\epsfysize=70pt \epsfbox{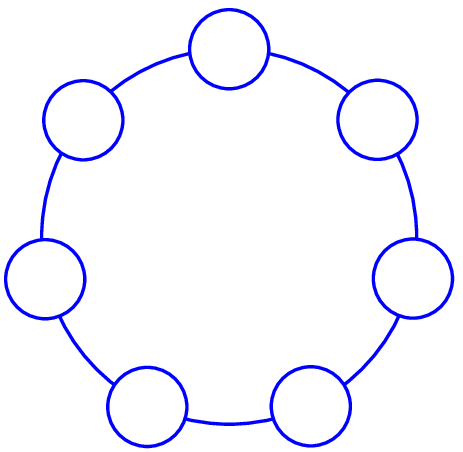}}
                  \put(32.5,62.5){$\kappa$}
     \put(11,51){$\iota$}         \put(57,51){$\iota$}
     \put(5,26.5){$\iota$}        \put(61,26.5){$\kappa$}
     \put(20.5,6){$\iota$}      \put(46.5,6){$\iota$}
   \end{picture}
   \qquad
   \begin{picture}(72,72)
     \put(0,2){\epsfysize=70pt \epsfbox{figures/necklace7.eps}}
                  \put(32.5,62.5){$\kappa$}
     \put(11,51){$\iota$}         \put(56,51){$\kappa$}
     \put(5,26.5){$\iota$}        \put(62,26.5){$\iota$}
     \put(20.5,6){$\iota$}      \put(46.5,6){$\iota$}
   \end{picture}
   \qquad
   \begin{picture}(72,72)
     \put(0,2){\epsfysize=70pt \epsfbox{figures/necklace7.eps}}
                   \put(32.5,62.5){$\kappa$}
     \put(11,51){$\iota$}         \put(57,51){$\iota$}     
    \put(5,26.5){$\iota$}        \put(62,26.5){$\iota$}
     \put(20.5,6){$\iota$}      \put(45.5,6){$\kappa$}
  \end{picture}
 \]
and three possible necklaces  with 3 cusps and 3 flexes.
 \[
   \setlength{\unitlength}{.9pt}
   \begin{picture}(92,82)
     \put(0,0){\epsfysize=72pt \epsfbox{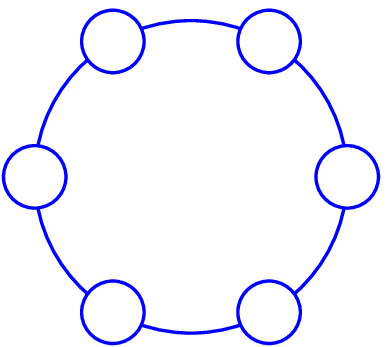}}
     \put(5,37){$\iota$} \put(79,37){$\kappa$}
     \put(23,5){$\kappa$} \put(23,69){$\kappa$}
     \put(61,5){$\iota$} \put(61,69){$\iota$}
   \end{picture}\qquad
   \begin{picture}(92,82)
     \put(0,0){\epsfysize=72pt \epsfbox{figures/necklace6.eps}}
     \put(5,37){$\iota$} \put(80,37){$\iota$}
     \put(23,5){$\kappa$} \put(23,69){$\kappa$}
     \put(61,5){$\iota$} \put(60,69){$\kappa$}
   \end{picture}\qquad
   \begin{picture}(92,82)
     \put(0,0){\epsfysize=72pt \epsfbox{figures/necklace6.eps}}
     \put(5,37){$\iota$} \put(79,37){$\kappa$}
     \put(24,5){$\iota$} \put(23,69){$\kappa$}
     \put(61,5){$\iota$} \put(60,69){$\kappa$}
   \end{picture}
 \]
The leftmost necklace shown for both values of $\kappa$ is that of
the corresponding quintic in Figure~\ref{fig:shustin}.
We can realize three of the remaining four necklaces using the variant of the
construction of Theorem~\ref{thm:classical} where we use a different local model
for the perturbation of a tacnode into a cusp, as explained in
Remark~\ref{rem:alt-model}.
This gives the three quintics displayed below
$$
  \epsfysize=90pt\epsfbox{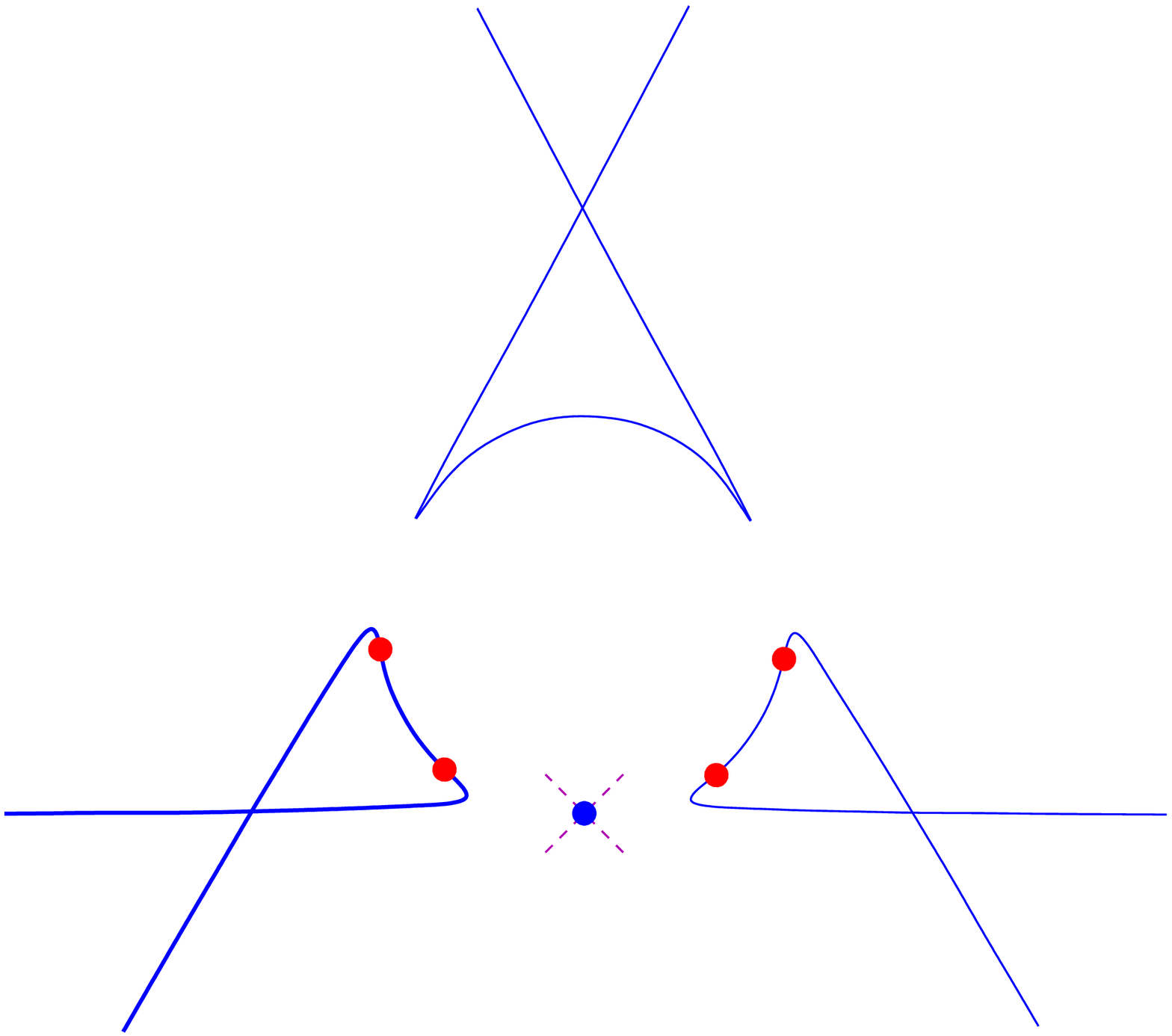}\qquad
  \epsfysize=90pt\epsfbox{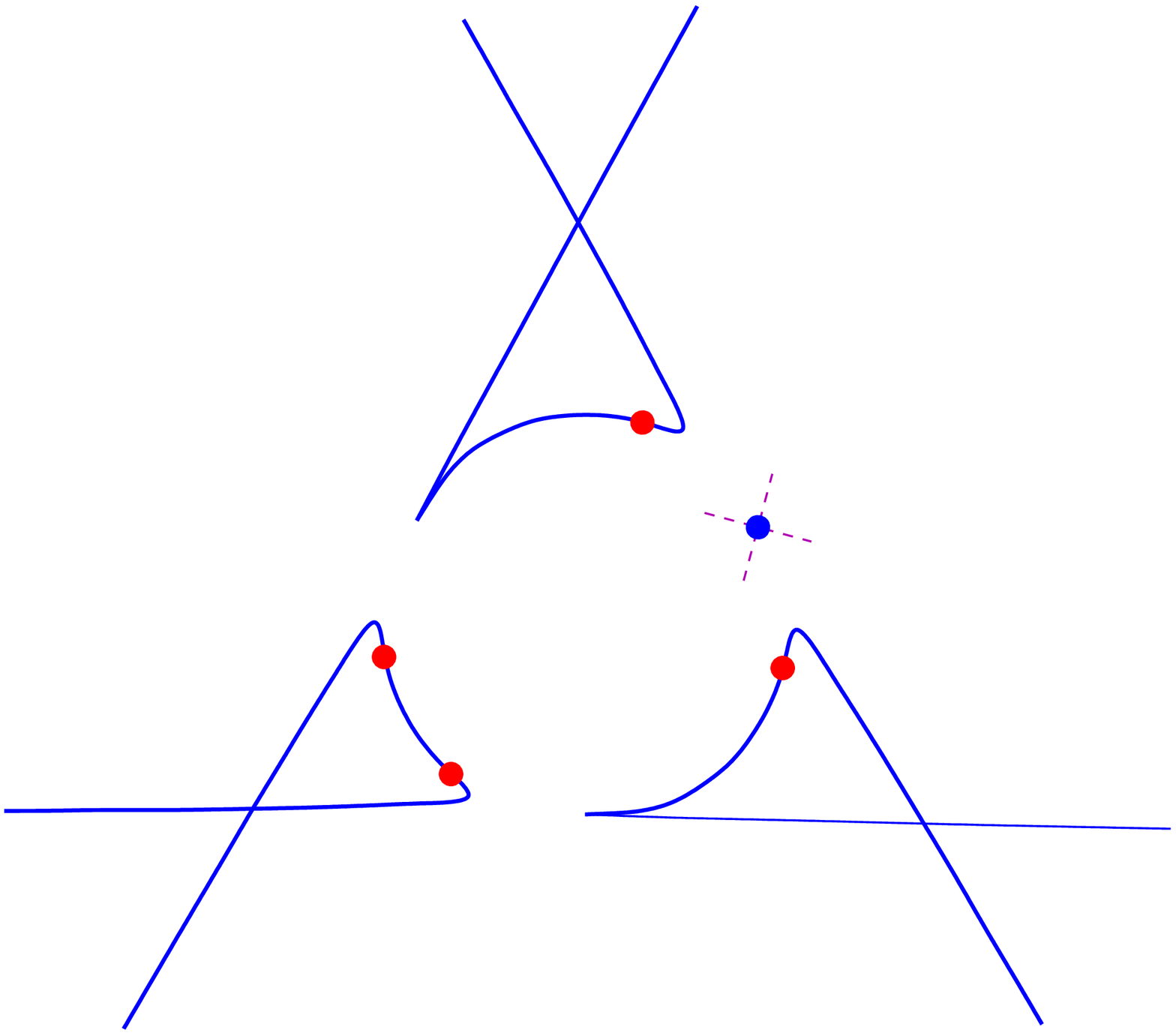}\qquad
  \epsfysize=90pt\epsfbox{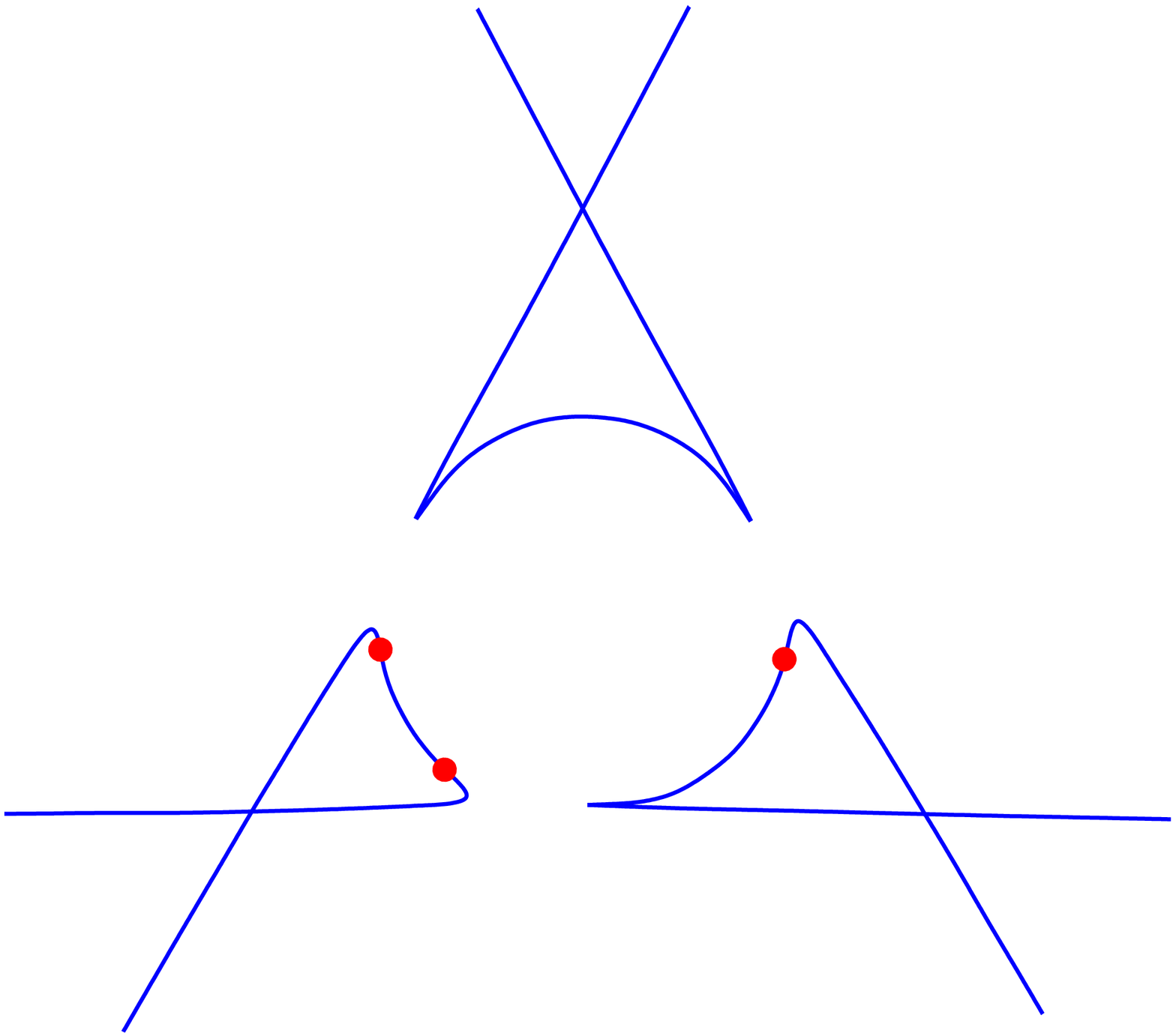}
$$
which correspond to the remaining two necklaces of curves with two cusps,
and the second necklace for curves with three cusps.

What is missing is a quintic with three consecutive cusps.
There are six maximally inflected plane quintics with cusps at $\infty,\pm 3$
and flexes at $0, \pm 1$, which we have computed using symbolic methods.
One of these six is particularly interesting.
The two pictures on the left below are two different views of this curve.
In the first, we put one flex at infinity, and in the second, one cusp at
infinity.
$$
  \epsfysize=90pt\epsfbox{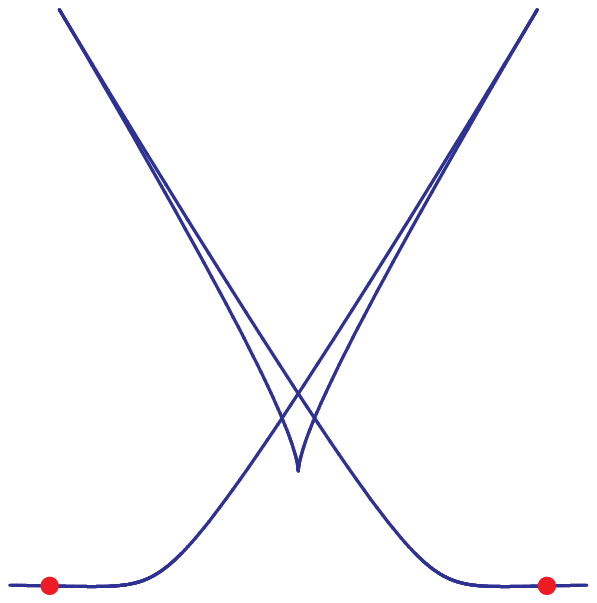}\qquad
  \epsfysize=90pt\epsfbox{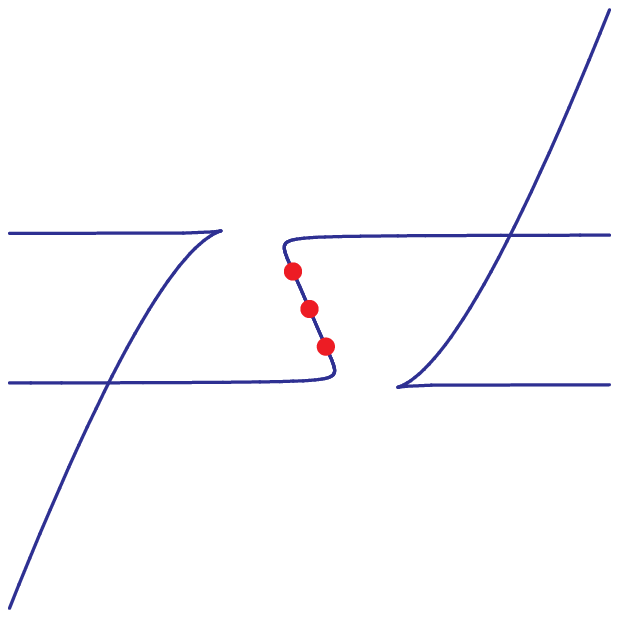}\qquad
  \epsfysize=90pt\epsfbox{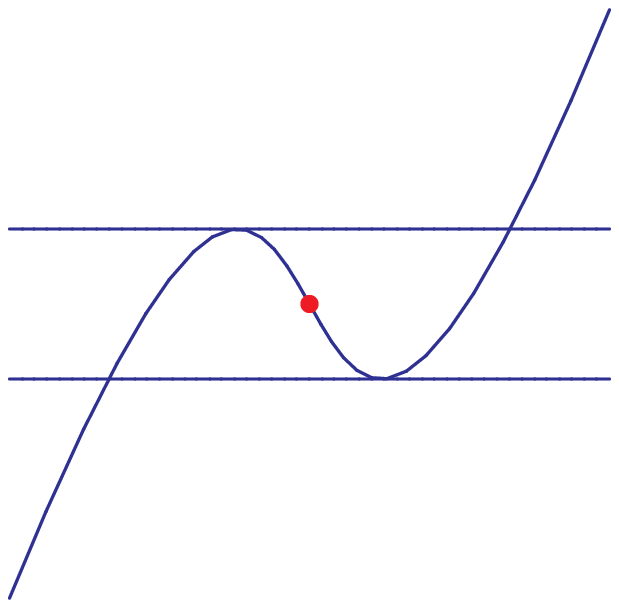}
$$
The second view of this curve suggests that it is a perturbation
of the singular
and reducible curve shown on the right.

\subsection{Solitary points of quintics with cusps and
flexes}\label{S:kis}
In this section, we
address the finer classification of Question~\ref{ques:two} concerning the
possible numbers of solitary points and nodes in a maximally inflected quintic
having flexes and cusps with a given necklace.
Here, the situation is quite intricate and our knowledge at present is not
definitive.
It is, however, based upon extensive experimental numerical evidence.
Briefly, for quintics with either zero or one cusp, all possibilities occur,
but for most necklaces with two or three cusps, we have not yet seen all
possibilities for the numbers of solitary points.
We have also observed a fascinating global phenomenon, related to (but weaker
than) the very strong classification result we obtained for quartics in
Theorem~\ref{T:quart-isotopy}.

Let us recall the results of Corollary~\ref{cor:bounds} for quintics with
flexes and cusps.

\begin{figure}[htb]
\[
    \epsfysize=70pt\epsfbox{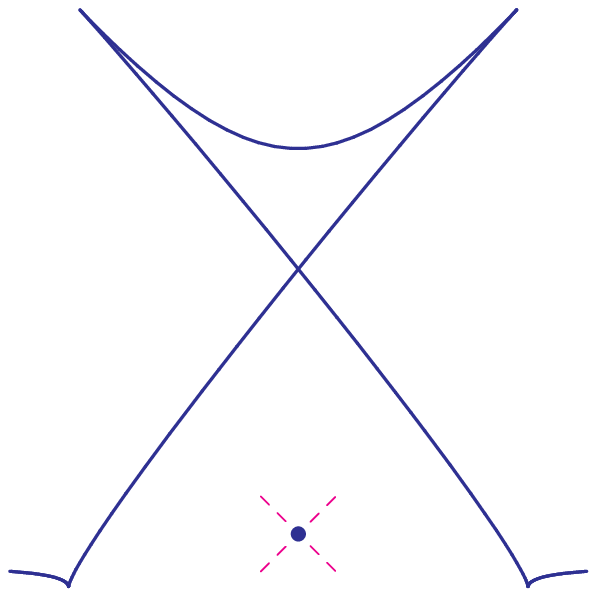}\qquad\qquad
    \epsfysize=70pt\epsfbox{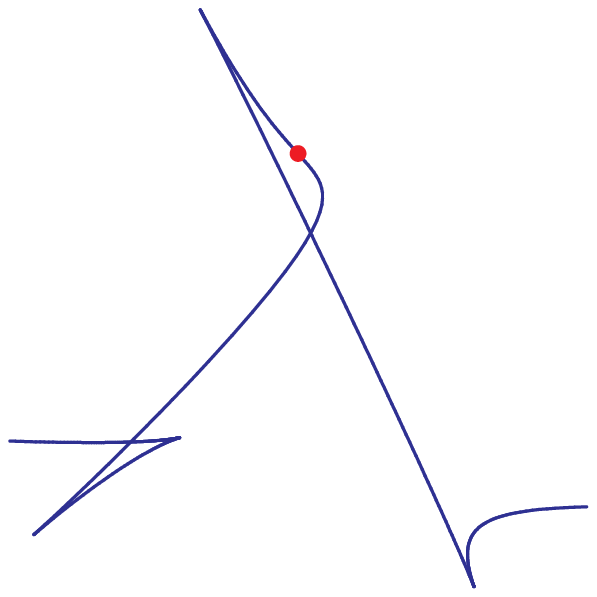}
\]
\caption{Quintics with four cusps}\label{F:4cusps}
\end{figure}

A rational quintic plane curve with $\kappa$ cusps and $9-2\kappa$ flexes has
genus $g_\kappa=6-\kappa$.
Maximally inflected quintics with four cusps ($\kappa=4$) are dual to
maximally inflected quartics with four flexes and one cusp, and were discussed
from different points of view
in Remarks~\ref{k4.1},~\ref{k4.2}, and~\ref{k4.3}.
In particular, there we showed the impossibility of such a quintic with no
solitary points and no nodes.
In Figure~\ref{F:4cusps}, we give examples of the two remaining possibilities
for maximally quintics with four cusps and one flex.
The quintic on the left has its flex at infinity.

Consider quintics with three or fewer cusps.
By Corollary~\ref{cor:bounds}, the number $\delta$ of solitary points of such
a quintic satisfies
 \[
    3-\kappa\ \leq \ \delta\ \leq\ 6-\kappa\ =:\ g_\kappa\,,
 \]
or, $3\leq \delta+\kappa\leq 6$.
The number $c$ of complex nodes is an even number between 0 and
$g_\kappa-\delta$, and the number $\eta$ of nodes satisfies
$\eta+c=g_\kappa-\delta$.
Thus the possibilities for these numbers for quintics are in bijection with
the cells of the diagram
 \[
   \begin{picture}(80,55)
    \put(29,30){3} \put(42,30){4} \put(55,30){5} \put(68,30){6}
    \put(12,16){0} \put(12, 3){2}
    \put(0,9){$c$}  \put(40,45){$\delta+\kappa$}

    \put(25,0){\begin{picture}(40,20)
       \put( 0,26){\line(1,0){52}}
       \put( 0,13){\line(1,0){52}}
       \put( 0, 0){\line(1,0){26}}
       \put( 0, 0){\line(0,1){26}}
       \put(13, 0){\line(0,1){26}}
       \put(26, 0){\line(0,1){26}}
       \put(39,13){\line(0,1){13}}
       \put(52,13){\line(0,1){13}}
      \end{picture}}
   \end{picture}
 \]
with the first and second rows corresponding to the values of 0 and 2 for $c$,
and the columns correspond to the possible vales for 
$\delta+\kappa$, from $3$ on the left to $6$ on the right.
The first column of Table~\ref{T:last} lists the possible necklaces for
quintics with at most three cusps.
The filled cells of the diagrams in its second column indicate the observed
cardinalities of solitary points, complex nodes, and nodes.

\def\k{\kappa}
\def\i{\iota}
\begin{table}[htb]
\begin{tabular}{|c||c|cccc|c|c|}\hline
\multirow{2}{45pt}{Necklace}&   Observed    &\multicolumn{4}{c|}{$\delta+\k$}&
\multirow{2}{10pt}{$N$}&Number\\
                       &$\delta,\eta,c$&3&4&5&6&&Tested\\\hline\hline
 $\k\k\k\i\i\i$&
%
  \raisebox{-9pt}{\epsffile{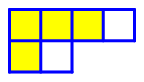}}&
 1$\!$&$\!$3$\!$&$\!$2$\!$&$\!$0\;& 6& 985 \\\hline
 $\k\k\i\k\i\i$&
%
  \raisebox{-9pt}{\epsffile{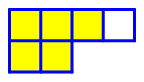}}&
 2$\!$&$\!$3$\!$&$\!$1$\!$&$\!$0& 6& 986\\\hline
 $\k\i\k\i\k\i$&
  \raisebox{-9pt}{\epsffile{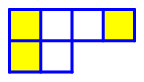}}&
 5$\!$&$\!$0$\!$&$\!$0$\!$&$\!$1& 6& 986\\\hline\hline
 $\k\k\i\i\i\i\i$&
%
  \raisebox{-9pt}{\epsffile{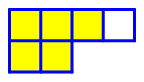}}&
 3$\!$&$\!$5$\!$&$\!$3$\!$&$\!$0& 11& 731\\\hline
 $\k\i\k\i\i\i\i$&
  \raisebox{-9pt}{\epsffile{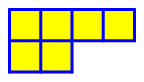}}&
 5$\!$&$\!$4$\!$&$\!$1$\!$&$\!$1& 11& 731 \\\hline
 $\k\i\i\k\i\i\i$&
%
  \raisebox{-9pt}{\epsffile{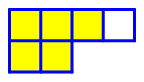}}&
 4$\!$&$\!$5$\!$&$\!$2$\!$&$\!$0& 11& 731\\\hline\hline
 $\k\i\i\i\i\i\i\i$&
  \raisebox{-9pt}{\epsffile{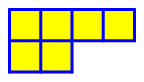}}&
 7$\!$&$\!$9$\!$&$\!$4$\!$&$\!$1& 21&  1000\\\hline\hline
 $\i\i\i\i\i\i\i\i\i$&
  \raisebox{-9pt}{\epsffile{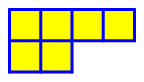}}&
 12&18&9&3& 42& 11,416\\\hline
\end{tabular}\smallskip
\caption{Observed numbers of nodes and solitary points of quintics}
\label{T:last}
\end{table}

The third column of Table~\ref{T:last} records an interesting global
phenomenon we have observed.
Recall from Section~\ref{S:1.2} that for a given collection of ramification
data $\alpha^1,\ldots,\alpha^n$, there will be the same number
$N:=N(\alpha^1,\ldots,\alpha^n)$ of rational
curves having that ramification at specified points.
We observe that that the number of these $N$ curves
having a given number of solitary points appears to depend only upon the
necklace, and not on the placement of the ramification.
We record this in the third column, with the columns corresponding to the
possible values of $\delta+\kappa$ between 3 and 6, from left to right.
The last two columns of Table~\ref{T:last} give the number
$N(\alpha^1,\ldots,\alpha^n)$ of rational quintics of that type having a given
choice of ramification, and the number of choices of ramification for which we
computed all $N(\alpha^1,\ldots,\alpha^n)$ quintics and determined their
numbers of solitary points, complex nodes, and real nodes.\medskip

We make the observation in the third column of Table~\ref{T:last} more precise
in the following conjecture.

\begin{conj}\label{C:quintic-solitary}
  The number of solitary points in a maximally inflected plane quintic is a
  deformation invariant.
\end{conj}

In particular, if Conjecture~\ref{conj:non-degen} (concerning
non-degeneracy) held for quintics, then the number of
curves having a given ramification and given number of solitary points depends
only on the necklace of the ramification, that is, only on the relative
positions of the points of ramification.
This is similar to the conclusion of Theorem~\ref{T:quart-isotopy} for
quartics, but it is a weaker phenomenon, as the topology of the embedding of a
quintic can change under an isotopy.  (We give an example below.)
If Conjecture~\ref{conj:non-degen} held for quintics, then similar arguments
as given in the proof of Theorem~\ref{T:quart-isotopy} concerning isolating
solitary points by  bitangents and B\'ezout's theorem may suffice to prove
Conjecture~\ref{C:quintic-solitary}.

An instructive example is provided by quintics whose ramification consists of
three cusps at $\pm1$ and $\infty$ and three flexes at 0 and $\pm t$, where
$1<t<\infty$.
For these the cusps and flexes alternate, and for a given choice of
ramification, there are 6 curves.
Letting $t$ vary, we have 6 families of curves, which are deformations of the
curves at any value of $t$.
All the curves in one family have three solitary points,
while those in the other five families have no solitary points.
Two of these 
five
families consist of curves with three nodes.
When $t=3$, each curve in the remaining three families has two smooth branches
with a point of contact of order three.
The other curves in two of these  three
families each have a single node and two
complex nodes.
The latter two families differ only by
reparameterization of their curves.
We display curves from one  of them at the values
$t=2.9, 3$, and $3.1$.
All three have a cusp at infinity.
The line joining the complex nodes is drawn in the first and third pictures.
Despite appearences, it is not tangent to the curve.  
If it were tangent, then it would have intersection number at least 6 with the
curve, which contradicts B\'ezout's theorem.
 \[
    \epsfysize=72pt\epsfbox{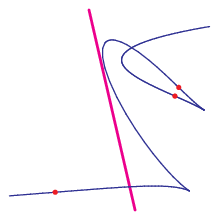}\qquad\qquad
    \epsfysize=68pt\epsfbox{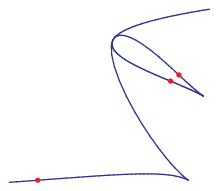}\qquad\qquad
    \epsfysize=72pt\epsfbox{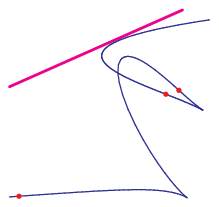}
 \]
The curves in the third
family have one node and two complex
nodes when $t>3$, but three nodes for $t<3$.
We display curves from this family at the values $t=\frac{5}{2}, 3$,
and $\frac{7}{2}$.
The horizontal line in the third picture is the real line meeting the two
complex nodes, and all three have a flex at infinity.
 \[
    \epsfysize=72pt\epsfbox{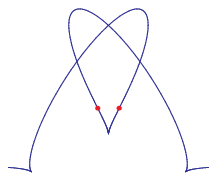}\qquad\qquad
    \epsfysize=68pt\epsfbox{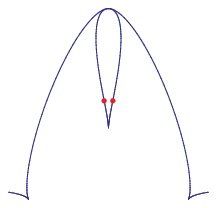}\qquad\qquad
    \epsfysize=72pt\epsfbox{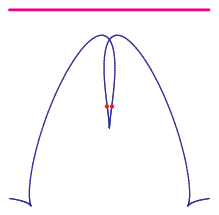}
 \]
We remark that these curves, like those shown in Section 4 and in
Figure~\ref{fig:quintics} below, were drawn using the computer algebra systems
Maple and Singular~\cite{Singular}.
We first computed the centers of projection giving all curves with a given
ramification, formulating and solving the problem in local coordinates for
the Grassmannian.
Given the centers, we computed parameterizations and then set up and solved
the equations for the double points.
Those whose embedding we needed to study, we plotted, and for those displayed
in this paper, we created postscript files which
were then further edited to enhance important features, such as flexes and
solitary points.

\subsection{Chord diagrams of quintics and embeddings}

We classify the possible chord diagrams of maximally inflected quintics and
discuss the possible topology of the image as a subset of $\mathbb{RP}^2$
(ignoring possible solitary points).
Since the proof we give does not use flexes, it also classifies chord
diagrams of quintics with at most three nodes.

\begin{thm}
 There are 6 possibilities for the chord diagram of a maximally inflected
 quintic in $\mathbb{RP}^2$.
 Here are the four with 2 or more chords.
$$
  \epsfysize=50pt\epsfbox{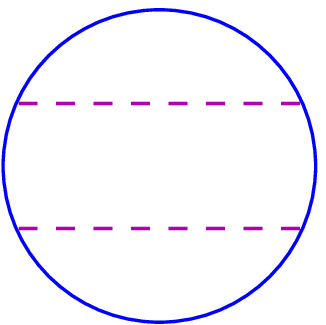}
  \qquad
  \epsfysize=50pt\epsfbox{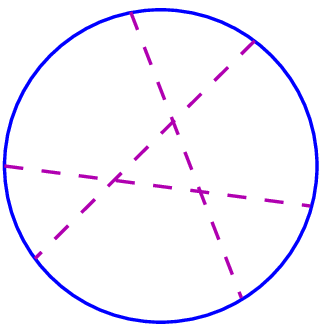}
  \qquad
  \epsfysize=50pt\epsfbox{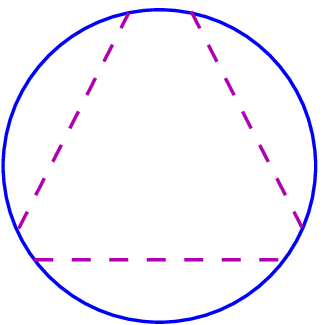}
  \qquad
  \epsfysize=50pt\epsfbox{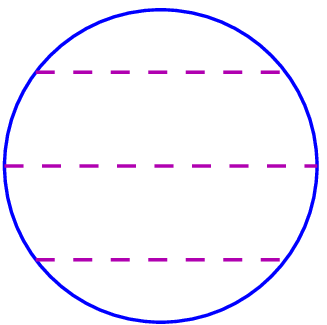}
$$
\end{thm}

{\bf Proof. }
 By Corollary~\ref{cor:bounds}, a maximally inflected quintic can have at most
 three nodes.
 Recall that the degree of any (piecewise smooth) closed curve in
 $\mathbb{RP}^2$ is well-defined modulo 2, and this degree is additive,
 again modulo 2.
 Consider an image of $\mathbb{RP}^1$ in $\mathbb{RP}^2$ with a node.
 Under an orientation from $\mathbb{RP}^1$, the node has two incoming
 arcs and two outgoing arcs.
 We may split the curve into two pieces at the node by joining an incoming
 arc of one branch with the outgoing arc of the other branch.
 We illustrate this splitting below.
$$
  \epsfysize=40pt\epsfbox{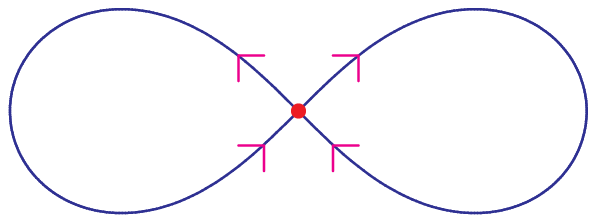}
  \quad \raisebox{18pt}{$\Longrightarrow$} \quad
  \epsfysize=40pt\epsfbox{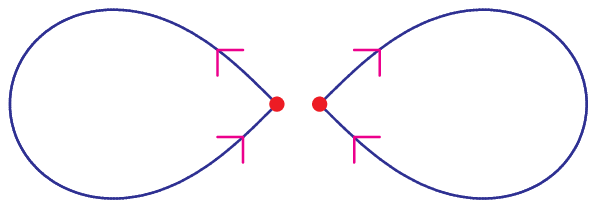}
$$

 The first case to consider is that of a quintic with two nodes.
 In a chord diagram of such a quintic, the two chords cannot cross.
 If they did, then split the quintic at a node.
 The resulting two closed curves have exactly one point of intersection
 (the other node), and so each piece necessarily has odd degree.
 It follows that their union, the original quintic, has even degree, a
 contradiction.

 Suppose now that we have a quintic with three nodes.
 The argument we just gave forbids any chord diagrams containing a chord that
 meets exactly one other chord, and thus the only possibilities are as claimed.
 Concluding, we get only one chord diagram with 2 chords
 (two non-intersecting chords) and three chord diagrams with 3
 chords (three pairwise intersecting chords, three sides of a hexagon,
 and three parallel chords).
\QED

Recall that a pseudoline is a closed curve in $\mathbb{RP}^2$ which
has odd degree 
and
whose complement is connected (it is a one-sided curve).
An oval is a two-sided closed curve which necessarily has even degree.
We enumerate the possible topological embeddings of $\mathbb{RP}^1$ given by a
maximally inflected quintic, beginning with the classification of chord
diagrams.

\begin{figure}[htb]
\[
  \begin{picture}(300,90)
   \put( 0,20){\epsfysize=70pt\epsfbox{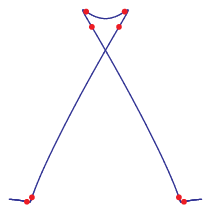}}
   \put( 30, 2){(a)}

   \put(100,20){\epsfysize=70pt\epsfbox{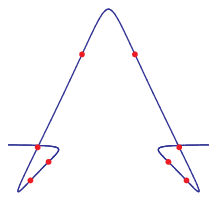}}
   \put(130, 2){(b)}

   \put(200,20){\epsfysize=70pt\epsfbox{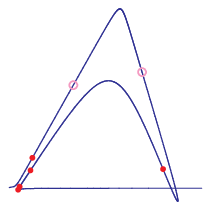}}
   \put(230, 2){(c)}
  \end{picture}
\]
\[
  \begin{picture}(300,100)
   \put(  0,20){\epsfysize=70pt\epsfbox{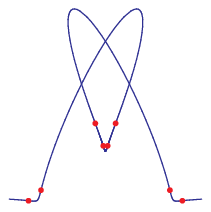}}
   \put( 30, 2){(d)}

   \put(100,20){\epsfysize=70pt\epsfbox{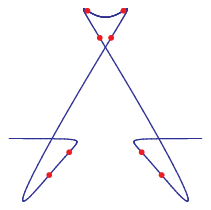}}
   \put(130, 2){(e)}

   \put(200,20){\epsfysize=70pt\epsfbox{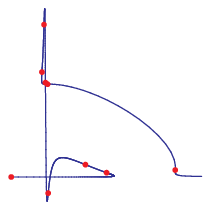}}
   \put(230, 2){(f)}
  \end{picture}
\]
\caption{Maximally inflected quintics realizing different embeddings}
\label{fig:quintics}
\end{figure}

If there is a single node, then the curve must look like a pseudoline with a
loop.
Splitting the curve at the node and applying the Jordan Curve Theorem to the
loop in the complement of the pseudoline shows that the loop is two-sided.
Thus there is a unique possible topological embedding.
Such a maximally inflected quintic is shown in
Figure~\ref{fig:quintics}(a).
We do not draw the solitary points.\smallskip

Now consider a quintic with two nodes whose chord diagram (necessarily)
consists of two non-intersecting chords.
We split the curve at both nodes to obtain three closed curves that meet only at
the nodes.
Deforming them slightly away from the nodes, we see that 
at most one is a pseudoline
(for any two pseudolines meet), and so exactly one is a
pseudoline, as the original curve was a quintic.
As before, the other components are then 2-sided ovals.
These ovals cannot be nested.
A line through a solitary point or a singular ramification point (a maximally
inflected quintic must have one such point) that meets the nest has intersection
number at least 6 with the quintic, contradicting  B\'ezout's theorem.

Thus there are two possibilities.
Either the pseudoline meets both ovals (see Figure~\ref{fig:quintics}(b)),
or the pseudoline meets only one oval, which then meets the other
(see Figure~\ref{fig:quintics}(c)).
As before the ovals cannot be nested.
Each open circle in Figure~\ref{fig:quintics}(c) represents two
flexes that have nearly merged to create a planar point.\smallskip

Suppose now that there are three pairwise intersecting chords.
If we split the curve at {\it one} node, we obtain a pseudoline and an oval,
which have two additional points of intersection.
There is only one possibility for this configuration, and we have already
seen it in the last pictures in each of Sections~\ref{S:kkkiii}
and~\ref{S:kis}.
We display yet another such curve below in
Figure~\ref{fig:quintics}(d).\smallskip

\begin{figure}[htb]
\[
  \begin{picture}(250,70)
   \put( 0,20){\epsfysize=50pt\epsfbox{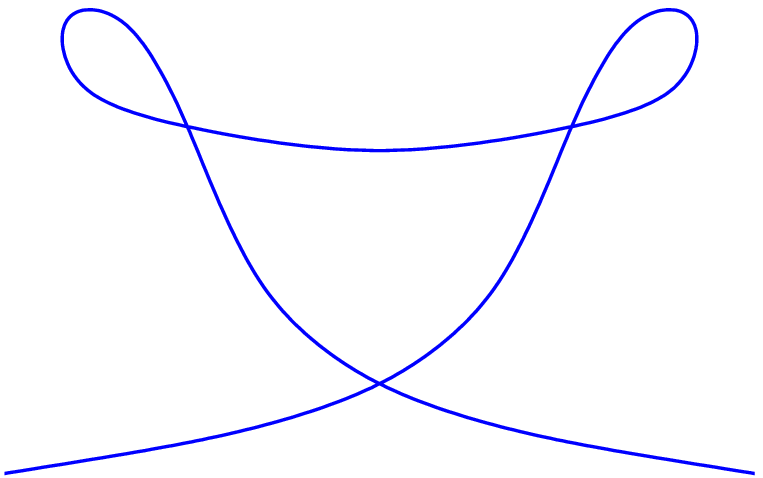}}
   \put(28, 2){(a)}

   \put(145,20){\epsfysize=50pt\epsfbox{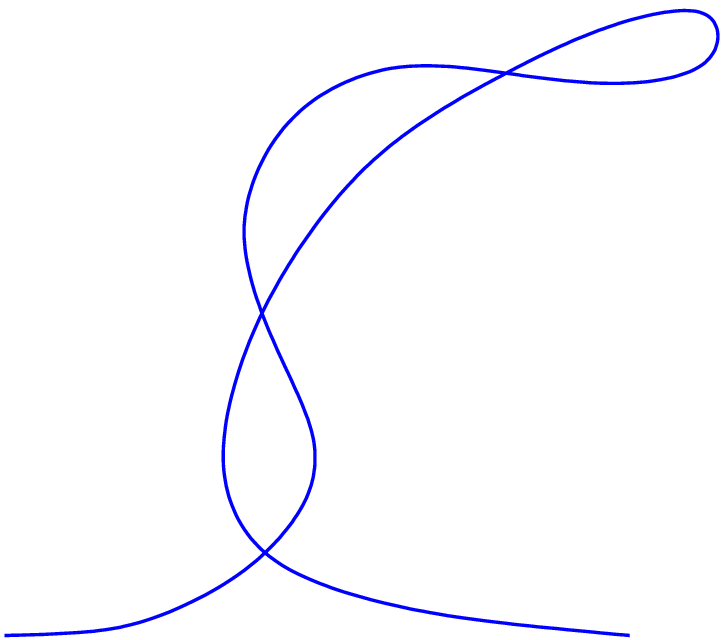}}
   \put(165, 2){(b)}
  \end{picture}
\]
\caption{Possible embeddings for quintics that have not been observed}
\label{fig:notfound}
\end{figure}

If the three chords are three sides of a hexagon, we may
split the curve at each of its nodes to obtain a pseudoline and three loops.
One piece meets the other three, and these three are disjoint from each other.
If the pseudoline meets the three loops, there are two possibilities for the
disposition of the three loops along the pseudoline.
Either they alternate (as shown in
Figure~\ref{fig:quintics}(e) or in Figure~\ref{fig:shustin}) or they do not.
We forbid the possibility of non-alternating loops in
Proposition~\ref{P:forbid} below.
We have not yet observed
such a quintic where the pseudoline meets one loop, which then meets the other
two loops.
A schematic for this is Figure~\ref{fig:notfound}(a).\smallskip

Lastly, we have the possibility of three parallel chords.
As before, we split the curve at each of its nodes to obtain one pseudoline and
three ovals, and the ovals cannot be nested.
Two components meet two others, and two meet exactly one other.
Either the pseudoline either meets two ovals
(See Figure~\ref{fig:quintics}(f)), or it meets only one.
This last case has not yet been observed, but we provide a schematic of it in
Figure~\ref{fig:notfound}(b)

\begin{prop}\label{P:forbid}
 A maximally inflected quintic consisting of a pseudoline with three loops
 must have the loops alternating.
\end{prop}

\noindent{\bf Proof. }
 Suppose that we have a maximally inflected quintic whose embedding
 consists of a pseudoline with three loops attached to the pseudoline.
 Such a curve has three nodes.
 By Corollary~\ref{cor:generalbounds}, a maximally inflected quintic with three
 nodes may have only flexes, 
 planar points, cusps, or a point with ramification sequence $(0,1,5)$.
 Let $\kappa$ be its number of cusps, which is at most 3---otherwise
 Corollary~\ref{cor:bounds} restricts the number of nodes to be at most
 2.
 Since it must have at least $3-\kappa$ solitary points and $6-\kappa$ double
 points in all, the curve has $3-\kappa$ solitary points, or 3 solitary points
 plus cusps.

 Gudkov's extension~\cite{Gu62} of Brusotti's Theorem~\cite{Br1921} allows us to
 independently smooth the cusps and nodes to obtain a smooth quintic.
 We smooth each cusp as in the local model of Figure~\ref{Fig:smooth}, smooth
 each node to detach its loop from the psuedoline, and smooth each
 solitary point to obtain an oval.
 Thus we obtain a smooth plane quintic with 6 ovals, which is an
 $M$-curve. 
 Furthermore, the three ovals arising from solitary points or cusps are
 distinguished.

Pick a point inside each of the three distinguished ovals, connect each pair by
a line and study the configuration of these lines with respect to the quintic.
By B\'ezout's theorem, each of these bisecants intersects the one-sided component
at one point and its two ovals at two points each. Hence, the configuration does not
depend on the choice of the points. Since  all $M$-quintics are deformation equivalent,
i.e., belong to one connected family of nonsingular quintics
(see \cite{Kh78}, a little correction of
the proof is given in \cite{DIK00}), the configuration does not depend on the choice
of the $M$-quintic either. Examining any of the various known examples, one finds
the configuration like the picture below, where we have drawn
the lines joining the solitary points of
 a rational quintic (the two pictures are for
 similar curves in different projections), 
 leaving the smoothing to the imagination of the reader.
 These solitary points are indicated by open circles.
 \[
    \epsfysize=90pt\epsfbox{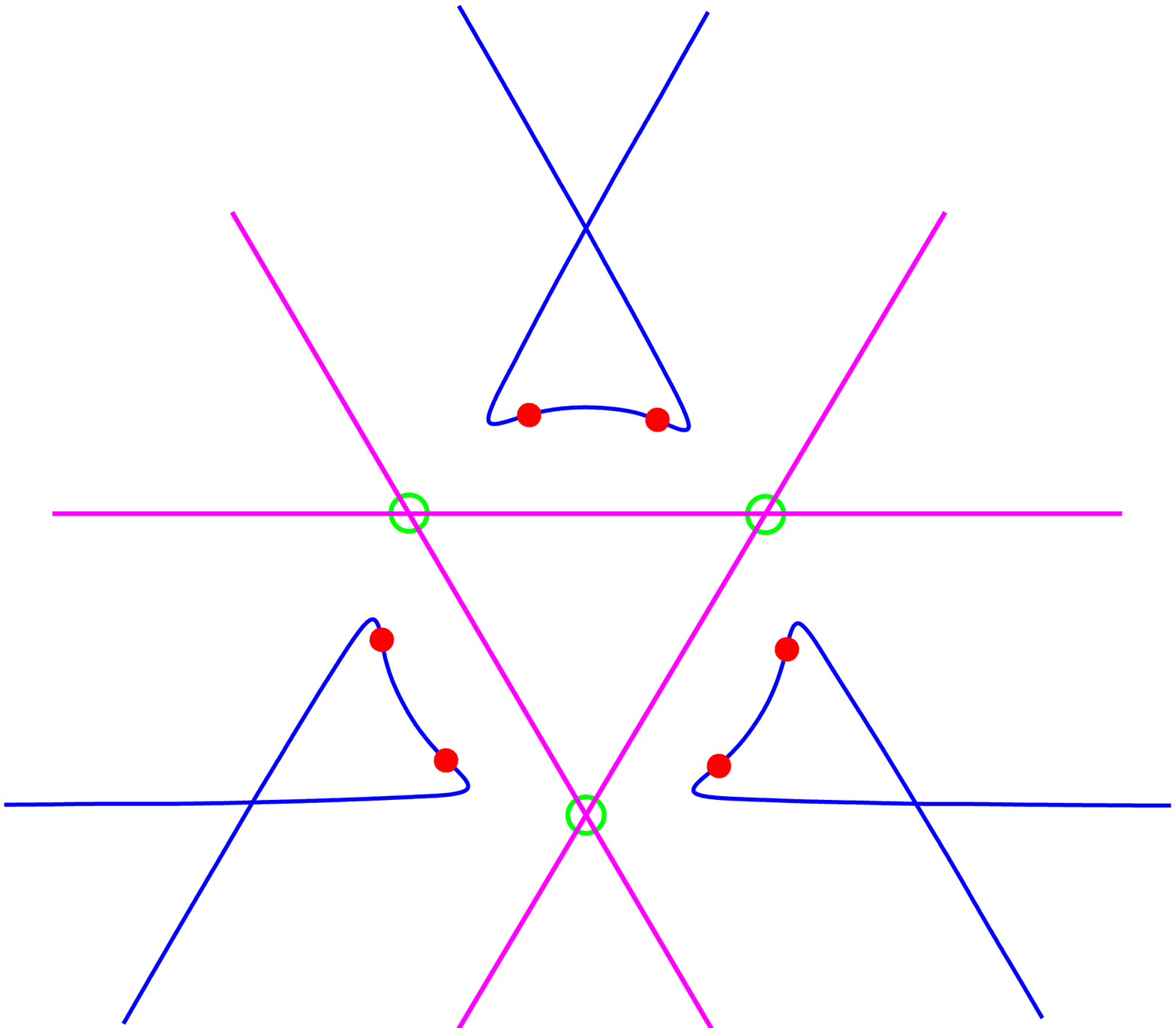}\qquad\qquad
    \epsfysize=90pt\epsfbox{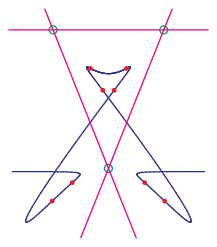}
 \]
 The bisecants form four triangles in $\mathbb{RP}^2$, with a distinguished
 triangle not meeting the pseudoline and containing
 no ovals (i.e., no ovals arising from the loops).
 The pseudoline divides each of the other triangles into two components, with
 the oval in the component adjacent to
 that edge of the distinguished triangle
 which is not intersected by the pseudoline.
 This configuration just described is equivalent to the statement that the
 loops alternate along the pseudoline.\QED

\def\cprime{$'$}

\end{document}